\documentclass[a4paper,reqno]{amsart}
\usepackage{amsmath}
\usepackage{amssymb}
\usepackage{tikz-cd}
\usepackage{color}
\usepackage{verbatim}

\newtheorem{innercustomth}{Theorem}
\newenvironment{customth}[1]
  {\renewcommand\theinnercustomth{#1}\innercustomth}
  {\endinnercustomth}

\usepackage{amscd}
\usepackage{hyperref}

\definecolor{Blue}{rgb}{0.3,0.3,0.9}
\definecolor{Red}{rgb}{0.9,0.3,0.3}

\usepackage[nobysame]{amsrefs}

\newtheorem{Lemma}{Lemma}[section]
\newtheorem{Th}[Lemma]{Theorem}
\newtheorem{PDTh}[Lemma]{Package Deal Theorem}
\newtheorem{Prop}[Lemma]{Proposition}

\newtheorem{Cor}[Lemma]{Corollary}
\newtheorem{Construction}[Lemma]{Construction}

\theoremstyle{definition}

\newtheorem{Def}[Lemma]{Definition}
\newtheorem{Ex}[Lemma]{Example}

\theoremstyle{remark}
\newtheorem{Rem}[Lemma]{Remark}

\newtheorem{Remark}[Lemma]{Remark}
\newenvironment{Proof}{{\sc Proof.}\ }{~\rule{1ex}{1ex}\vspace{0.5truecm}}

\newcommand{\add}{\mbox{\rm add}}
\newcommand{\Add}{\mbox{\rm Add}}

\newcommand{\im}{\mbox{\rm Im}}

\newcommand{\N}{\mathbb N}
\newcommand{\No}{{\mathbb N}_0}

\newcommand{\fm}{\mathfrak{m}}
\newcommand{\fn}{\mathfrak{n}}

\newcommand{\Mod}{\mbox{\rm Mod-}}

\DeclareMathOperator{\mSpec}{mSpec}

\DeclareMathOperator{\Ker}{Ker}
\DeclareMathOperator{\Hom}{Hom}
\DeclareMathOperator{\End}{End}
\DeclareMathOperator{\Tr}{Tr}
\DeclareMathOperator{\rank}{rank}

\title[Torsion-free modules over one dimensional domains] {Torsion-free Modules over Commutative Domains of Krull Dimension One}

\begin{document}

\author{Rom\'an \'Alvarez}

\address{Departament de Matem\`atiques, 
Universitat Aut\`onoma de Barcelona, 08193 Bellaterra
(Barcelona), Spain\newline 
Charles University, Faculty of Mathematics and Physics \\Department
of Algebra, Sokolovsk\'a~83,
18675 Praha 8, Czech Republic}
\email{roman9496@gmail.com}

\thanks{The research of the first author was supported by the pre-doctoral grant 2021FI-B 00913 of the Generalitat de Catalunya}
\author{Dolors Herbera}

\address{Departament de Matem\`atiques,  
Universitat Aut\`onoma de Barcelona, 08193 Bellaterra
(Barcelona), Spain \newline
Centre de Recerca Matemàtica,  08193 Bellaterra
(Barcelona), Spain}
\email{dolors.herbera@uab.cat}

\thanks{The second author was supported by the Spanish State Research Agency, through the Severo Ochoa and María de Maeztu Program for Centers and Units of Excellence in R\& D (CEX2020-001084-M). The first and second authors were partially supported by the projects MIMECO  PID2020-113047GB-I00 financed by the Spanish Government and  \emph{Laboratori d'Interaccions entre Geometria, \`Algebra i Topologia} (LIGAT) with reference number 2021 SGR 01015 financed by the Generalitat de Catalunya. }
 
 \author{Pavel P\v r\'\i hoda}
\address{Charles University, Faculty of Mathematics and Physics \\Department
of Algebra, Sokolovsk\'a~83,
18675 Praha 8, Czech Republic}
\email{prihoda@karlin.mff.cuni.cz}
\thanks{The third author was supported by Czech Science Foundation grant GA\v CR 23-05148S}

\date{\today}

\begin{abstract} 

Let $R$ be a domain of Krull dimension one. We study when  the class $\mathcal{F}$  of modules over $R$ that are arbitrary  direct sums  of finitely generated torsion-free modules is closed under direct summands. 

If  $R$ is local, we show that $\mathcal{F}$ is closed under direct summands if and only if any indecomposable, finitely generated, torsion-free module has local endomorphism ring. If, in addition, $R$ is noetherian this is equivalent to saying that the normalization of $R$ is a local ring.

If $R$ is an $h$-local domain of Krull dimension $1$ and $\mathcal{F}_R$  is closed under direct summands,  then the property is inherited by the localizations of $R$ at maximal ideals. Moreover,  any localization of $R$ at a maximal ideal, except maybe one, satisfies that any finitely generated ideal is $2$-generated. The converse  is true when the domain $R$ is, in addition, integrally closed, or noetherian semilocal, or   noetherian    with module-finite normalization.

Finally, over a commutative domain of finite character and with no restriction on the Krull dimension, we show that the isomorphism classes of countably generated modules in $\mathcal{F}$ are determined by their genus.
\end{abstract}

\maketitle

\setcounter{tocdepth}{1}
\tableofcontents

\section*{Introduction}

Let $R$ be a commutative domain, $\Lambda$ an $R$-algebra, and let $\mathcal{F}_\Lambda$ denote the class of $\Lambda$-modules that are direct sums of finitely generated $\Lambda$-modules which are torsion-free over $R$ (we should write $\mathcal{F}$ if the ring is clear). The class $\mathcal{F}_\Lambda$ is closed under arbitrary direct sums, and we ask whether $\mathcal{F}_\Lambda$ is closed under direct summands. When $\Lambda$ is noetherian, this  is equivalent to the question whether $\mathcal{F}_{\Lambda}$ coincides with the class of pure projective $\Lambda$-modules which are torsion-free over $R$. 

If $R$ is a semihereditary domain, then the modules in $\mathcal{F}_{R}$ are projective modules. It is well known that, in this case, all projective modules are direct sums of finitely generated ideals of $R$. In particular, $\mathcal{F}_{R}$ is the class of all projective modules, and then it is closed under direct summands with no extra assumption on the Krull dimension of $R$.

The interest in these kinds of questions was awakened with the study of the so-called, \emph{generalized lattices}. Let  $R$ be a Dedekind domain with field of fractions $Q$, and let $\Lambda$ be an $R$-order in a separable $Q$-algebra (hence, $\Lambda$ is finitely generated and projective as $R$-module). Classically, a  right $\Lambda$-module is a lattice if it is finitely generated and projective as an $R$-module. Dropping the finitely generated condition from the definition of lattice, we encounter the  generalized $\Lambda$-lattices. In this context, the class $\mathcal{F}_\Lambda$ coincides with the 
(arbitrary) direct sums of $\Lambda$-lattices. If $\Lambda$ is locally lattice-finite, the  generalized lattices are precisely the direct summands of modules in $\mathcal{F}_\Lambda$. See \cite[Theorem~4]{butler} for the proof in the lattice-finite case and  \cite[p. 112, Corollary]{Rump} for its extension to the locally lattice-finite case. 

In \cite[Theorem~4]{Rump}, W.~Rump proved that if generalized $\Lambda$-lattices are direct sums of lattices, then $\Lambda$ is locally lattice-finite.  Moreover, he answered the question of when $\mathcal{F}_\Lambda$ is  closed under direct summands, finding that the answer depended on representation-theoretical data of $\Lambda$. 

P\v r\'\i hoda in the paper~\cite{P3} translated the question to $\Lambda = R$ being a noetherian domain of Krull dimension 1 with module-finite normalization. In the local case, he proved that $\mathcal{F}_R$ is closed under direct summands if and only if its normalization (which is assumed to be finitely generated) is local.
In the global case, and under the assumption that $R$ is, in addition,  lattice-finite (i.e., there are only finitely many indecomposable finitely generated torsion-free $R$-modules up to an isomorphism) he characterized when $\mathcal{F}_R$ is closed  
under direct summands  by certain ring-theoretical properties of $R$ that  will also appear in Theorem~\ref{B}.

In the present paper, we  consider the case  $\Lambda = R$ being a commutative domain of Krull dimension $1$, and we also follow the strategy of studying first the case of a local domain, and then extending to a suitable global situation.  For the case of local domains, we get quite a satisfactory characterization  that is summarized in the following theorem,

\begin{customth}{A} \label{A} \emph{(Corollary~\ref{generalrank}, Corollary~\ref{local})}
Let $R$ be a   local domain of Krull dimension $1$. The class $\mathcal{F}$ is closed under direct summands if and only if any finitely generated, indecomposable, torsion-free module has local endomorphism ring. If $R$, in addition, is noetherian this is equivalent to having local integral closure.

Therefore, in this situation, any module $M$ in $\mathcal{F}$ can be written as $M=\bigoplus _{i\in I}N_i$, where each $N_i$ is a finitely generated module with local endomorphism ring, and such decomposition is unique. In addition, any direct summand of $M$ is isomorphic to a direct sum of the modules $\{N_i\}_{i\in I}$.
\end{customth}

Our global setting will be the one of $h$-local domains. Matlis introduced $h$-local domains in the 60s in \cite{matlis3} as commutative domains satisfying that:
\begin{itemize}
    \item[(i)] they have   finite character. That is, every non-zero ideal is contained only in finitely many maximal ideals.
    \item[(ii)]  every non-zero prime ideal is contained only in one maximal ideal. 
\end{itemize}
However, Jaffard appears to be the first author to give an equivalent definition of $h$-local domains in~\cite{jaffard}, but under the name of \textit{domains of Dedekind type}. For the story and characterizations of $h$-local domains, the reader is referred to Fuchs and Salce's monograph \cite[\S IV.3]{fuchssalce} and to the paper by Olberding \cite[Theorem~2.1]{olberding}. 

A domain of Krull dimension 1 is $h$-local if and only if it has finite character, and  noetherian domains of Krull dimension $1$ are always $h$-local (equivalently, always have finite character).

If $R$ is a domain such that $\mathcal{F}_R$ is closed under direct summands then, in particular, all   projective modules are direct sums of finitely generated projective modules. Not so much is known about projective modules over arbitrary domains. However,  Hinohara   \cite{hinohara} proved that, over an $h$-local domain, an infinitely generated projective module, which is not finitely generated,  is always free. 

The following theorem summarizes our results for this class of domains,

\begin{customth}{B} \label{B}\emph{(Theorem~\ref{th2generated}, Corollary~\ref{integrallyclosed}, Corollary~\ref{integralclosure}, Proposition~\ref{semilocal}, Theorem~\ref{converse})} 
Let $R$ be an $h$-local domain   of Krull dimension $1$. If the class $\mathcal{F}_R$ is closed under direct summands, then $R$ satisfies that
\begin{enumerate}
    \item[(1)] $\mathcal{F}_{R_\fm}$ is closed under direct summands for any maximal ideal $\fm$ of $R$;
    \item[(2)] for any maximal ideal $\fm$ of $R$, except maybe one, all finitely generated ideals of $R_\fm$ are at most two-generated;
    \item[(3)] the integral closure $\overline{R}$ of $R$ in its field of fractions satisfies that $\mathcal{F}_{\overline{R}}$ is closed under direct summands.
\end{enumerate}
 Conversely, if $R$ is an $h$-local domain of Krull dimension $1$ satisfying $(1)$ and $(2)$, then $\mathcal{F}_R$ is closed under direct summands in the following situations:
 \begin{itemize}
     \item[(a)] $R$ is integrally closed, or 
     \item[(b)] $R$ is semilocal noetherian,  or  
     \item[(c)] $R$ is noetherian and has module-finite normalization.
\end{itemize}
\end{customth}

That the closure of $\mathcal{F}$ under direct summands implies conditions (1) and (2) was already proven by P\v r\'\i hoda for noetherian domains of Krull dimension $1$  with module-finite normalization \cite{P3}. The fact that a condition like $(2)$ still holds in the generality we are working  has been an interesting surprise, and its proof is one of the more involved in the paper.  

A good description for torsion-free modules over a local ring with the $2$-generated property for finitely generated ideals is only available in the noetherian case: if $R$ has module-finite integral closure, then it is a Bass domain, and the theory goes back to Bass' fundamental paper \cite{bass}. The general case of noetherian domains with two-generated ideals remained open for a long time and was finally  settled by Rush in \cite{rush}. In both cases, it is proven that finitely generated torsion-free modules are isomorphic to a direct sum of (finitely generated) ideals. As far as we know, it is an open question whether such a result could still hold outside the noetherian case. 

We also study the isomorphism classes of modules in  $\mathcal{F}$. For finitely generated modules, the isomorphism class of a module $M$ is usually determined by its genus and some extra information on the \emph{class group}. This is the case, for example, of Bass domains as it was shown in  
\cite[Theorem 4.2]{levywiegand} and for general noetherian domains with $2$-generated ideals \cite{rush}. For countable direct  sums in $\mathcal{F}$ we show that the genus is enough to determine the isomorphism class. Now we state the precise result we prove, notice that it is for finitely generated $R$-algebras and that there is no assumption on the Krull dimension of the domain $R$.

\begin{customth}{C} \emph{(Theorem~\ref{countablesum})}\label{C}
Let $R$ be a commutative domain of finite character with field of fractions $Q$. Let $\Lambda$ be a module-finite $R$-algebra such that $\Lambda \otimes _R Q$ is a simple artinian ring. Let $M=\bigoplus_{i\in\No} A_i$ and $N=\bigoplus_{i\in\No} B_i$ be direct sums of non-zero finitely generated right $\Lambda$-modules which are torsion-free as $R$-modules. If $M$ and $N$ are in the same genus, then there are decompositions
    \[\textstyle M=\bigoplus_{i\in\No} M_i\qquad\text{and}\qquad N=\bigoplus_{i\in\No} N_i\]
such that both $M_i$ and $N_i$ are finitely generated, and $M_i\cong N_i$ for every $i\in\No$. In particular, $M$ and $N$ are isomorphic.
\end{customth}

Theorem~\ref{C} extends a result proven by Rump in \cite{Rump} for $\Lambda$  an  order in a separable algebra over a Dedekind domain $R$.

Now we briefly describe the content of the different sections of the paper, taking the chance to highlight some of the main ideas.

To work in a problem like the closure of  $\mathcal{F}$ under direct summands, we need some method to construct \emph{interesting and possibly infinitely generated} direct summands. This is the main topic of Section~\ref{s:lifting}, in which we  exploit  that projective modules can be lifted modulo the  trace ideal of a projective module (cf. Theorem~\ref{liftingproj}). This machinery for projective modules was   developed in the noetherian setting in \cite{P2},  extended to general rings in \cite{traces} and the connection with the study of the category $\Add(M)$, of direct summands of direct sums of copies of a finitely generated module $M$,  was developed in  \cite{wiegand} (cf. Theorem~\ref{equivalencia}). 

Section~\ref{s:endo} is quite technical but essential for the rest of the paper. In it,   we survey the results needed on finitely generated (torsion-free) and/or finite rank modules over a domain $R$,  as well as properties of their endomorphism rings. Given a domain $R$ we also need to treat the case of  finitely generated modules over an $R$-algebra $\Lambda$ that is torsion-free as $R$-modules. We also introduce the concepts of domain of finite character and of $h$-local domain, developing the first basic properties of endomorphism rings of finitely generated modules over these classes of rings.

In Section~\ref{s:local} we characterize local domains of Krull dimension $1$ satisfying  that $\mathcal{F}$ is closed under direct summands. The results in this section prove Theorem~\ref{A}. 

In Section~\ref{s:dsummands}  we develop tools to relate the property of being locally a direct summand  with being a direct summand  in the setting of finitely generated torsion-free modules. If $R$ is a domain of finite character, and $N$ and $M$ are finitely generated torsion-free modules such that $N\otimes_R Q$ is a direct summand of $M\otimes_R Q$, then $N_\fm$ is a direct summand of $M_\fm$ for almost all maximal ideals $\fm$ of $R$. This means that we only have a finite number of ``problematic maximal ideals.'' In the remainder of the section, we prove several results on how to proceed for these problematic maximal ideals. Analogous results were proven in the one-dimensional noetherian setting in~\cite{GL} and were inspired by the representation theory of orders. 

In Section~\ref{s:deal} we study finitely generated modules over an algebra $\Lambda$ over an $h$-local domain $R$. We prove the results that will serve as a bridge between the local and the global case. 
A key  observation is the following: if $M$ and $N$ are finitely generated $\Lambda$-modules that are torsion-free as $R$-modules and satisfy $M\otimes _RQ\cong N\otimes _R Q$, where $Q$ is the field of fractions of $R$, then the localizations of $M$ and $N$ at maximal ideals of $R$ are isomorphic for almost all maximal ideals of $R$. 

Then the problem is, given $\{X(\fm)\mid \fm\in\mSpec R\}$  a family of right $\Lambda_\fm$-modules, when there is a right $\Lambda$-module $N$ such that $N_\fm\cong X(\fm)$ for every maximal ideal $\fm$ of $R$. 
This question was already studied by Levy and Odenthal in \cite{levyodenthal2} for $R$ being a noetherian domain of Krull dimension 1, and the answer is that mild compatible conditions between the localizations are enough to warrant the existence of $N$.  Here we give a generalized version of Levy-Odenthal's Package Deal Theorem for localizations over $h$-local domains (cf. Theorem~\ref{dealsubmodules}), and give another one for localizations of trace ideals of countably generated projective $\Lambda$-modules (cf. Theorem~\ref{dealtraces}). 

The case $\Lambda =R$ is particularly neat. The localization at a maximal ideal of a finitely generated torsion-free module  is  free for almost all maximal ideals (cf. Corollary~\ref{almost-free}).  Conversely, a family $\{X(\fm)\mid \fm\in\mSpec R\}$  of  $R_\fm$-modules can be \emph{glued together} in a finitely generated torsion-free $R$-module $N$ provided all the modules $X_{\fm}$ have the same rank and  they are free for almost all maximal ideals $\fm$ (cf. Corollary~\ref{existence}).

In Section~\ref{s:global} we exploit the theory developed in the previous sections, proving the first part of Theorem~\ref{B}. An important intermediate result is Theorem~\ref{coprime}, which shows that if $\mathcal{F}_R$ is closed under direct summands, then finitely generated indecomposable torsion-free modules over different localizations at maximal ideals must have coprime rank. This type of result was first spotted by P\v r\'\i hoda in \cite{P3}, and our proof follows the same strategy as the original result. 

In Section~\ref{s:ic}, the remaining part of the characterization of integrally closed domains with $\mathcal{F}$ closed under direct summands, contained in Theorem~\ref{B}, is proven. 

In Section~\ref{s:genus} we give the proof of Theorem~\ref{C}. This section also includes examples showing that over a semilocal domain, being a direct summand of an infinite direct sum of finitely generated torsion-free modules could be satisfied locally but not globally. This contrasts with the case of finitely generated modules, cf.~ Corollary~\ref{semilocalsummands}.

The noetherian part of Theorem~\ref{B} is proven in Section~\ref{s:noeth}. We stress that Theorem~\ref{converse} gives the converse of one of the main results in \cite{P3}. 

We close the paper in Section~\ref{s:example},  with a family of examples    of $h$-local domains of Krull dimension $1$ satisfying Theorem~\ref{B}.  These examples are of the form $R=K+xL[x]$ where  $K\subseteq L$ is a field extension. If such an extension is finitely generated, then $R$ is noetherian. If the extension is transcendental, then $R$ is integrally closed in its field of fractions (and it is not a Pr\"ufer domain!). Moreover, under mild restrictions over the dimension of the field extension, $R$ has indecomposable finitely generated torsion-free modules of all ranks.

\bigskip

Throughout the paper, rings are associative with $1$, and morphisms are unital. In each section we try to be very precise about the setting we are working with. Our main topic is (torsion-free) modules over commutative domains, so most of the time $R$ is a commutative domain with field of fractions $Q$. We reserve $\Lambda$ to denote an algebra over a domain $R$.  If $M$ is a $\Lambda$-module (usually finitely generated and torsion-free as an $R$-module), we denote by  $S$ its endomorphism ring.

\bigskip

We are very grateful to the anonymous referees for the careful reading of the paper, and for a number of interesting observations that helped to improve its readability.

\section{Lifting direct summands modulo the stable category}\label{s:lifting} 

\begin{Def}
Let $R$ be a ring and let $M$ be a right $R$-module. Then the \emph{trace} of $M$ in $R$ is the two-sided ideal of $R$ defined as $$\Tr _R(M)=\sum _{f\in \Hom_R(M,R)} f(M).$$
\end{Def}

To simplify the notation, we will just use $\Tr (M)$ to denote the trace in $R$ of the module $M$ when the ring $R$ is clear. Our main interest will be in traces of projective modules. Recall that the trace of a projective module is always an idempotent ideal.

The following fact will be used throughout the paper.

\begin{Remark} \label{ideal-end}
Let $R$ be a ring and $M$ a right $R$-module. Let $L$ be a two-sided ideal of $R$. Then $I=\{f\in\End_R(M)\mid f(M)\subseteq ML\}$ is a two-sided ideal of $\End_R(M)$.

In addition, the $R$-$R$ bimodule structure of $R$ gives a structure of left $R$-module to $\mathrm{Hom}_R(M,R)$ where $r\cdot \omega$, for $r\in R$ and $\omega \in \mathrm{Hom}_R(M,R)$, is defined as the module homomorphism $r\cdot \omega \colon M\to R$ such that $r\cdot \omega (x)=r\omega (x)$ for any $x\in M$. This gives a canonical morphism
\[\varphi \colon M\otimes _R \mathrm{Hom}_R(M,R)\to \mathrm{End}_R(M)\]
where, for $m\in M$ and $\omega \in \mathrm{Hom}_R(M,R)$, $\varphi (m\otimes \omega)$ is the endomorphism of $M$ defined by $\varphi (m\otimes \omega)(x)=m\omega (x)$, for any $x\in M$. It is not difficult to show that $\varphi (M\otimes _R \Hom_R (M,R))= \{f\in S\mid \mbox{$f$ factors through $R^n$ for some $n\ge 1$} \}$. 

Throughout the paper we will use the notation $$\varphi (M\otimes _R \Hom_R (M,R))=M \Hom_R (M,R).$$
\end{Remark}

\begin{Lemma} \label{stableideal}
Let $R$ be a ring. Let $M$ be a right $R$-module with endomorphism ring $S=\End_R(M)$. Let $I=\Tr_R(M)$, and $J=\{f\in S\mid f(M)\subseteq MI\}$ (which is a two-sided ideal of $S$ by Remark~\ref{ideal-end}). Then 
$M \Hom_R (M,R)\subseteq J.$
\end{Lemma}

\begin{Proof} In the notation of Remark~\ref{ideal-end} and because of the definition of $\varphi$, it follows that $\varphi (m\otimes \omega )(M)\subseteq MI$. So that, also $\varphi (\sum _{i=1}^n m_i\otimes \omega _i )(M)\subseteq MI$ for $m_i\in M$ and $\omega _i\in \mathrm{Hom}_R(M,R)$.   
\end{Proof}

\begin{Lemma}\label{tracequotients}
Let $R$ be a ring and let $P$ be a right $R$-module. Then, for any two-sided ideal $I$ of $R$,
\begin{itemize}
	\item[(i)] $(\Tr_{R} \left(P\right)+I)/I\subseteq \Tr_{R/I} \left(P/PI\right)$;
	\item[(ii)]  If, in addition, $P$ is projective, $\Tr_{R/I} \left(P/PI\right)=(\Tr_{R} \left(P\right)+I)/I$.
\end{itemize}	
\end{Lemma}
\begin{Proof}
$(i)$. Let $x\in\Tr_R(P)$, then there exist $f_1,\dotsc,f_n\in\Hom_R(P,R)$ and $p_1,\dotsc,p_n\in P$ such that $x=\sum_{i=1}^n f_i(p_i)$. Then, for each $i=1,\dotsc,n$, there is an induced homomorphism $\overline f_i\colon P/PI\to R/I$ given by $\overline f_i(p+PI)=f(p)+I$. Since $\overline f_i(P/PI)\subseteq \Tr_{R/I}(P/PI)$, in particular, $f_i(p_i)+I\in\Tr_{R/I}(P/PI)$. Since $x=\sum_{i=1}^n f_i(p_i)$, it follows that $x+I\in\Tr_{R/I}(P/PI)$.

Statement $(ii)$ follows from the lifting property of projective modules, since every $R/I$-module homomorphism $\overline f\colon P/PI\to R/I$ lifts to an $R$-module homomorphism $f\colon P\to R$.
\end{Proof}

\begin{Remark}
\label{Eilenberg} In this section we will make repeated use of the so-called \emph{Eilenberg trick}, which implies, for example, that if $P$ is a countably generated projective module and $X$ is a direct summand of $P$, then $P^{(\omega)}\oplus X\cong P^{(\omega)}$.
\end{Remark}

We recall the definition of the class of modules generated by a projective module.

\begin{Def}
Let $R$ be a ring. Let $M$ be a right $R$-module. We define the \textit{class $\mathrm{Gen}(M)$ of modules generated by $M$} as the class of right $R$-modules that are a homomorphic image of a direct sum of copies of $M$.
\end{Def}

For further quoting, we recall the following characterization of the class $\mathrm{Gen}(P)$, when $P$ is a projective module.

\begin{Lemma} \label{genP} {\rm \cite[Lemma~2.10]{wiegand}}
Let $R$ be a ring, and let $P$ be a projective right $R$-module with trace ideal $I$. Then $\mathrm{Gen}\, (P)$ 
coincides with the class of right $R$-modules $M$ such that $MI=M$. 
\end{Lemma}

For further use, we note the following property of direct sums of finitely generated modules.

\begin{Lemma} \label{finitesummand}
Let $R$ be a ring, and assume there is  a direct sum decomposition of right $R$-modules $M'\oplus M=X\oplus Y$ with $M$ finitely generated and $X$, $Y$ direct sums of finitely generated modules. Then there exist $A$ and $B$  direct summands of $X$ and $Y$, respectively, such that $M\oplus Z=A\oplus B$ with $Z$ a  finitely generated direct summand of $M'$.
\end{Lemma}

\begin{Proof}
Since $M$ is finitely generated, by the hypothesis, $M\subseteq A\oplus B$ with $A$ and $B$ being finitely generated direct summands of $X$ and $Y$, respectively. Therefore, by the modular law, $A\oplus B=M\oplus Z$, where $Z=M'\cap (A\oplus B)$.
\end{Proof}

\begin{Th} \label{liftingproj}
Let $R$ be a ring, and let $I$ be an ideal of $R$ that is the trace of a countably generated projective right $R$-module $Q$. 
\begin{itemize}
    \item[(i)] \emph{\cite[Theorem 3.1]{traces}} Let $P'$ be a countably generated projective right module over $R/I$. Then there exists a countably generated projective right $R$-module $P$ such that  $P/PI \cong P'$.
    \item[(ii)] Let $P_1$ and $P_2$ be countably generated projective modules such that $P_1/P_1I\cong P_2/P_2I$ then  $Q^{(\omega)}\oplus P_1\cong Q^{(\omega)}\oplus P_2$.
    \item[(iii)] Let $Q'$ be a countably generated projective $R$-module such that  $\mathrm{Tr}_R(Q')=I$. Then $Q^{(\omega)}\cong (Q')^{(\omega)}$.
\end{itemize}
\end{Th}
\begin{Proof} 
Statement $(ii)$ follows from the proof of \cite[Lemma~2.5]{P2}. We give a self-contained proof for the reader's convenience.

For each $i=1,2$, let $\pi_i\colon P_i \to P_i/P_iI$ denote the canonical projection. Let $f\colon P_1/P_1I\to P_2/P_2I$ be an isomorphism. Then there exists $g_1\colon P_1 \to P_2$ such that the diagram
    \[\begin{tikzcd}[column sep=huge,row sep=huge]
        &P_1\arrow[dl,"g_1"']\arrow[d, "f \circ \pi _1"]\\
        P_2\arrow[r, "\pi _2"]& P_2/P_2I
    \end{tikzcd}\]
is commutative. Hence, $P_2=g_1(P_1)+P_2I$.

Since $I$ is the trace ideal of the countably generated projective module $Q$, it is idempotent and, by Lemma~\ref{genP}, there exists an onto module homomorphism $g_2\colon Q^{(\kappa)}\to P_2I$, where $\kappa$ is an infinite cardinal. Let $h\colon Q^{(\kappa)}\oplus P_1 \to P_2/P_2I$ be the module homomorphism defined by $h(q,p)=f \circ \pi _1(p)$ for any $(q,p)\in Q^{(\kappa)}\oplus P_1$. Then there is a commutative diagram
    \[\begin{tikzcd}[column sep=huge,row sep=huge]
        &Q^{(\kappa)}\oplus P_1\arrow[dl,"g"']\arrow[d, "h"]\\
        P_2\arrow[r, "\pi _2"]& P_2/P_2I
    \end{tikzcd}\]
where $g$ is the onto module homomorphism defined by $g(q,p)= g_2(q)+g_1(p)$ for any $(q,p)\in Q^{(\kappa)}\oplus P_1$. Since $P_2$ is countably generated, there exists $C \subseteq \kappa$ countably infinite such that $P_2 = g(Q^{(C)} \oplus P_1)$. Let $g'\colon Q^{(C)}\oplus P_1 \to P_2$ be the restriction of $g$ to $Q^{(C)} \oplus P_1$. Since $P_2$ is projective, $Q^{(C)}\oplus P_1= P'_2\oplus \Ker g'$ where $P_2'$ is a submodule of  $Q^{(C)}\oplus P_1$ isomorphic to $P_2$. 

As $\Ker g' \subseteq \Ker h=Q^{(\kappa)}\oplus P_1I$, we deduce that 
    $$Q^{(C)}\oplus P_1I=\Ker g'\oplus (P_2'\cap (Q^{(C)}\oplus P_1I)).$$
 By Lemma~\ref{genP}, $P_1I\in \mathrm{Gen}\, (Q)$. Hence,  the countably generated module $\Ker g'$ is a homomorphic image of $Q^{(\omega)}$. Since it is also projective, it is a direct summand of $Q^{(\omega)}$. By the Eilenberg trick (cf. Remark~\ref{Eilenberg}), $Q^{(\omega)}\oplus \Ker g' \cong Q^{(\omega)}$. Therefore,
    $$Q^{(\omega)}\oplus P_1\cong Q^{(\omega)}\oplus P_2\oplus \Ker g'\cong Q^{(\omega)}\oplus P_2, $$
as we wanted to prove.

Statement $(iii)$ is a consequence of $(ii)$. Taking $P_1= Q^{(\omega)}$ and $P_2= (Q')^{(\omega)}$, we deduce from $(ii)$,
\[Q^{(\omega)}\cong Q^{(\omega)} \oplus Q^{(\omega)}\cong Q^{(\omega)}\oplus (Q')^{(\omega)}.\]
Exchanging the roles of $Q$ and $Q'$, we deduce from $(ii)$ that $(Q')^{(\omega)}\cong (Q')^{(\omega)}\oplus  Q^{(\omega)} $, which yields the claimed isomorphism.
\end{Proof}

As a corollary of Theorem~\ref{liftingproj} we show that direct sums of projective modules can be lifted modulo a trace ideal. This result is implicit in \cite{wiegand}. 

\begin{Cor} \label{liftingds}  
Let $R$ be a ring, and let $P$ be a countably generated projective right $R$-module. Let $I$ be the trace of a countably generated projective right $R$-module $Q$. Assume that $P/PI= X'\oplus Y'$. Then there exist countably generated projective right $R$-modules $X$ and $Y$ such that  $Q^{(\omega)}\oplus P\cong X\oplus Y$ and  $X/XI\cong X'$, $Y/YI\cong Y'$ as $R/I$-modules.

If, in addition, $P$ is finitely generated and $X$, $Y$ are direct sums of finitely generated modules, then there exist $A$ and $B$ isomorphic to direct summands of $X$ and $Y$, respectively, such that $P\oplus Z=A\oplus B$, with $Z$ a  finitely generated direct summand of $Q^{(\omega)}$, and satisfying that $A/AI\cong X'$, $B/BI\cong Y'$. 
\end{Cor}
\begin{Proof}
By Theorem~\ref{liftingproj}, there exist countably generated projective right $R$-modules $X_1$ and $Y_1$, such that $X_1/X_1I\cong X'$, $Y_1/Y_1I\cong Y'$. Hence, $P/PI\cong X_1/X_1I\oplus Y_1/Y_1I$.

Set $X=Q^{(\omega)}\oplus X_1$ and $Y=Q^{(\omega)}\oplus Y_1$. By Theorem~\ref{liftingproj}, $Q^{(\omega)}\oplus P\cong X\oplus Y$.  

To prove the final part of the statement, assume for simplicity that $Q^{(\omega)} \oplus P = X \oplus Y$. Notice that the existence of $A$, $B$ and $Z$ follows from Lemma~\ref{finitesummand}. Therefore, there are direct sum decompositions $X=A\oplus A'$ and $Y=B\oplus B'$. Now
    $$Q^{(\omega)}\oplus P = A\oplus A'\oplus B \oplus B' = P \oplus (Z\oplus A'\oplus B').$$
so that the canonical projection $\pi \colon Q^{(\omega)}\oplus P\to Z\oplus A'\oplus B'$ induces an isomorphism $Q^{(\omega)}\to Z\oplus A'\oplus B'$. Therefore $Z$, $A'$ and $B'$ are modules in $\mathrm{Gen}(Q)$. By Lemma~\ref{genP}, we deduce that $X'\cong X/XI\cong A/AI$ and $Y' \cong Y/YI\cong B/BI$. This concludes the proof of the statement.
\end{Proof}

\begin{Prop} \label{equivalencia} \emph{\cite[Theorem 4.7]{libro}}
Let $R$ be a ring. Let $M$ be a right $R$-module. Let $S=\End_R(M)$. Then the functor $\Hom_R(M,-)$ induces a category equivalence between $\add (M)$ and the category of finitely generated projective right $S$-modules. The inverse of this equivalence is given by the functor $-\otimes _S M$.

Assume, in addition, that $M_R$ is finitely generated. Then the functor $\Hom_R(M,-)$ induces a category equivalence between $\Add (M)$ and the category of projective right $S$-modules. 
\end{Prop}

\begin{Lemma} {\rm \cite[Theorem~2.11(i)]{wiegand}} \label{chartraceend} 
Let $R$ be a ring. Let $M$ be a finitely generated right $R$-module with endomorphism ring $S$. Let $X$ be an object of $\Add (M)$,  and let  $P_S=\Hom_R(M,X)$. Set $I=\Tr _S(P)$. Then,
\begin{itemize}
    \item[(i)] $I =\{f \in S\mid \mbox{$f$ factors through $X^n$ for some $n\in \mathbb{N}$}\},$ that is 
        $$I=\Hom_R(X,M)\Hom_R(M,X).$$
    \item[(ii)] For any $Y\in \Add (M)$, 
        $$\Hom_R(M,Y)I=\{f \in \Hom_R(M,Y)\mid \mbox{$f$ factors through $X^n$ for some $n\in \mathbb{N}$}\}.$$
\end{itemize}
\end{Lemma}

\begin{Remark} \label{mplusM}
Let $M$ and $M'$ be two right modules over a ring $R$, and with endomorphism rings $S$ and $S'$, respectively.
Let	
	$$T:=\End_R(M'\oplus M)=\begin{pmatrix} S'&\Hom_R (M,M')\\ \Hom_R (M',M)&S \end{pmatrix}$$

Consider the finitely generated projective right $T$-module $P=\left(\begin{smallmatrix}1&0\\ 0&0\end{smallmatrix}\right)T\cong \Hom_R(M'\oplus M,M')$.  Its trace ideal is 
    $$I=T\begin{pmatrix}1&0\\ 0&0\end{pmatrix}T=\begin{pmatrix}S'&\Hom_R (M,M')\\ \Hom_R (M',M)&\Hom_R (M',M)\Hom_R (M,M')\end{pmatrix}.$$ 
Then $T/I\cong S/\Hom_R (M',M) \Hom_R (M,M')$.

We will use the case $M'=R$, so that $T/I\cong S/M  \Hom_R (M,R)$. That is, it is the endomorphism ring of $M_R$ in the stable category.
\end{Remark}

\begin{Cor} \label{stable}
Let $R$ be a ring. Let $M$ be a finitely generated right $R$-module with endomorphism ring $S$. Let 
    \[J=\{f\in S\mid \mbox{$f$ factors through $R^n$ for some $n\in \mathbb{N}$}\}=M\Hom _R(M,R).\]
Let $e$ be an idempotent of $S/J$. Then there is a direct sum decomposition $R^{(\omega)}\oplus M=X\oplus Y$,  such that 
    \[\Hom_R(M,X)/X\Hom_R(M,R)\cong \Hom_R(R\oplus M,X)/X\Hom_R(R\oplus M,R)\cong e(S/J)\]
and 
    \[\Hom_R(M,Y)/Y\Hom_R(M,R)\cong \Hom_R(R\oplus M,Y)/Y\Hom_R(R\oplus M,R) \cong (1-e)(S/J).\]

If, in addition, $X$ and $Y$ are direct sums of finitely generated modules, then there exist $A$ and $B$ direct summands of $X$ and $Y$, respectively, such that $M\oplus P=A\oplus B$ with $P$ a finitely generated projective right $R$-module and satisfying that
    \[\Hom_R(M,A)/A\Hom_R(M,R) \cong \Hom_R(R\oplus M,A)/A\Hom_R(R\oplus M,R)\cong e(S/J)\]
and
    \[\Hom_R(M,B)/B\Hom_R(M,R)\cong \Hom_R(R\oplus M,B)/B\Hom_R(R\oplus M,R)\cong (1-e)(S/J).\]
\end{Cor}

\begin{Proof}
Let $T= \End_R (R\oplus M)$. Taking into account Remark~\ref{mplusM}, $S/J \cong T/I$ with $I$ the trace ideal of the projective right $T$-module $Q=\Hom_R (R\oplus M, R)$. That is, following the notation in Remark~\ref{mplusM}, $I=T\begin{pmatrix}1&0\\ 0&0\end{pmatrix}T$.

By Corollary~\ref{liftingds} and its proof, there is a direct sum decomposition $Q^{(\omega)}\oplus \Hom (R\oplus M, M)= X'\oplus Y'$ with $X'/X'I \cong e (S/J)$ and $Y'/Y'I \cong (1-e) (S/J)$ and such that $X'\oplus Q^{(\omega)}\cong X'$ and $Y'\oplus Q^{(\omega)}\cong Y'$. Notice that, by Theorem~\ref{liftingproj}, these properties determine $X'$ and $Y'$ up to endomorphism.

Set $X=X'\otimes _T (R\oplus M)$ and $Y=Y'\otimes _T (R\oplus M)$. By Proposition~\ref{equivalencia}, $R^{(\omega)}\oplus M\cong X\oplus Y$,  $R^{(\omega)}\oplus X\cong X$,  $R^{(\omega)}\oplus Y\cong Y$.
By Lemma~\ref{chartraceend} and using the notation explained in Remark~\ref{ideal-end}, $X$ and $Y$ also satisfy \[\Hom_R(R\oplus M,X)/X\Hom_R(R\oplus M,R)\cong e(S/J)\] and
\[\Hom_R(R\oplus M,Y)/Y\Hom_R(R\oplus M,R) \cong (1-e)(S/J).\]
These properties determine $X$ and $Y$ up to isomorphism.

Let $\pi\colon R\oplus M\to M$ denote the canonical projection. The assignment $g\mapsto g\circ \pi$, for any $g\in \mathrm{Hom}_R(M,X) $, defines a morphism of abelian groups $\mathrm{Hom}_R(M,X) \to \mathrm{Hom}_R(R\oplus M,X)$. Composed with the projection $\mathrm{Hom}_R(R\oplus M,X)\to \Hom_R(R\oplus M,X)/X\Hom_R(R\oplus M,R)$ gives an onto morphism of right $S$-modules $$\Phi \colon \mathrm{Hom}_R(M,X)\to \Hom_R(R\oplus M,X)/X\Hom_R(R\oplus M,R).$$ Recall that, since $X\in \mathrm{Add}\, (R\oplus M)$,  $$\Hom_R(R\oplus M,X)/X\Hom_R(R\oplus M,R)=\Hom_R(R\oplus M,X)/\Hom_R(R\oplus M,X)I.$$ 

Clearly,
\[\mathrm{Ker}\,\Phi =\{g\in \mathrm{Hom}_R(M,X) \mid g \mbox{ factors through $R^n$ for some $n\in \N$}\}.\]
That is, $\mathrm{Ker}\,\Phi=X\Hom_R(M,R)$. This finishes the proof of the isomorphism $\Hom_R(M,X)/X\Hom_R(M,R)\cong \Hom_R(R\oplus M,X)/X\Hom_R(R\oplus M,R)$ as right $S$-modules.

In a similar way,  $$\Hom_R(M,Y)/Y\Hom_R(M,R)\cong \Hom_R(R\oplus M,Y)/Y\Hom_R(R\oplus M,R)$$ as right $S$-modules.

If $X$ and $Y$ are sums of finitely generated modules, then Lemma~\ref{finitesummand} allows us to proceed as above to prove the final claim of the statement.
\end{Proof}

\begin{Lemma} \label{exchange} 
Let $R$ be a ring. Let $M$ be a right $R$-module with semiperfect endomorphism ring. If $X$ is a right $R$-module such that $M\oplus X=Y_1\oplus Y_2$ for suitable right $R$-modules $Y_1$, $Y_2$, then for $i=1,2$, $Y_i\cong N _i\oplus Y_i'$ where $N_1\oplus N_2\cong M$ and $Y'_1\oplus Y'_2\cong X$.
\end{Lemma}
\begin{Proof}
This is an easy consequence of the fact that modules with a semiperfect endomorphism ring are a finite direct sum of modules with a local endomorphism ring, and that modules with local endomorphism ring satisfy the exchange property (cf. \cite[Lemma 2.7 and Theorem 2.8]{libro}). Namely, if $M\oplus X=Y_1\oplus Y_2$, and $M$ has local endomorphism ring, then  we may assume that $M$ is isomorphic to a direct summand of  $Y_1$ so that $Y_1\cong M\oplus Y'_1$ and, moreover, $X\cong Y'_1 \oplus Y_2$. Now if $M=M_1\oplus \cdots \oplus M_n$ with each $M_i$ having local endomorphism ring, the statement follows by induction on $n$.
\end{Proof}

\begin{Lemma} \label{tracelocal}
Let $R$ be a semiperfect ring. Let $M$ be a right $R$-module. Then $\Tr_R(M)\subseteq J(R)$ if and only if $M$ has no non-zero projective direct summands.
\end{Lemma}

\begin{Proof}
Assume that $\Tr_R(M)\not \subseteq J(R)$. Hence, there exists $f\colon M\to R$ and $e^2=e\in R$, with $eR/eJ(R)$ simple, such that $ef(M)\neq eJ(R)$. Therefore $eR= ef(M)+eJ(R)$. By Nakayama's Lemma, $eR=ef(M)$ so $eR$ is a projective direct summand of $M$. 

The converse statement is clear because the trace of a non-zero projective module is never contained in $J(R)$. 
\end{Proof}

\begin{Lemma} \label{projdirectsummand}
Let $R$ be a semiperfect ring. Then:
\begin{itemize}
    \item[(i)] If $M$ is a finitely generated module, $M=P\oplus M'$ with $P$ a finitely generated projective module and $M'$ with no non-zero projective direct summands. In addition, this decomposition is unique up to isomorphism.
    \item[(ii)] If $M_1$ and $M_2$ are $R$-modules with no  non-zero projective direct summand then $M_1\oplus M_2$ has the same property.
\end{itemize}
\end{Lemma}
\begin{Proof} 
Statement $(i)$ is \cite[Theorem~3.15]{libro}.

To prove statement $(ii)$, notice that if $P$ is an indecomposable projective $R$-module, then its endomorphism ring is local because $R$ is semiperfect. Therefore, if such $P$ is a direct summand of $M_1\oplus M_2$, then it is a direct summand either of $M_1$ or of $M_2$, which is not possible by our assumptions.
\end{Proof}

\begin{Prop} \label{criteria}
Let $R$ be a semiperfect ring. Let $M$ be a finitely generated  right $R$-module with endomorphism ring $S$. Set $J=M \Hom_R (M,R)$. If every direct summand of $R^{(\omega)}\oplus M$ is a direct sum of finitely generated modules, then the following statements hold:
\begin{itemize}
    \item[(i)] Let $e\in S/J$ be an idempotent, then there exists $f^2=f\in S$ such that  $(f+J)S/J\cong eS/J$ and $(1-f +J)S/J\cong (1-e)S/J$;
    \item[(ii)] if $M$ is indecomposable, then so is the right $S$-module $S/J$;
    \item[(iii)] every direct summand of $R^{(\omega)}\oplus M$ is of the form $P\oplus X$, where $X$ is a direct summand of $M$ and $P$ is a direct summand of $R^{(\omega)}$ (that is, a countably generated projective module).
\end{itemize}
\end{Prop}
\begin{Proof} 
$(i)$. We consider the setting of   Remark~\ref{mplusM} with $M'=R$, and follow the notation on that remark. 

By Corollary~\ref{stable}, there is a decomposition $P\oplus M= A\oplus B$ with $P$ finitely generated projective and $A$ and $B$ finitely generated such that \[\Hom_R(M,A)/A\Hom_R(M,R)\cong e(S/J)\]
and
    \[\Hom_R(M,B)/B\Hom_R(M,R)\cong (1-e)(S/J).\]

By Lemma~\ref{exchange}, $A=A_1\oplus A_2$ and $B=B_1\oplus B_2$ with $M\cong A_1\oplus B_1$ and $A_2$ and $B_2$ projective $R$-modules. Therefore, $\Hom_R(M,A_1)\cong fS$ and $\Hom_R(M,B_1)\cong (1-f)S$ for  $f^2=f\in S$.

By Lemma~\ref{chartraceend}, 
    \begin{align*}
        e(S/J)&\cong\Hom_R(M,A)/A\Hom_R(M,R) \cong \Hom_R(M,A_1)/A_1\Hom_R(M,R) \\
        &= \Hom_R(M,A_1)/\Hom_R(M,A_1)J\cong (f+J)S/J.
    \end{align*}
Similarly, $(1-e) S/J\cong (1-f+J) S/J$.

Statement $(ii)$ is a consequence of $(i)$.

$(iii)$. Since $R$ is semiperfect, $M=P\oplus M'$ where $P$ is finitely generated projective and $M'$ has no non-zero projective direct summands (cf. Lemma~\ref{projdirectsummand}). Since $R^{(\omega)}\oplus P\cong R^{(\omega)}$, to prove our statement, we may assume that $M=M'$ has no non-zero projective direct summands. 

Let $N$ be a finitely generated direct summand of $R^{(\omega)}\oplus M$. Then, there exists $n\ge 0$ such that $N$ is a direct summand of $R^n\oplus M$. By Lemma~\ref{exchange}, we can deduce that $N$ is of the required form.

Assume now that $N$ is a direct summand of $R^{(\omega)}\oplus M$ that is not finitely generated. By hypothesis, $N=\bigoplus _{i\in \mathcal{A}} N_i$ where $N_i$ are finitely generated modules. By the previous step, for each $i\in \mathcal{A}$, $N_i= P_i\oplus M_i$ with $P_i$ finitely generated projective and $M_i$ isomorphic to a direct summand of $M$. Notice that each $M_i$ has no non-zero projective direct summands, by our assumptions.

Since the ring $R$ is semiperfect, then $M/MJ(R)$ is a semisimple module of length    $\ell$. Then the maximal number of non-zero direct summands in a direct sum decomposition of $M$  is $\ell$.  We claim that this implies that only  $\ell$ of the direct summands $M_i$ are different from zero. 

Indeed, take $i_1, \dots , i_{\ell +1}$ to be  $\ell+1$ different elements of the index set $\mathcal{A}$. Then $N'=\bigoplus _{j=1}^{\ell+1} M_{i_j}$ is a direct summand of $R^{(\omega)}\oplus M$. By Lemma~\ref{projdirectsummand}, $N'$ has no non-zero projective direct summands. By the first step, $N'$ is a direct summand of $M$. Hence, there is some $M_{i_j}=\{0\}$. 

Since only finitely many $M_i$'s are different from zero, as before, we can deduce that its   direct sum is a direct summand of $M$. Since $\bigoplus _{i\in \mathcal{A}} P_i$ is a countably generated projective module, it is a direct summand of $R^{(\omega)}$. This finishes the proof of the claim.
\end{Proof}

Now we want to interpret these results in terms of \emph{lifting of idempotents} modulo an ideal. The following fact characterizing when two idempotents are conjugated by a unit will be useful.

\begin{Lemma} \cite[Exercise~21.15]{lam}\label{ferran} Let $S$ be a ring, let $e$ and $g$ be idempotents of $S$. Then there exists a  unit $u\in S$ such that $ueu^{-1}=g$ if and only if $eS\cong gS$ and $(1-e)S\cong (1-g)S$.
\end{Lemma}

\begin{Lemma} \label{liftingidempotents} Let $S$ be a ring, and let $J$ be an ideal of $S$ contained in $J(S)$. Let $e=e^2\in S/J$. Assume that $S_S\cong X_1\oplus X_2$ and $X_1/X_1J\cong e(S/J)$ and $X_2/X_2J\cong (1-e)(S/J)$.

Then there exists $f^2=f\in S$ such that $f+J=e$ and $fS\cong X_1$ and $(1-f)S\cong X_2$.
\end{Lemma}

\begin{Proof} Let $g^2=g\in S$ be such that there are isomorphisms $f_1\colon gS \to X_1$ and $f_2\colon (1-g)S\to X_2$. Then restriction of $f_1$ and $f_2$ give isomorphisms $gJ\cong X_1J$ and $(1-g)J\cong X_2J$, so that $$(g+J)S/J\cong gS/gJ\cong X_1/X_1J$$ and $$(1-g+J)S/J\cong (1-g)S/(1-g)J\cong X_2/X_2J.$$

By Lemma~\ref{ferran}, there exists $v$ a unit of $S/J$ such that $v(g+J)v^{-1}=e$. Since $J$ is contained in $J(S)$, any $u\in S$ such that $u+J=v$ is also a unit of $S$. Pick such an invertible element $u$, and set $f=ugu^{-1}$. Then $f$ is an idempotent element of $S$ with $f+J=v(g+J)v^{-1}=e$. Moreover, as $f$ is conjugated of $g$ we can apply  Lemma~\ref{ferran} to deduce  that $fS\cong X_1$ and $(1-f)S\cong X_2$.
\end{Proof}

\begin{Cor} \label{localontobij} 
Let $R$ be a semiperfect ring. Let $M$ be a finitely generated right $R$-module with endomorphism ring $S$. Assume that  onto endomorphisms of $M$ are bijective and that $M$ has no non-zero projective direct summands.  Then:
\begin{itemize}
    \item[(i)] $S$ is semilocal and the ideal $J=M \Hom_R (M,R)$ is contained in $J(S)$.     
    \item[(ii)] If the direct summands of $R^{(\omega)}\oplus M$ are direct sum of finitely generated modules, then idempotents of $S/J$ can be lifted to idempotents of $S$. If, in addition, $J(S)/J$ is a nil-ideal then $S$ is semiperfect.
 \end{itemize}   
\end{Cor}
\begin{Proof}
$(i)$. By hypothesis and Lemma~\ref{tracelocal}, $f\in J$ implies that $f(M)\subseteq MJ(R)$, so $f(M)$ is superfluous in $M$. Then $1-f$ is onto, hence bijective. This implies that $J\subseteq J(S)$. The ring $S$ is semilocal because $M/MJ(R)$ is a semisimple module of finite length, so we can apply the results in \cite{HS}.

$(ii)$. By Proposition~\ref{criteria}, any direct summand of $S/J$ can be lifted to a direct summand of $S$. By $(i)$, $J\subseteq J(S)$, and we deduce that  idempotents can be lifted modulo $J$ from Lemma~\ref{liftingidempotents}.

If $J(S)/J$ is nil, then idempotents of $S/J(S)$ can be lifted to $S/J$. Hence, to $S$. Since, by $(i)$, $S$ is semilocal, we deduce that $S$ is semiperfect.
\end{Proof}

\section{The endomorphism ring of a finitely generated module}\label{s:endo}

As implicit in Corollary~\ref{localontobij}, the condition that onto endomorphisms of a finitely generated module are bijective is going to play an important role in our study. In this section, we recall some more or less well-known facts on when this happens, and some of the consequences it has. 

We start with some well-known facts about the endomorphism ring of a finitely generated module over a commutative ring. We will use them throughout the paper, sometimes without previous acknowledgment.

\begin{Lemma} \label{basicfg}
Let $R$ be a commutative ring. Let $M$ be a finitely generated $R$-module with endomorphism ring $S$. Then:
\begin{itemize}
    \item[(i)] \emph{(Vasconcelos \cite[Proposition~1.2]{Vasconcelos})} Any onto endomorphism of $M$ is bijective.
    \item[(ii)] The central ring extension $R/\mathrm{ann}_R(M)\hookrightarrow S$ is integral (that is, any element of $S$ satisfies a monic polynomial with coefficients in $R/\mathrm{ann}_R(M)$).
    \item[(iii)] Let $L$ be a two-sided ideal of $R$ containing $\mathrm{ann}_R(M)$ and let $I=\{f\in S\mid f(M)\subseteq ML\}$. If $M$ can be generated by $\ell$ elements, then for every $f\in I$, there exist $b_1,\dotsc,b_\ell\in L$ such that $f^\ell+f^{\ell-1}\overline b_\ell+\dots+f\overline b_2+\overline b_1=0$, where $\overline b_i$, $i=1,\dotsc,\ell$, means the elements $R/\mathrm{ann}_R(M)$ viewed inside $S$.
\end{itemize}
\end{Lemma}
\begin{Proof}
Statements $(ii)$ and $(iii)$ are applications of the determinant trick. To prove $(ii)$, let $M=m_1R+\cdots +m_\ell R$ and let $f\in S$. Then there exists a matrix $A\in M_\ell(R)$ such that $(f \mathrm{Id}_\ell-A)\left(\begin{smallmatrix}
m_1\\ \vdots \\m_\ell\end{smallmatrix}\right) =0$. Then we deduce that $\mathrm{det}(f \mathrm{Id}_\ell-A)m_i=0$ for any $i=1,\dots ,\ell$. Hence, $ \overline{\mathrm{det}(f \mathrm{Id}_\ell-A)} \in (R/\mathrm{ann} _R(M))[f]$ is zero. Hence, $f$ is an integral element over $R/\mathrm{ann} _R(M)$. 

To prove $(iii)$, just observe that the coefficients of the polynomial $\overline{\mathrm{det}(f \mathrm{Id}_\ell-A)}$ satisfy the claimed properties if $f\in I$.
\end{Proof}

\begin{Remark} \label{rem-basicfg}
Let $R$ be a commutative ring, and let $\Lambda$ be a (non-necessarily commutative) module-finite $R$-algebra. Let $M$ be a finitely generated right $\Lambda$-module. Then the conclusions of Lemma~\ref{basicfg} also hold for $S=\End_\Lambda(M)$. $M_\Lambda$ being finitely generated and $\Lambda$ being module-finite implies that $M_R$ is also finitely generated as an $R$-module. Moreover, $\End_\Lambda(M)$ is a subring of $\End_R(M)$.
\end{Remark}

Goodearl in \cite{repetitive} characterized which rings satisfy that any onto endomorphism of a finitely generated module is bijective. The following result follows easily from such characterization, and it includes an extension of Lemma~\ref{basicfg} and Remark~\ref{rem-basicfg}.

\begin{Prop} \label{fgendomorphism}
Let $R$ be a commutative ring, and let $\Lambda$ be a (non-necessarily commutative) module-finite $R$-algebra. Let $M$ be a finitely generated right $\Lambda$-module with endomorphism ring $S=\End_\Lambda(M)$. Then
\begin{enumerate}
    \item[(i)] Every finitely generated right $S$-module $X$ satisfies that any onto endomorphism of $X$ is bijective. 
    \item[(ii)] If $S$ is semilocal, then so is the endomorphism ring of any finitely generated right (or left) $S$-module. In particular, if $\Lambda$ is semilocal, then so is $S$.
    \item[(iii)] $I=\{f\in\End_S(X)\mid f(X)\subseteq XJ(S)\}$ is a two-sided ideal of $\End_S(X)$ contained in $J(\End_S(X))$. In particular, $\{f\in S\mid f(M)\subseteq MJ(\Lambda)\}$ is a two-sided ideal of $S$ contained in $J(S)$.
\end{enumerate}
\end{Prop}
\begin{Proof} 
$(i)$. Since  $M_n(S)\cong \End_\Lambda(M^n)$ for any $n\ge 1$, we deduce from Lemma~\ref{basicfg}(ii) and Remark~\ref{rem-basicfg} that the extension $R/\mathrm{ann}_R(M)\to M_n(S)$ is integral for any $n\ge 1$. Now, the conclusion follows from \cite[Corollary~6]{repetitive}.

$(ii)$. The first part of the result is a particular case of \cite[Proposition~3.2]{FH}. Taking $M=\Lambda$ gives the second part of the statement.

$(iii)$. By Remark~\ref{ideal-end}, $I$ is a two-sided ideal of $S$, so it is enough to show that $1-f$ is bijective for every $f\in I$. 

Let $f\in I$. Then $S=f(X)+(1-f)(X)=XJ(S)+(1-f)(X)$. By Nakayama's Lemma, $X=(1-f)(X)$. By (i), $1-f$ is bijective.
\end{Proof}

\begin{Def}
A ring homomorphism $f\colon R\to S$ is said to be \emph{local} if, for any $r\in R$, $f(r)$ is invertible in $S$ if and only if $r$ is invertible in $R$.
\end{Def}

\begin{Lemma}\label{radical}
Let $R$ be a commutative ring. Let $R\subseteq S$ be a central integral extension of rings, where $S$ is not necessarily commutative. Then, 
\begin{itemize}
    \item[(i)] For any $s\in S$, the inclusion $R[s]\hookrightarrow S$ is a local ring homomorphism. Therefore $J(S)\cap R[s]\subseteq J(R[s])$.
    \item[(ii)] If $R$ is a ring of Krull dimension $0$, then $J(S)$ is a nil-ideal.
\end{itemize}
\end{Lemma}
\begin{Proof}
$(i)$. Let $s\in S$. To show that the inclusion $R[s]\subseteq S$ is a local homomorphism, we need to prove that if  $x\in R[s]$ is invertible in $S$, then the inverse is in $R[s]$. This is a well-known fact for integral extensions of commutative rings. We recall the argument to see that it also works in our situation.

Since $R\subseteq S$ is an integral extension, the element $x^{-1}$ satisfies   $x^{-n}-r_1x^{-n+1}-\cdots -r_n=0$ for suitable $r_1,\dots ,r_n \in R$. Since $x^{-1}$ is invertible, we may assume $r_n\neq 0$. Therefore, multiplying by $x^n$ we obtain the equality $1=(r_1+r_2x+\cdots +r_nx^{n-1})x$. Therefore, $x^{-1}=r_1+r_2x+\cdots +r_nx^{n-1}\in R[s]$.
	
The second part of the statement follows immediately from the first part.

$(ii)$. Assume now that $R$ has Krull dimension $0$. Let $s\in J(S)$. By $(i)$, $s\in J(R[s])$. Since $R\subseteq R[s]$ is an integral extension of commutative rings, also $R[s]$ is a ring of Krull dimension $0$. Hence, $J(R[s])$ coincides with the nilradical of $R[s]$. So we deduce that $s$ is nilpotent. This shows that $J(S)$ is a nil-ideal.
\end{Proof}

\begin{Lemma}\label{jacobson} 
Let $R$ be a commutative ring, and let $\Lambda$ be a module-finite $R$-algebra. If $M$ is a finitely generated right $\Lambda$-module, then $\End_\Lambda(M)J(R)\subseteq J(\End_\Lambda(M))$.
\end{Lemma}
\begin{Proof} 
Let $f_1,\dots,f_k\in \End_\Lambda(M)$ and let $r_1,\dots,r_k\in J(R)$. Since $M_R$ is finitely generated, by Lemma~\ref{basicfg}, to show that $1-\sum_{i = 1}^k f_ir_i$ is invertible it suffices 
to show that $1 - \sum_{i = 1}^kf_ir_i$ is an onto $\Lambda$-endomorphism of $M$. To this aim, notice that $M=(1-\sum_{i = 1}^kf_ir_i)(M)+(\sum_{i = 1}^kf_ir_i)(M)=(1-\sum_{i = 1}^kf_ir_i)(M)+MJ(R)$. Since $MJ(R)$ is a small $R$-submodule of $M$, we deduce that $M=(1-\sum_{i = 1}^kf_ir_i)(M)$.
\end{Proof}

\begin{Lemma} \label{maximal}
Let $R$ be a commutative ring, and let $\Lambda$ be a module-finite $R$-algebra. Let $M$ be a non-zero finitely generated right $\Lambda$-module with endomorphism ring $S=\End_\Lambda(M)$. Let $\varphi\colon R\to S$ denote the canonical homomorphism. Let $\fn$ be a two-sided maximal ideal of $S$. Then
\begin{enumerate}
    \item[(i)] $\fm=\varphi^{-1}(\fn)$ is a maximal ideal of $R$.
    \item[(ii)] $I=\{f\in S\mid f(M)\subseteq  M\fm\}$ is a two-sided ideal of $S$ contained in $\fn$ and satisfying $\fm=\varphi^{-1}(I)$.
\end{enumerate}
\end{Lemma}
\begin{Proof}
$(i)$. By Lemma~\ref{basicfg} and Remark~\ref{rem-basicfg}, the morphism $\varphi\colon R\to S$ induces a central integral extension $R/\mathrm{ann}_R(M)\to S$. Then $S/\fn$ is a simple ring and $R/\varphi^{-1}(\fn)\to S/\fn$ is also a central integral extension. The center of the simple ring $S/\fn$ is always a field $K$ (the two-sided ideal generated by a non-zero element $x$ in the center must be the total and, being the element central, this ideal is just the set of multiples of $x$; this is to say that $x$ is invertible in $S/\fn$ because $x$ is central, so is its inverse). We deduce that $R/\varphi^{-1}(\fn)\to K$ is an integral extension of commutative rings. Hence, $R/\varphi^{-1}(\fn)$ is a field (cf. Lemma~\ref{radical}) and $\fm=\varphi^{-1}(\fn)$ is a maximal ideal of $R$.

$(ii)$. By Remark~\ref{ideal-end}, $I$ is a two-sided ideal of $S$. Clearly, $\fm\subseteq\varphi^{-1}(I)$. Since $\mathrm{ann}_R(M)$ is contained in the maximal ideal $\fm$ of $R$ and $M_R$ is finitely generated, $M_\fm\neq\{0\}$. By Nakayama's Lemma, $M\neq M\fm$, so we conclude that $1\notin\varphi^{-1}(I)$. Since $\fm$ is a maximal ideal contained in $\varphi^{-1}(I)$, we deduce that $\fm=\varphi^{-1}(I)$. 

By Lemma~\ref{basicfg}(iii), $(I+\fn)/\fn$ is a two-sided nil-ideal of $S/\fn$. Since $S/\fn$ is simple, we deduce that $(I+\fn)/\fn$ is the zero ideal, so $I\subseteq\fn$, as claimed.
\end{Proof}

\subsection{Finitely generated torsion-free modules.}

Let $R$ be a domain with field of fractions $Q$. An $R$-module $M$ is \emph{torsion-free} if the natural map $M\to M_Q$ is injective, where $M_Q$ denotes the localization of an $R$-module $M$ at $R\setminus\{0\}$. Sometimes we will also use the tensor product notation $M\otimes_R Q$.

It is well known that, over commutative rings, $\mathrm{Hom}$ and localization commute when the first variable of the $\mathrm{Hom}$ is a finitely presented module. The following lemma shows that this is true  for finitely generated torsion-free modules. This result is crucial throughout the paper. For a similar result (and proof) the reader can check Warfield's paper \cite[Lemma~3.5]{warfield}.

\begin{Lemma}\label{isofg}
Let $R$ be a domain, and let $\Lambda$ be an $R$-algebra. Let $\Sigma$ be a multiplicative subset of $R$. If $M$ is a finitely generated right $\Lambda$-module which is torsion-free as an $R$-module, then the canonical injective homomorphism
    \[\varphi\colon\Hom_\Lambda(M,N)\otimes_R R_\Sigma\to\Hom_{\Lambda\otimes_R R_\Sigma}(M\otimes_R R_\Sigma,N\otimes_R R_\Sigma)\]
is an isomorphism for every $\Lambda$-module $N$ which is torsion-free as an $R$-module.
\end{Lemma}
\begin{Proof} 
Note that,  the $\Lambda$-modules that are torsion-free as $R$-modules can be seen as submodules of its localization at $\Sigma$. 

Let $f\in \Hom_{\Lambda\otimes_{R}R_\Sigma}(M\otimes _RR_\Sigma,N\otimes _RR_\Sigma)$. Let $m_1,\dots ,m_n$ be a set of generators of $M$. Then, there exists $r\in \Sigma$, such that, for any $i\in \{1,\dots ,n\}$, $f(m_i)=n_i \otimes \frac 1r$ for $n_i\in N$. Therefore $g=fr|_{M} \in \Hom_\Lambda (M,N)$, so that $f=\varphi (g\otimes \frac 1r)$. This shows that $\varphi$ is onto, as claimed in the statement.
\end{Proof}

Let $\Lambda$ be an algebra over a commutative domain $R$. We say that a short exact sequence 
    \begin{equation}\label{eq:spl}
        \begin{tikzcd} 0\rar & L\rar & M\rar & N\rar & 0 \end{tikzcd}\tag{$*$}
    \end{equation}
of right $\Lambda$-modules is \emph{locally split} if it is split when we apply the functor $R_\fm \otimes _R-$ for any maximal ideal $\fm$ of $R$. Therefore, if $X$ is a left $\Lambda$-module we obtain an exact sequence of $R$-modules
    \[\begin{tikzcd} 0\rar & K\rar & L\otimes_\Lambda X\rar & M\otimes_\Lambda X\rar & N\otimes_\Lambda X\rar & 0. \end{tikzcd}\]
As \eqref{eq:spl} is locally split, $R_\fm \otimes _R K=0$ for any maximal ideal $\fm$ of $R$. Therefore $K=0$. This proves that a locally split exact sequence of right $\Lambda$-modules is pure.

Now we show that the finitely generated $\Lambda$-modules that are torsion-free as $R$-modules are projective with respect to the class of locally split exact sequences of right $\Lambda$-modules that are torsion-free as $R$-modules.

\begin{Lemma} \label{locallysplit}
Let $R$ be a commutative domain, and let $\Lambda$ be an $R$-algebra. Let $L,M,N, X$ be right $\Lambda$-modules which are torsion-free as $R$-modules. Assume that there is a short exact sequence,
    \begin{equation}\label{eq:ses}
        \begin{tikzcd} 0\rar & L\rar{f} & M\rar{g} & N\rar & 0. \end{tikzcd}
    \end{equation}
Then:
\begin{enumerate}
    \item[(i)] If $X_\Lambda$ is finitely generated, and \emph{(\ref{eq:ses})} is a locally split exact sequence, then the map $\mathrm{Hom}_\Lambda (X,g)\colon \Hom_\Lambda(X,M)\to \Hom_\Lambda(X,N)$ is  onto. Equivalently, finitely generated torsion-free modules are projective with respect to the class of locally split exact sequences of torsion-free modules.
    \item[(ii)] If $N$ is a direct summand of a direct sum of finitely generated $\Lambda$-modules that are torsion-free as $R$-modules, then the exact sequence \emph{(\ref{eq:ses})} splits if and only if it is locally split.
\end{enumerate}
\end{Lemma}
\begin{Proof}
$(i)$. If $g_\fm\colon M_\fm\to N_\fm$ is a splitting epimorphism for every maximal ideal $\fm$ of $R$. Then, $\Hom_{\Lambda_\fm}(X_\fm, g_\fm)\colon\Hom_{\Lambda_\fm}(X_\fm,M_\fm)\to\Hom_{\Lambda_\fm}(X_\fm, N_\fm)$ is onto for every maximal ideal $\fm$ of $R$. By Lemma~\ref{isofg} we may identify $\Hom_{\Lambda_\fm}(X_\fm, g_\fm)$ with $R_\fm\otimes _R \mathrm{Hom}_\Lambda (X,g)$. Therefore, $R_\fm\otimes _R \mathrm{Hom}_\Lambda (X,g)$ is onto for every maximal ideal $\fm$ of $R$, so $\Hom_{\Lambda}(X, g)$ is onto, as desired.

$(ii)$.  As the class of modules that are projective with respect to the class of locally split exact sequences is closed under direct sums and direct summands, we  deduce the statement $(2)$ from $(1)$. 
\end{Proof}

\begin{Cor}\label{loc.proj}
Let $R$ be a commutative domain, and let $\Lambda$ be an $R$-algebra. Let $N$ be a direct summand of a direct sum of finitely generated right $\Lambda$-modules which are torsion-free as  $R$-modules. Then $N$ is a projective $\Lambda$-module if and only if $N_\fm$ is a projective $\Lambda_\fm$-module for every maximal ideal $\fm$ of $R$.
\end{Cor}
\begin{Proof}
Let $f\colon F\to N$ be an epimorphism of $\Lambda$-modules for some  free $\Lambda$-module $F$. By Lemma~\ref{locallysplit}(ii), $f$ is a splitting epimorphism if and only if $f_\fm$ is a splitting epimorphism for every maximal ideal $\fm$ of $R$. This is equivalent to the statement.
\end{Proof}

The following easy lemma will be useful throughout the paper. 

\begin{Lemma}\label{d}
Let $R$ be a commutative domain, and let $\Lambda$ be an $R$-algebra. Let $\Sigma$ be a multiplicative subset of $R$. Let $N$ be a right $\Lambda$-module which is torsion-free as an $R$-module. If $M$ is a finitely generated $\Lambda$-submodule of $N_\Sigma$, then $dM$ is isomorphic to a $\Lambda$-submodule of $N$ for  some $d\in \Sigma$.
\end{Lemma}
\begin{Proof}
Note that, since  $N$ is torsion-free as an $R$-module, it can be seen as a $\Lambda$-submodule of $N_\Sigma$.

Let $m\in M$ and let $\{n_1/s_1,\dotsc,n_t/s_t\}\subseteq N_\Sigma$ be a finite set of $\Lambda$-generators of $M$. Then,
    \[m=\frac{n_1}{1}\frac{a_1}{s_1}+\dotsb+\frac{n_t}{1}\frac{a_t}{s_t}\in N_\Sigma\]
for some $a_1,\dotsc,a_t\in \Lambda$, and  $s_1,\dotsc,s_t\in\Sigma$. Multiplying by $d=s_1\dotsb s_t\in\Sigma$, 
    \[dm=\frac{n_1 a_1'+\dotsb+n_ta_t'}{1}\in\lambda(N),\]
where $\lambda\colon N\to N_\Sigma$ denotes the localization map. Therefore, $dM\subseteq\lambda(N)\cong N$. 
\end{Proof}

\begin{Remark} We note that Lemmas~\ref{isofg} and \ref{d} could be stated for general commutative rings, provided that the multiplicative set $\Sigma$ consists of non-zero divisors of $R$. In this general context, we mean that a module $M$ is torsion-free if $md=0$, for $m\in M$ and $d$ a nonzero divisor of $R$, implies that $m=0$.
    \end{Remark}

Let $R$ be a commutative domain, and let $Q$ denote the field of fractions of $R$. The \textit{rank} of an $R$-module $M$ is the dimension of $M_Q$ as a $Q$-vector space.

\begin{Remark} \label{rem-almost-free}
\textit{Let $R$ be a commutative domain with field of fractions $Q$. Let $N$ be a torsion-free $R$-module of rank $n$. Then $N$ contains a free $R$-module of rank $n$.} Note that, since $N$ is torsion-free of rank $n$, we can see $N$ as an essential submodule of its injective hull $E(M)\cong Q^n$. On the other hand, $R^n$ has the same injective hull $Q^n$. Therefore, by Lemma~\ref{d}, since $R^n$ is a finitely generated $R$-submodule of $N_Q=Q^n$, $dR^n\subseteq N$ for some non-zero element $d\in R$.
\end{Remark}

\begin{Lemma} \label{samerank} 
Let $R$ be a commutative domain with field of fractions $Q$. Let $M$ and $N$ be non-zero torsion-free $R$-modules of the same finite rank $n$. Let $f\colon M\to N$ be an $R$-module homomorphism. Then,
\begin{itemize}
    \item[(i)] $f$ is onto if and only if it is bijective;
    \item[(ii)] $f$ is injective if and only if $\im f$ is an essential submodule of $N$.
\end{itemize}
\end{Lemma}
\begin{Proof}
Applying the functor $-\otimes _RQ$ to the exact sequence,
    \[\begin{tikzcd}
        0\rar & \Ker f\rar & M\rar{f} & N\rar & 0
    \end{tikzcd}\]
we get the exact sequence
    \[\begin{tikzcd}
        0\rar & \Ker f\otimes _RQ\rar & M\otimes _RQ\rar{f\otimes _R\mathrm{Id}_Q} & N\otimes _RQ\rar & 0
    \end{tikzcd}\]
in which the $Q$-module homomorphism $f\otimes _R\mathrm{Id}_Q$ is an isomorphism because it is an onto morphism between $Q$-vector spaces of the same finite dimension. Therefore, $\Ker f\otimes _RQ=0$ and $\Ker f$ being torsion-free, it must be zero. Hence, $f$ is bijective.

The statement $(ii)$ is easy to prove.
\end{Proof}

\begin{Prop} \label{algebrafiniterank}
Let $R$ be a commutative domain with field of fractions $Q$, and let $\Lambda$ be a torsion-free $R$-algebra of finite rank. Let $M_\Lambda$ be a non-zero finitely generated right $\Lambda$-module which is torsion-free as an $R$-module. Then,
\begin{itemize}
    \item[(i)] any onto endomorphism of $M$ is bijective;
    \item[(ii)] if $J(R) \subseteq J(\Lambda)$ then $\End_\Lambda(M)J(R)\subseteq J(\End_\Lambda (M))$;
    \item[(iii)] if $\Lambda$ is semilocal then so is $\End_\Lambda (M)$.
\end{itemize}
\end{Prop}
\begin{Proof} 
$(i)$. Since $M$ is a homomorphic image of a finite number of copies of $\Lambda$, $M$ is a torsion-free $R$-module of finite rank. Since a $\Lambda$-module homomorphism is also an $R$-module homomorphism, the statement follows from Lemma~\ref{samerank}(i).

$(ii)$. Our hypothesis allows us to repeat the proof of Lemma~\ref{jacobson} to conclude the statement.

$(iii)$. If $\Lambda$ is semilocal, then $MJ(\Lambda)$ is a small submodule of $M$ and $M/MJ(\Lambda)$ is semisimple artinian of finite length. By $(i)$, any onto endomorphism ring of $M$ is bijective. As an application of \cite{HS}, we deduce that $\End_\Lambda (M)$ is semilocal.
\end{Proof}

\subsection{Domains of finite character and $h$-local domains.}

\begin{Def}
A commutative domain $R$ is said to have \emph{finite character} if any non-zero element is contained only in a finite number of maximal ideals. A commutative domain of finite character is said to be \emph{$h$-local} if, in addition, any non-zero prime ideal of $R$ is contained in a unique maximal ideal.
\end{Def}

The domain $R$ has finite character if $R/I$ is a semilocal ring for every non-zero ideal $I$ of $R$. If $R$ is $h$-local then it also satisfies that  $R/\mathfrak{p}$ is a local domain for every non-zero prime ideal $\mathfrak{p}$ of $R$.

If $R$ is a domain of Krull dimension $1$, then the notion of $h$-locality and being of finite character coincide.  They are equivalent to saying that, for every non-zero ideal $I$ of $R$, $R/I$ is a semiperfect ring whose Jacobson radical is a nil ideal. Because $R/I$ is a ring of Krull dimension $0$ with only a finite number of maximal ideals, $J(R/I)$ coincides with the nilradical of $R/I$.

Matlis introduced the notion of $h$-local domain in \cite{matlis3} as a generalization of local domains and noetherian domains of Krull dimension $1$. In the next result, we recall the key property of $h$-local domains that we will use throughout the paper.

\begin{Lemma}\label{hlocal} 
A commutative domain $R$ is $h$-local if and only if $R/I$ is semiperfect for every non-zero proper ideal $I$ of $R$. In particular, if $R$ is $h$-local, $I$ is a non-zero proper ideal of $R$, and $\{\fm_1,\dots,\fm_\ell\}$ is the finite set of maximal ideals of $R$ containing $I$, then
\begin{enumerate}
    \item[(i)] $I=I_1\dotsb I_\ell$, where $I_i=I_{\fm_i}\cap R$, so $I_1,\dotsc,I_\ell$ are pairwise comaximal ideals of R such that each $I_i$ is contained in exactly one maximal ideal of $R$.
    \item[(ii)] the canonical map $R/I\to (R/I)_{\fm_1}\times \cdots \times (R/I)_{\fm_\ell}$, given by the localization at each component, is an isomorphism, and $(R/I)_{\fm_i}=(R/I_i)_{\fm_i}$ for $i=1,\dots, \ell$. 
\end{enumerate}
\end{Lemma} 
\begin{Proof}
The first part is \cite[Theorem~4.9]{bazzoni}. Statement (i) is proven in \cite[Lemma 5.1]{olberding2}. To prove (ii), notice that, by the Chinese Remainder Theorem, there is an isomorphism
    \[\begin{tikzcd}
        \varphi\colon R/I\rar & R/I_1\times \cdots \times R/I_\ell.
    \end{tikzcd}\]
Since each $I_i$ is contained only in the maximal ideal $\fm _i$, we deduce that $R/I_i$ is a local ring and that $R/I_i=(R/I_i)_{\fm _i}$. So that we have the isomorphisms claimed.
\end{Proof}

\begin{Lemma} \label{endfg}
Let $R$ be a commutative domain with field of fractions $Q$. Let $M$ be a non-zero finitely generated torsion-free module with endomorphism ring $S$. Then
\begin{itemize}
    \item[(i)] $S$ is a torsion-free $R$-module containing $R$, and $S$ is integral over $R$. The restriction of the embedding $R\subseteq S$ gives an embedding $J(R)\subseteq J(S)$.
    \item[(ii)] For any $f\in S$, there is some non-zero $d\in R$ such that $df\in M \Hom_R(M,R)$. In particular, $M \Hom_R(M,R)\neq \{0\}$.
    \item[(iii)] $S_Q\cong M_n(Q)$ where $n$ is the rank of $M$. So $S$ is a subring of the simple artinian algebra $M_n(Q)$. 
    \item[(iv)] For any non-zero two-sided ideal $I$ of $S$, $I\cap R\neq \{0\}$.
    \item[(v)] $S_\fm\cong \End_{R_\fm}(M_\fm)$ is a semilocal ring for all maximal ideals $\fm$ of $R$.
\end{itemize}
If, in addition, $R$ is $h$-local, then
\begin{itemize}
    \item [(vi)] Let $I$ be any non-zero two-sided ideal of $S$. Then $S/I$ is a semilocal ring. Moreover, $I_{\fm}=S_\fm$ for almost all maximal ideals $\fm$ of $R$. If $R$ has Krull dimension $1$ then $J(S/I)$ is a nilideal, so $S/I$ is semiperfect.
    \item[(vii)] If $I$ is the trace ideal of a non-zero projective right $S$-module, then $I$ is contained only in a finite number of trace ideals of projective right $S$-modules. 
\end{itemize}
\end{Lemma}
\begin{Proof} 
Recall that, since $R$ is commutative, there is a ring morphism $\varphi\colon R\to S$ where, for any $r\in S$, $\varphi (r)$ is the morphism given by multiplication of $r$. Since $M$ is a faithful $R$-module, $\varphi$ is injective, and we may identify $R$ as a subring of $S$. We will use this identification throughout the proof.
	
$(i)$. Let $f\in S$ be such that there is some non-zero $r\in R$ with $rf=0$. Equivalently, $rf(M)=\{0\}$ so that $f(M)=\{0\}$ because $M$ is torsion-free. 

The extension $R\subseteq S$ is integral because of Lemma~\ref{basicfg}. Repeating the arguments in Lemma~\ref{jacobson}, we deduce that $J(R)\subseteq J(S)$.

$(ii)$. Since $M$ is finitely generated and torsion-free, we can see $M$ as an essential submodule of its injective hull $E(M)\cong Q^n$ for some $n\ge 1$. On the other hand, $R^n$ has the same injective hull $Q^n$. Therefore, since $M$ is finitely generated, by Lemma~\ref{d}, $gM\subseteq R^n$ for some non-zero $g\in R\subseteq S$.

Now we proceed with the proof of the statement. The claim is clear for the zero endomorphism.  Pick $f\in S\setminus \{0\}$. Since $Q^n$ is an injective $R$-module, there exists $h\in \Hom_R(R^n, Q^n)\cong M_n(Q)$ such that the diagram
    \[\begin{tikzcd}[column sep=huge,row sep=huge]
        M \arrow[swap,dr, "g"] \arrow[r, "f"]& M\arrow[r, hook] &Q^n \\
        & R^n  \arrow[swap,ur, "h"]
    \end{tikzcd}\]
commutes. Since $R^n$ is finitely generated and $M$ is essential in $Q^n$, by Lemma~\ref{d}, $dh(R^n)\subseteq M$ for some non-zero $d\in R$, that is, $dh \in \Hom_R(R^n, M)$. Therefore, $df=(dh)\circ g\in M \Hom_R(M,R)$.

$(iii)$. By Lemma~\ref{isofg}, $S_Q\cong \End_Q(M_Q)\cong M_n(Q)$, where $n$ is the rank of $M$. By~$(i)$, $S$ is torsion-free as an $R$-module, so the localization map $S\to S_Q$ is injective.

$(iv)$. By $(iii)$ we can identify $S$ with a subring of its central localization at the non-zero divisors of $R$ which is the simple ring $M_n(Q)$. If $I$ is a non-zero ideal of $S$ then $I_Q=IM_n(Q)$ is  a non-zero ideal of $M_n(Q)$. Therefore $IM_n(Q)=M_n(Q)$. Therefore, there is some $f\in I$ and some non-zero $d\in R$ such that $\frac fd=\mathrm{Id}_n$. This implies that $f$ is the endomorphism of $M$ given by multiplication by $d$. So $f\in R\cap I$.

$(v)$. By Lemma~\ref{isofg}, $S_\fm\cong \End_{R_\fm}(M_\fm)$. By Lemma~\ref{basicfg}, $S_\fm$ is semilocal.

From now on, we are assuming that $R$ is $h$-local.

$(vi)$. Let $I$ be a non-zero two-sided ideal of $S$. The inclusion $\varphi \colon R\to S$ induces an injective ring homomorphism $R/(I\cap R)\to S/I$. By $(iv)$, $I\cap R\neq \{0\}$ and, since $R$ is $h$-local there are only finitely many maximal ideals of $R$, say $\fm_1,\dots,\fm_t$, containing $I\cap R$. By Lemma~\ref{hlocal}, the canonical homomorphisms
    \[R/(I\cap R)\longrightarrow (R/(I\cap R))_{\fm_1}\times\dots\times (R/(I\cap R))_{\fm_t}\]
and
    \[S/I\longrightarrow (S/I)_{\fm_1}\times\dots\times (S/I)_{\fm_t}\]
are isomorphisms. Therefore, we have a commutative diagram
    \[\begin{tikzcd}[column sep=huge,row sep=huge]
        R/(I\cap R)\arrow[d,"\cong"]\arrow[hookrightarrow]{r}&S/I\arrow[d,"\cong"]\\
        \left(	R/(I\cap R)\right)_{\fm_1}\times \cdots \times 	\left(	R/(I\cap R)\right)_{\fm_\ell}\arrow[hookrightarrow]{r}& \left(	S/I\right)_{\fm_1}\times \cdots \times 	\left(	S/I\right)_{\fm_\ell}
    \end{tikzcd}\]
By $(v)$, $S/I$ is a finite product of semilocal rings, so it is a semilocal ring.

Since $I\cap R$ is only contained in a finite number of maximal ideals of $R$, $(I\cap R)_\fm=R_\fm$ for almost all maximal ideals $\fm$ of $R$. Therefore, $I_\fm =S_\fm$ for almost all maximal ideals $\fm$ of $R$.

Assume now that $R$ has Krull dimension $1$. As $R/(I\cap R)\subseteq S/I$ is an integral extension, and $R/(I\cap R)$ is a ring of Krull dimension $0$,  Lemma~\ref{radical} implies that $J(S/I)$ is a nil-ideal. As idempotents can be lifted modulo a nilideal we deduce that $S/I$ is a semiperfect ring. 

The statement $(vii)$ follows from $(vi)$ because semilocal rings only have a finite number of trace ideals of projective modules.	
\end{Proof}

\section{Local domains of Krull dimension one}\label{s:local}

Now we are ready to start our study of direct summands of direct sums of torsion-free modules. In this section, we will specialize to the case of local domains of Krull dimension one. The next lemma points out the key properties of endomorphism rings of torsion-free modules over such domains that we will need.

\begin{Lemma}  \label{localnil} 
Let $R$ be a commutative local domain of Krull dimension $1$, and with field of fractions $Q$. Let $M$ be a non-zero finitely generated torsion-free $R$-module with endomorphism ring $S$. If $M$ does not contain a non-zero free direct summand, then:
\begin{enumerate}
    \item[(i)] $J=M\Hom_R(M,R)\subseteq J(S)$ and  $J(S)/M\Hom_R(M,R)$ is a nil-ideal of $S/J$.
    \item[(ii)] If every direct summand  of $R^{(\omega)}\oplus M$ is a direct sum of finitely generated modules, then $S$ is semiperfect.
\end{enumerate}
\end{Lemma}
\begin{Proof} 
By Lemma~\ref{basicfg}, $M$ satisfies that any onto endomorphism ring of $M$ is bijective.

$(i)$. By Corollary~\ref{localontobij}, $J=M\Hom_R(M,R)\subseteq J(S)$. The rest of the statement is included in Lemma~\ref{endfg}.

$(ii)$. Follows from $(i)$ and Corollary~\ref{localontobij}.
\end{Proof}

\begin{Prop} \label{localAddRM}
Let $R$ be a commutative local domain of Krull dimension $1$, and with field of fractions $Q$. Let $M$ be a finitely generated torsion-free $R$-module. The following statements are equivalent:
\begin{itemize}
    \item [(i)] every module in $\Add(R\oplus M)$ is a direct sum of finitely generated modules;
    \item[(ii)] every direct summand of $R^{(\omega)}\oplus M$ is a direct sum of finitely generated modules;
    \item [(iii)] $\End_R(M)$ is semiperfect; and
    \item[(iv)] every finitely generated indecomposable module in $\Add(R\oplus M)$ has local endomorphism ring.
\end{itemize}
\end{Prop}
\begin{Proof}
$(i) \Rightarrow (ii)$ is clear.

$(ii) \Rightarrow (iii)$  Write $M=F\oplus M'$, with $F$ finitely generated free and $M'$ has no projective direct summands. By Lemma~\ref{localnil}(ii), $\End_R(M')$ is semiperfect. As $R$ is local, so is $\End_R(M)$.

Statement $(iii)$ is equivalent to saying that $M=Y_1\oplus \dotsb\oplus Y_m$, where each $Y_i$ has local endomorphism ring for $i\in\{1,\dotsc,n\}$ (cf. \cite[Theorem 3.14]{libro}). 

In general, over a commutative ring, direct summands of direct sums of finitely generated modules with local endomorphism rings  are also direct sums of finitely generated modules with local endomorphism rings, and such decomposition is unique (cf. \cite[Corollary 2.55]{libro}). Therefore, every element in $\Add(R\oplus M)$ is isomorphic to a direct sum of copies of $R,Y_1,\dotsc,Y_m$. This proves $(i)$.

In particular, since every indecomposable in $\Add(R\oplus M)$ is isomorphic to one in $\{R,Y_1,\dotsc,Y_m\}$, we deduce $(iv)$.

Note that $(iv)$ also implies that $M$ is a direct sum of indecomposable modules with local endomorphism rings, so this implies $(iii)$.
\end{Proof}

\begin{Cor} \label{generalrank}
Let $R$ be a commutative local domain of Krull dimension $1$, and with field of fractions $Q$. The following statements are equivalent:
\begin{itemize}
    \item[(i)] the class of direct sums of finitely generated torsion-free $R$-modules is closed under direct summands;
    \item[(ii)] for any finitely generated torsion-free $R$-module $M$, every direct summand of $R^{(\omega)}\oplus M$ is a direct sum of finitely generated modules; 
    \item[(iii)] for any finitely generated torsion-free $R$-module $M$, every module in $\Add(M)$ is a direct sum of finitely generated modules; and
    \item[(iv)] every finitely generated, indecomposable, torsion-free $R$-module has local endomorphism ring.
\end{itemize}
In addition, the above equivalent conditions are fulfilled by any intermediate ring $S$ between $R$ and its integral closure.
\end{Cor}
\begin{Proof}
The equivalence of the four statements is an immediate consequence of Proposition~\ref{localAddRM}. 

To finish the proof of the statement, let $S$ be a ring such that $R\subseteq S\subseteq \overline{R}$. Notice that if $X_S$ is a finitely generated torsion-free $S$-module then it has finite rank $n$. Identifying $X$ with a submodule of $Q^n$,
$$\mathrm{End}_S(X)=\{A\in \mathrm{End}_Q(Q^n)\mid AX\subseteq X\}= \mathrm{End}_R(X).$$
In particular, $X$ is indecomposable as an $S$-module if and only if it is indecomposable as an $R$-module. 

If $S$ is finitely generated over $R$ then, also $X$ is finitely generated over $R$, so by $(iv)$, $X_S$ is indecomposable if and only if $\mathrm{End}_S(X)$ is local. 

Assume $S$ is not finitely generated over $R$, and fix $x_1,\dots ,x_r$ a family of generators of $X_S$. Then $\mathrm{End}_S(X)$ is the directed union of $\mathrm{End}_T (\sum _{i=1}^rx_iT)$ where $T$ varies between all subrings of $S$ that are finite extensions of $S$. If $X_S$ is indecomposable then so is the $T$-module $\sum _{i=1}^rx_iT$. By the previous step, $\mathrm{End}_T (\sum _{i=1}^rx_iT)$ is a local ring for any $T$. Therefore, $\mathrm{End}_S(X)$ is also a local ring.
\end{Proof}

We do not know of a characterization of  local domains (of Krull dimension $1$) such that their indecomposable finitely generated torsion-free modules have local endomorphism ring. If we restrict the condition to finitely generated, torsion-free, rank one modules, then it is equivalent to having a local integral closure, as we show in the next result.

\begin{Prop} \label{rankone}
Let $R$ be a commutative local domain of Krull dimension $1$, and with field of fractions $Q$. The following statements are equivalent:
\begin{itemize}
    \item [(i)] the class of direct sums of finitely generated, rank-one, torsion-free $R$-modules is closed under direct summands;
    \item[(ii)] for any finitely generated, rank-one, torsion-free $R$-module $M$, every direct summand  of $R^{(\omega)}\oplus M$ is a direct sum of finitely generated modules;
    \item[(iii)] every finitely generated, rank-one, torsion-free module has local endomorphism ring; and
    \item[(iv)] the integral closure of $R$ (into its field of fractions) is a local ring.
\end{itemize}
\end{Prop}
\begin{Proof} 
The equivalence of $(i)$, $(ii)$ and $(iii)$ follows from Proposition~\ref{localAddRM}.

$(iii) \Rightarrow (iv)$.  By $(iii)$ any finitely generated, integral extension of $R$ inside $Q$ is local. Let $S\subseteq Q$ be any integral extension of $R$ inside $Q$. If $S$ has two different maximal ideals $\fm_1$ and $\fm_2$, then for any $i=1,2$, choose $s_i\in \fm_i$ such that $s_1+s_2=1$. Then $R\subseteq R[s_1,s_2]$ is a finitely generated integral extension of $R$. By the first part of the proof, $R[s_1,s_2]$ is a local ring, but its maximal ideal should contain $s_1$ and $s_2$, which contradicts our assumptions. Therefore, $S$ is a local ring. 	

$(iv) \Rightarrow (iii)$. If $M$ is a finitely generated, rank-one, torsion-free module, then its endomorphism ring $S$ is a subring of $Q$, integral over $S$. Hence, it is contained in the integral closure $\overline{R}$ of $R$. Since $S\subseteq \overline{R}$ is integral and $\overline{R}$ is local, so is $S$ (cf. \cite[Theorem~9.3]{matsumura}).
\end{Proof}

In the case of noetherian domains of Krull dimension $1$ with finitely generated integral closure, P\v r\'\i hoda in \cite{P3} proved that having local integral closure (so the integral closure is a discrete valuation ring) already implies that all indecomposable finitely generated torsion-free modules have local endomorphism ring. 

The following criteria to ensure local endomorphism ring  encodes P\v r\'\i hoda's ideas in \cite[Proposition~1]{P3}. It will allow us to prove that for any noetherian local domain of Krull dimension $1$ having local integral closure is equivalent to the fact that indecomposable, finitely generated, torsion-free modules have local endomorphism ring.

\begin{Lemma} \label{tecnical_extension}
Let $R$ be a commutative local domain of Krull dimension $1$, and with field of fractions $Q$. Let $M_R$ be a finitely generated, indecomposable, torsion-free module. Assume that:
\begin{itemize}
    \item[(i)] there exists a local  over ring $T$ of $R$, of Krull dimension $1$, contained in  $Q$, such that  $MT_T$ is a direct sum of $T$-modules with local endomorphism ring;
    \item[(ii)] the conductor ideal $\mathfrak{c}=\{r\in R\mid rT\subseteq R\}$ is non-zero.
\end{itemize}
Then $M_R$ has local endomorphism ring.
\end{Lemma}
\begin{Proof}
Since $MT_T$ is a (finite) direct sum of modules with local endomorphism ring, then $\End_T (MT)$ is a semiperfect ring. Since $\mathfrak{c}$ is a non-zero proper ideal of both the rings $R$ and $T$, it is contained in their corresponding maximal ideals. Since $R$ and $T$ have Krull dimension $1$, $R/\mathfrak{c}$ and   $T/\mathfrak{c}$ are rings of Krull dimension zero.

Let $S=\End_R(M)$ and $J=\{f\in S\mid f(M)\subseteq M\mathfrak{c}\}$. By Proposition~\ref{fgendomorphism}, $J$ is a two-sided ideal of $S$ contained in $J(S)$. So it is enough to show that $S/J$ is a local ring.

Since $R\subseteq T\subseteq Q=E(R)$, $M$ can be identified with an essential submodule of $Q^n$ for some $n$. Via this identification, $S=\{A\in M_n(Q)\mid AM\subseteq M\}$ and $\End_T (MT)= \{A\in M_n(Q)\mid AMT\subseteq MT\}$, so that $S$ is a subring of the semiperfect ring $\End_T(MT_T)$. Moreover $J=\{A\in M_n(Q)\mid AM\subseteq M\mathfrak{c}\}= \{A\in M_n(Q)\mid AMT\subseteq M\mathfrak{c}\}$, so it is also a two-sided ideal of $\End_T(MT_T)$ contained in its Jacobson radical.

By Lemma~\ref{basicfg}, $R/(J \cap R) \subseteq S/J$ is an integral extension. Since $\mathfrak{c} \subseteq J \cap R$,  $R/(J \cap R)$ has Krull dimension $0$, by Lemma~\ref{radical}, we deduce that $J(S)/J$ is a nil-ideal. So the idempotents of the semisimple artinian ring $S/J(S)$ can be lifted to $S/J$. Let $e\in S$ be such that $e+J$ is an idempotent of $S/J\subseteq \End_T(MT_T)/J$. Since $\End_T(MT_T)$ is a semiperfect ring, by \cite[Proposition~27.4, Theorem~27.6]{andersonfuller}, there exists an idempotent $e'\in \End_T(MT_T)$ such that $e-e' \in J$. But then $e'M\subseteq (e'-e)M+eM\subseteq M\mathfrak{c}+M=M$. Hence $e'\in S$. Since $M$ is indecomposable, $e'$ is a trivial idempotent. Therefore, $e+J$ is also a trivial idempotent of $S/J$ and we deduce that $S/J$ is a local ring.
\end{Proof}

\begin{Lemma} \label{finite_extension}
Let $R$ be a commutative local domain of Krull dimension $1$, and with local integral closure. Let $M_R$ be a finitely generated indecomposable torsion-free module. Assume that there exists a finitely generated integral extension $T$ of $R$ such that $MT_T$ is a direct sum of $T$-modules with local endomorphism ring. Then $M_R$ has local endomorphism ring.
\end{Lemma}
\begin{Proof}
Since the integral closure of $R$ is local, so is $T$. Moreover, $T$ is also a ring of Krull dimension $1$, and since $T_R$ is finitely generated, the conductor ideal $\mathfrak{c}=\{r\in R\mid rT\subseteq R\}$ is different from zero. Therefore, the statement follows from Lemma~\ref{tecnical_extension}.
\end{Proof}

\begin{Lemma} \label{valuationring-lemma}
Let $R$ be a commutative local domain of Krull dimension $1$, and with field of fractions $Q$. Assume that the integral closure $\overline R$ of $R$ is a valuation domain. For any non-zero finitely generated torsion-free module $M$, there exists a finite integral extension $T$ of $R$ inside $Q$ such that $MT\cong T\oplus M'$ for a suitable $T$-module $M'$.
\end{Lemma}
\begin{Proof}
Let $n>0$ be the rank of $M$.   We may assume that $M$ is an essential submodule of $Q^n$ that contains $R^n$ as an (essential) submodule. Since $\overline R$ is a valuation domain, the module $M\overline{R}_{\overline{R}}$ is free (see, for example, \cite[Corollary~V.2.8]{fuchssalce}), so that $\Tr_{\overline R}(M\overline{R})=\overline R$. Therefore, there exist $f_1,\dots ,f_\ell\in \Hom_{\overline{R}}(M\overline{R}, \overline{R})$  and $m_1,\dots ,m_\ell \in M$ such that $1=f_1(m_1)+\dots+ f_\ell (m_\ell)$.

Since $\overline{R}^n$ is a submodule of $M\overline{R}$, 
    $$\Hom_{\overline{R}}(M\overline{R}, \overline{R})=\{(s_1,\dots ,s_n)\in\overline{R}^n\mid (s_1,\dots ,s_n)M\subseteq \overline{R}\}.$$ 
So that, for $i=1,\dots, \ell$, $f_i=(s^i_1,\dots ,s^i_n)$. Let $T$ be the finite integral extension of $R$ generated by the elements $\{s^i_j\}_{\substack{i=1,\dots ,\ell \\ j=1,\dots ,n }}$ 
and $\cup_{i = 1}^{\ell} f_i(G)$, where $G$ is a finite set of generators of $M$.  

By construction, $\Tr_T(MT)=T$ and then, by Lemma~\ref{tracelocal} we deduce that there exists $M'$ such that $MT\cong T\oplus M'$. This finishes the proof of the lemma.    
\end{Proof}

\begin{Prop} \label{valuationring}
Let $R$ be a commutative local domain of Krull dimension $1$, and with field of fractions $Q$. Assume that the integral closure $\overline R$ of $R$ is a valuation domain. Let $M_R$ be a finitely generated torsion-free module. Then $M$ is a direct sum of modules with local endomorphism ring.
\end{Prop}
\begin{Proof}
We are going to prove the statement by induction on the rank $n$ of $M$. If $n=1$, Lemma~\ref{valuationring-lemma} implies that $MT\cong T$. Then, since $T$ must be also a local ring, we may conclude by Lemma~\ref{finite_extension}. Now assume that $n>1$ and that the statement is proven for  modules of smaller rank. By the previous claim, there exists a finite integral extension of $T$ such that $MT\cong T\oplus M'$ for a suitable $T$-module $M'$ of smaller rank. By the inductive hypothesis, $M'$ is a $T$-module that is a direct sum of modules with local endomorphism ring. By Lemma~\ref{finite_extension}, we conclude that $M$ is a direct sum of modules with local endomorphism ring.		
\end{Proof}

\begin{Cor} \label{local} 
Let $R$ be a commutative local noetherian domain of Krull dimension $1$, and with field of fractions $Q$. Then the following statements are equivalent:
\begin{itemize}
    \item [(i)] every direct summand of a direct sum of finitely generated, torsion-free modules of rank one is a direct sum of finitely generated modules;
    \item[(ii)] for any finitely generated, torsion-free $R$-module $M$ of rank one, every direct summand of $R^{(\omega)}\oplus M$ is a direct sum of finitely generated modules;
    \item [(iii)] every finitely generated, torsion-free module of rank one has local endomorphism ring;
    \item[(iv)] the integral closure of $R$ is local;
    \item[(v)] every direct summand of a direct sum of finitely generated, torsion-free modules is a direct sum of finitely generated modules;
    \item[(vi)] for any finitely generated, torsion-free $R$-module $M$, every direct summand of $R^{(\omega)}\oplus M$ is a direct sum of finitely generated modules;
    \item[(vii)] every finitely generated, indecomposable, torsion-free module has local endomorphism ring.
\end{itemize}
\end{Cor}
\begin{Proof}
The first four statements are equivalent because of Proposition~\ref{rankone}. Since $R$ is a commutative local noetherian ring of Krull dimension $1$, its integral closure is a discrete valuation ring. Hence, by Proposition~\ref{valuationring}, all finitely generated torsion-free modules are direct sums of modules with local endomorphism rings. Now, the equivalence of the remaining statements follows from Corollary~\ref{generalrank}. This finishes the proof of the equivalence of the statements.
\end{Proof}

The reader is referred to Lemma~\ref{nummonoids} for examples that satisfy the conclusion of Corollary~\ref{local}.

\section{Local versus global direct summands}\label{s:dsummands}

The following result, for finitely generated torsion-free modules over a domain $R$, gives a general relation between the property of being locally a direct summand and being a direct summand. Afterwards, we will see how the result can be refined when we assume, in addition, that the domain is of finite character.

\begin{Prop}\label{summandk}
Let $R$ be a commutative domain, and $\Lambda$ be an $R$-algebra. Let $M,N$ be finitely generated right $\Lambda$-modules which are torsion-free as $R$-modules. If $M_\fm$ is a direct summand of $N_\fm$ for every maximal ideal $\fm$ of $R$, then there exists $k>0$ such that $M$ is a direct summand of $N^k$.
\end{Prop}
\begin{Proof}
Let $\{\fm_\alpha \} _{\alpha \in \Omega}$ denote the set of all maximal ideals of $R$ and fix $\alpha \in \Omega$. Since $M_{\fm_\alpha}$ is a direct summand of $N_{\fm_\alpha}$, there are $\Lambda_{\fm_\alpha} $-module homomorphisms $\tilde f_\alpha \colon N_{\fm_\alpha}\to M_{\fm_\alpha}$ and $\tilde g_\alpha \colon  M_{\fm_\alpha}\to N_{\fm_\alpha}$ such that $\tilde f_\alpha \circ \tilde g_\alpha =\mathrm{Id}_{M_{\fm_\alpha}}$. By Lemma~\ref{isofg}, there are $\Lambda$-homomorphisms $f_\alpha\colon N\to M$ and $g_\alpha\colon M\to N$ such that $\tilde f_\alpha=f_\alpha/s_\alpha$, $\tilde g_\alpha=g_\alpha/s_\alpha$ and $f_\alpha \circ  g_\alpha =s^2_\alpha \mathrm{Id}_M$ for some $s_\alpha\in R\setminus \fm_\alpha$.

Note that, by the definition of the elements $s_\alpha$,  $I=\sum_{\alpha \in \Omega} s^2_\alpha R=R$ because no maximal ideal of $R$ contains $I$. Therefore, there exist $k>0$, $\alpha _1,\dots, \alpha _k \in \Omega$ and $r_1,\dots ,r_k\in R$ such that $1=\sum_{i=1}^k s^2_{\alpha _i}r_i$. Hence
    $$\sum_{i=1}^kr_i f_{\alpha_i} \circ  g_{\alpha_i}=\left(\sum_{i=1}^k s^2_{\alpha_i}r_i\right)\mathrm{Id}_M=\mathrm{Id}_M,$$
and we conclude that $M$ is a direct summand of $N^k$.
\end{Proof}

The following lemma is an extension of \cite[Lemma 2.1]{GL} to the case of finitely generated torsion-free modules over commutative domains of finite character.

\begin{Lemma}\label{almostallmaximals}
Let $R$ be a commutative domain of finite character with field of fractions $Q$. Let $\Lambda$ be an $R$-algebra, and let $M,N$ be right $\Lambda$-modules which are torsion-free as $R$-modules. Assume that there exists a homomorphism of $\Lambda_Q$-modules $F\colon M_Q\to N_Q$. If $M$ is finitely generated, then
\begin{enumerate}
    \item[(i)] If $F$ is a monomorphism (resp. an epimorphism), then there is a $\Lambda$-module homomorphism $f\colon M\to N$ such that the induced homomorphism $f_\fm\colon M_\fm\to N_\fm$ is a monomorphism (resp. an epimorphism) for almost all maximal ideals $\fm$ of $R$.
\end{enumerate}
Moreover, if $N$ is also finitely generated, then
\begin{enumerate}
    \item[(ii)] If $F$ is a splitting monomorphism (resp. a splitting epimorphism), then there is a $\Lambda$-module homomorphism $f\colon M\to N$ such that the induced homomorphism $f_\fm\colon M_\fm\to N_\fm$ is a splitting monomorphism (resp. a splitting epimorphism) for almost all maximal ideals $\fm$ of $R$.
    \item[(iii)] If $F$ is an isomorphism, then there is a $\Lambda$-module homomorphism $f\colon M\to N$ such that the induced homomorphism $f_\fm\colon M_\fm\to N_\fm$ is an isomorphism for almost all maximal ideals $\fm$ of $R$.
\end{enumerate}
\end{Lemma}
\begin{Proof}
$(i)$. By Lemma~\ref{isofg}, there is a $\Lambda$-module homomorphism $f\colon M\to N$ such that $F=f/s$ for some non-zero $s\in R$. Since $R$ is of finite character, $s$ is contained only in finitely many maximal ideals of $R$. For any other maximal ideal $\fm$, $s/1$ is a unit in $R_\fm$. Therefore, $f_\fm$ is a monomorphism (resp. an epimorphism) for every maximal ideal not containing $s$.

$(ii)$. If $F$ is a splitting monomorphism, then there is a $\Lambda_Q$-module homomorphism $G\colon N_Q\to M_Q$ such that $G\circ F=\mathrm{Id}_{M_Q}$. In particular, $G$ is a splitting epimorphism. By Lemma~\ref{isofg}, there are $\Lambda$-module homomorphisms $f\colon M\to N$ and $g\colon N\to M$ such that $F=f/s$, $G=g/s$ and $g\circ f=s^2\mathrm{Id}_M$ for some non-zero $s\in R$. For any maximal ideal $\fm$ not containing $s$, $s/1$ is a unit in $R_\fm$. Therefore, $f_\fm$ is a splitting monomorphism and $g_\fm$ is a splitting epimorphism for every maximal ideal not containing $s$. The proof in the case where $F$ is a splitting epimorphism follows by changing the roles of $F$ and $G$.

The proof of $(iii)$ follows from $(i)$ and $(ii)$ since $F$ is an isomorphism if and only if it is a monomorphism and a splitting epimorphism.
\end{Proof}

\begin{Cor}
Let $R$ be a commutative domain of finite character with field of fractions $Q$. Let $\Lambda$ be an $R$-algebra, and let $N$ be a finitely generated right $\Lambda$-module which is torsion-free as an $R$-module. If $N_Q$ is projective, then $N_\fm$ is projective for almost all maximal ideals $\fm$ of $R$.
\end{Cor}
\begin{Proof}
Suppose $N_Q$ is projective and let $\tilde f\colon F_Q\to N_Q$ be a splitting epimorphism for some finitely generated free $\Lambda$-module $F$. By Lemma~\ref{almostallmaximals}, there is a $\Lambda$-module homomorphism $f\colon F\to N$ such that the induced homomorphism $f_\fm\colon F_\fm\to N_\fm$ is a splitting epimorphism for almost all maximal ideals $\fm$ of $R$. Thus, $N_\fm$ is projective for almost all maximal ideals $\fm$ of $R$.
\end{Proof}

\begin{Lemma}\label{finitesetofmaximals} 
Let $R$ be a commutative domain, and let $\Lambda$ be a module-finite $R$-algebra. Let $M,N$ be finitely generated right $\Lambda$-modules which are torsion-free as $R$-modules. If $N_\fm$ is a direct summand of $M_\fm$ for each $\fm$ in some finite subset $\mathcal{M}\subseteq\mSpec R$, then there is a $\Lambda$-module homomorphism $f\colon M\to N$ such that the induced homomorphism $f_\fm\colon M_\fm\to N_\fm$ is a splitting epimorphism for every $\fm\in\mathcal{M}$. 
\end{Lemma}
\begin{Proof}
Let $\mathcal{M}=\{\fm_1,\dots,\fm_k\}$ and fix $i\in \{1,\dots ,k\}$. Since $N_{\fm_i}$ is a direct summand of $M_{\fm_i}$, there are $\Lambda_{\fm_i}$-module homomorphisms $\tilde f_i\colon  M_{\fm_i}\to N_{\fm_i}$ and $\tilde g_i\colon  N_{\fm_i}\to M_{\fm_i}$ such that ${\tilde{f_i}\circ \tilde{g_i}=\mathrm{Id}_{N_{\fm_i}}}$. By Lemma~\ref{isofg}, there are $\Lambda$-module homomorphisms $f_i\colon M\to N$ and $g_i\colon N\to M$ such that $\tilde f_i=f_i/s_i$, $\tilde g_i=g_i/s_i$ and $f_i\circ g_i=s_i^2\mathrm{Id}_N$ for some $s_i\notin\fm_i$. Note that, for any $r\notin\fm_i$, $(rf_i)_{\fm_i}$ is a splitting epimorphism. 

For any $i\in \{1,\dots ,k\} $, because  $R=\fm _i+ \bigcap _{j\neq i} \fm _j$, $1=s_i+r_i$ with $s_i\in \fm _i$ and $r_i \in \bigcap _{j\neq i} \fm _j$. 
Let $f=\sum_{i=1}^kr_if_i$. We claim that $f_{\fm_i}$ is a splitting epimorphism. 

Since  $(r_if_i)_{\fm_i}$ is a splitting epimorphism, there is a $\Lambda_{\fm_i}$-module homomorphism $h\colon N_{\fm_i}\to M_{\fm_i}$ such that $(r_if_i)_{\fm_i}\circ h= \mathrm{Id}_{N_{\fm_i}}$. Hence $f_{\fm_i}\circ h=\mathrm{Id}_{N_{\fm_i}}+h'$, where $h'\in \End_{\Lambda_{\fm_i}}(N_{\fm_i})\fm_i R_{\fm_i}$. By Lemma~\ref{jacobson}, $h'\in J(\End_{\Lambda_{\fm_i}}(N_{\fm_i}))$. Therefore, $f_{\fm_i}\circ h$ is invertible and  $f_{\fm_i}$ is a splitting epimorphism, as claimed.
\end{Proof}

\begin{Cor}\label{semilocalsummands}
Let $R$ be a commutative semilocal domain, and let $\Lambda$ be a module-finite $R$-algebra. Let $M,N$ be finitely generated right $\Lambda$-modules which are torsion-free as $R$-modules. Then
\begin{enumerate}
    \item[(i)] $N$ is a direct summand of $M$ if  and only if $N_\fm$ is a direct summand of $M_\fm$ for all maximal ideals $\fm$ of $R$.
    \item[(ii)] $N$ is isomorphic to $M$ if  and only if $N_\fm$ is isomorphic to $M_\fm$ for all maximal ideals $\fm$ of $R$.
\end{enumerate}
\end{Cor}
\begin{Proof} 
$(i)$ follows immediately from Lemma~\ref{finitesetofmaximals}. To prove $(ii)$, observe that by $(i)$ there exist $K$ and $K'$ finitely generated $\Lambda$-modules, such that $M\cong N\oplus K$ and $N\cong M\oplus K'$. Since $N$ and $M$ are torsion-free $R$-modules of the same rank, we deduce that $K=K'=0$.
\end{Proof}

The following is an extension of \cite[Lemma 6.1]{GL} to the case of torsion-free modules over  domains of finite character. 

\begin{Prop}\label{directNK}
Let $R$ be a commutative domain of finite character with field of fractions $Q$. Let $\Lambda$ be a module-finite $R$-algebra, and let $M,N,K$ be finitely generated right $\Lambda$-modules which are torsion-free as $R$-modules. Assume that:
\begin{enumerate}
    \item[(i)] $M_\fm$ is a direct summand of $N_\fm$ for all maximal ideals $\fm$ of $R$, and
    \item[(ii)] $M_Q$ is a direct summand of $K_Q$.
\end{enumerate}
Then $M$ is a direct summand of $N\oplus K$.
\end{Prop}

\begin{Proof} 
It suffices to find $\Lambda$-module homomorphisms $f\colon N\to M$, $g\colon K\to M$ such that, for any $\fm\in \mSpec R$, either $f_\fm$ or $g_\fm$ is a splitting epimorphism. For then $(f,g)\colon N\oplus K\to M$ is locally a splitting epimorphism and, by Lemma~\ref{locallysplit}, is a splitting epimorphism.

By (ii), since $M_Q$ is a direct summand of $K_Q$, there is a splitting epimorphism $\tilde g\colon K_Q\to M_Q$. By Lemma~\ref{almostallmaximals}, there is a $\Lambda$-module homomorphism $g\colon K\to M$ such that the induced homomorphism $g_\fm\colon K_\fm\to M_\fm$ is a splitting epimorphism for almost all maximal ideals $\fm$ of $R$.

By (i), since $M_\fm$ is a direct summand of $N_\fm$, there is a splitting epimorphism $\tilde f_\fm\colon N_\fm\to M_\fm$ for all maximal ideals $\fm$ of $R$. Let $\mathcal{M}\subseteq\mSpec R$ be the finite set of maximal ideals such that $g_\fm$ is not a splitting epimorphism. By Lemma~\ref{finitesetofmaximals}, there is a $\Lambda$-module homomorphism $f\colon N\to M$ such that the induced homomorphism $f_\fm:N_\fm\to M_\fm$ is a splitting epimorphism for every $\fm\in\mathcal{M}$. This finishes the proof of the proposition.
\end{Proof}

\section{Package deal theorems for localizations over \texorpdfstring{$h$}{h}-local domains}\label{s:deal}

In this section we will develop tools to deal with modules over an algebra $\Lambda$ over an $h$-local domain  $R$. A particular instance of this situation is when  $\Lambda$ is the endomorphism ring of a finitely generated, torsion-free, $R$-module (recall Lemma~\ref{endfg}). We are interested in knowing what modules over such endomorphism rings are.  

In \S~\ref{PDmodules}  we prove basic properties of the localization at maximal ideals of $R$ of  finitely  generated   modules  over  an algebra $\Lambda$ over a domain of finite character.

In \S~\ref{PDresults} we enter the study of constructing $\Lambda$-modules having prescribed localization at the maximal ideals of $R$, and proving the so-called \emph{Package Deal Theorems}. To do that, we extend the methods of Levy and Odenthal from algebras over noetherian rings of Krull dimension $1$ to algebras over $h$-local domains.

Finally, in \S~\ref{PDtraces} we prove a package deal theorem for traces of countably generated projective modules. Trace ideals are not, in general, finitely generated but, as we will see, they satisfy enough finiteness conditions to be able to extend our methods to this case. 

\subsection{Finitely generated modules over domains of finite character} \label{PDmodules}

\begin{Lemma}\label{almostQ}
Let $R$ be a commutative domain of finite character with field of fractions $Q$. Let $\Lambda$ be an $R$-algebra, and let $M$ be a finitely generated right $\Lambda$-module, which is torsion-free as an $R$-module. Let $N$ be a $\Lambda$-submodule of $M$ such that $M_Q=N_Q$. Then $M_\fm=N_\fm$ for almost all maximal ideals $\fm$ of $R$.

In particular, $N_\fm$ is finitely generated as $\Lambda _\fm$-module for almost all maximal ideals $\fm$ of $R$.
\end{Lemma}
\begin{Proof}
Since $M$ is torsion-free as an $R$-module, it can be seen as a $\Lambda$-submodule of $M_Q=N_Q$. Therefore, $M\subseteq N_Q$, and since $M$ is finitely generated, by Lemma~\ref{d}, $dM\subseteq N$ for some non-zero $d\in R$. Since $R$ is of finite character and $d$ is non-zero, $d$ is contained only in finitely many maximal ideals of $R$. For any other maximal ideal $\fm$, $d/1$ is a unit in $R_\fm$, hence, $M_\fm=N_\fm$, and $N_\fm$ is finitely generated.
\end{Proof}

\begin{Lemma}\label{almostfg}
Let $R$ be a commutative domain of finite character with field of fractions $Q$. Let $Y$ be a submodule of a finitely generated torsion-free $R$-module $M$. Then there exist $y_1,\dotsc,y_r\in Y$ such that $Y_\fm$ is a finitely generated free $\Lambda _\fm$-module with basis $\frac{y_1}{1},\dotsc,\frac{y_r}{1}$ for almost all maximal ideals $\fm$ of $R$.
\end{Lemma}
\begin{Proof} 
By Remark~\ref{rem-almost-free}, we may assume that $M\le R^s$ for $s=\rank M$, and then $Y$ is a submodule of $R^s$.

Let $r=\rank Y$, and notice that $r\le s$. Let $y_1,\dotsc,y_r\in Y$ such that they are $R$-linearly independent (that is, they form a basis of $Y_Q$). Then, there is an embedding $f\colon R^r\to R^s$ given by $f(a_1,\dots a_r)=\sum _{i=1}^r y_ia_i$. Since $f_Q$ is a splitting monomorphism, by Lemma~\ref{almostallmaximals}, $f_\fm$ is a splitting monomorphism for almost all maximal ideals of $R$, which implies that $Z_\fm= \mathrm {Im} f_\fm =\sum _{i=1}^r \frac{y_i}1 R_\fm$ is a direct summand of $R_\fm^s$ for almost all maximal ideals $\fm$ of $R$. Since $Z_\fm\subseteq Y_\fm$ and they have the same rank, we deduce, from the modular law, that $Z_\fm = Y_\fm$ for almost all maximal ideals $\fm$. Therefore,  $Y_\fm$ has $\frac{y_1}1,\dots , \frac{y_r}1$ as a basis for almost all maximal ideals $\fm$ of $R$.
\end{Proof}

\begin{Cor}\label{almost-free}
Let $R$ be a commutative domain of finite character with field of fractions $Q$. Let $M$ be a finitely generated $R$-module. Then, 
\begin{enumerate}
    \item[(i)] If $M$ is torsion-free then $M_\fm$ is a free $R_\fm$-module for almost all maximal ideals $\fm$ of $R$;
    \item[(ii)] $M_\fm$ is a finitely presented $R_\fm$-module module of projective dimension at most one for almost all maximal ideals $\fm$ of $R$.
\end{enumerate}

\end{Cor}
\begin{Proof}
$(i)$. Since $M$ is torsion-free as an $R$-module, it can be seen as an essential submodule of its injective hull $Q^n$. By Remark~\ref{rem-almost-free}, $dR^n\subseteq M$ for some non-zero $d\in R$. By Lemma~\ref{almostQ} with $\Lambda=R$ and $N=dR^n$, $M_\fm=R_\fm^n$ for almost all maximal ideals $\fm$ of $R$, as claimed.

$(ii)$. Assume now that $M$ is a finitely generated $R$-module. Consider a presentation of $M$
    \[\begin{tikzcd}
        0\rar & N\rar & R^n \rar & M\rar & 0.
    \end{tikzcd}\]
By Lemma~\ref{almostfg}, $N_\fm$ is a finitely generated free module for almost all maximal ideals $\fm$ of $R$. This proves the claim.
\end{Proof}

\begin{Cor}\label{almost-free-count-ds}
Let $R$ be a commutative domain of finite character. Let $M$ be a countable direct sum of finitely generated torsion-free $R$-modules. Then $M_\fm$ is a free $R_\fm$-module for all but countably many maximal ideals $\fm$ of $R$.
\end{Cor}
\begin{Proof}
Let $M=\bigoplus_{i\in\N} M_i$, where each $M_i$ is a finitely generated torsion-free module. By Corollary~\ref{almost-free}, $(M_i)_\fm$ is free for almost all maximal ideals $\fm$ of $R$. Therefore, $M_\fm$ is free for all but countably many maximal ideals $\fm$ of $R$.
\end{Proof}

\subsection{Gluing localizations over \texorpdfstring{$h$}{h}-local domains} \label{PDresults}

It is very important to keep in mind the following structure result of torsion modules over $h$-local domains. It was already proved by Matlis \cite[Theorem~22]{matlis}, but since the result and its proof are quite important in our work, we also include a proof of it.

\begin{Lemma}\label{criteriafghlocaltorsion} 
Let $R$ be an $h$-local domain, and let $\Lambda$ be an $R$-algebra. Let $M$ be a $\Lambda$-module that is torsion as an $R$-module. Then $M\cong \bigoplus _{\fm\in \mSpec (R)} M_\fm $, each $M_\fm$ is a homomorphic image of $M$ and its structure as $\Lambda _\fm$-module is the same as the structure as $\Lambda$-module. 

Moreover, $M$ is finitely generated as $\Lambda$-module if and only if there exists a finite set $\mathcal{S}$ of maximal ideals of $R$ such that   $M_\fm =\{0\}$ for any $\fm \not \in \mathcal{S}$ and $M_\fm$ is a finitely generated $\Lambda _\fm$-module for any $\fm\in \mathcal{S}$    
\end{Lemma} 
\begin{Proof} 
For any maximal ideal $\fm$ of $R$,  let  $\lambda _\fm \colon M\to M_\fm$ denote the localization map. Then there is  an inclusion $\lambda\colon M\hookrightarrow \prod _{\fm \in \mSpec(R)} M_\fm$ defined by $m\mapsto (\lambda _\fm (m))_{\fm \in \mSpec(R)}$ for any $m\in M$.

For any $m\in M$, $m\Lambda$ is an $R/\mathrm{ann}_R(m)$-module. Being $M$ torsion as an $R$-module, $\mathrm{ann}_R(m)\neq \{0\}$. Since $R$ is $h$-local, $\mathrm{ann}_R(m)$ is contained only in finitely many maximal ideals of $R$. Therefore, $m\Lambda _\fm=\{0\}$ for almost all maximal ideals $\fm$ of $R$. Hence,  $\lambda$ has its image in $\bigoplus _{\fm \in \mSpec(R)} M_\fm$. 

Let $\mathcal{S}$ denote the finite set of maximal ideals containing $\mathrm{ann}_R(m)$. By Lemma~\ref{hlocal}, $m\Lambda\cong \bigoplus _{\fm \in \mathcal{S}}\left(m\Lambda\right)_\fm$,  and also $\Lambda /\mathrm{ann}_R(m)\cong \prod _{\fm \in \mathcal{S}}\left(\Lambda /\mathrm{ann}_R(m)\right)_\fm $. 

This allows us to conclude that, for any $m\in M$ and $\fm \in \mSpec (R)$, there exists $m'\in M$ such that $\lambda _\fn (m')=0$ for any $\fn \neq \fm$ and such that $\lambda _\fm (m')=\lambda _\fm (m)$. This implies that $\lambda$ induces an isomorphism $M\cong \bigoplus _{\fm \in \mSpec(R)} M_\fm$.

In particular, for any maximal ideal $\fm$ of $R$, $M_\fm$ is a homomorphic image of $M$  and the structure as $\Lambda$-module is the same as the structure as $\Lambda _\fm$-module. So the statement about the finite generation of $M$ easily follows from our previous discussion.
\end{Proof}

Next lemma will be used to determine when  we get finitely generated modules in our package deal results.

\begin{Lemma}\label{lem:fg}
Let $R$ be an $h$-local domain with field of fractions $Q$. Let $\Lambda$ be an $R$-algebra, and let $M$ be a finitely generated right $\Lambda$-module, which is torsion-free as an $R$-module. Let $N$ be a $\Lambda$-submodule of $M$ such that $M_Q=N_Q$. Then $N$ is finitely generated as a $\Lambda$-module if and only if $N_\fm$ is finitely generated as a $\Lambda_\fm$-module for every maximal ideal $\fm$ of $R$.
\end{Lemma}
\begin{Proof}
It is clear that if $N$ is finitely generated as a $\Lambda$-module, then $N_\fm$ is finitely generated as a $\Lambda_\fm$-module for every maximal ideal $\fm$ of $R$.

To prove the converse implication, assume that $N_\fm$ is finitely generated as a $\Lambda_\fm$-module for every maximal ideal $\fm$ of $R$.
Since $M$ is torsion-free as an $R$-module, it can be seen as a $\Lambda$-submodule of $M_Q=N_Q$. Therefore, $M\subseteq N_Q$, and since $M$ is finitely generated as a $\Lambda$-module, by Lemma~\ref{d}, $dM\subseteq N$ for some non-zero $d\in R$. Note that $dM$ is also finitely generated. So we only need to prove that $N/dM$ is a finitely generated $\Lambda/d\Lambda$-module.

Since $R$ is $h$-local and $d$ is non-zero, $d$ is contained only in finitely many maximal ideals of $R$, say $\mathcal M=\{\fn_1,\dotsc,\fn_t\}$. By Lemma~\ref{hlocal}, there is an isomorphism
    \[\varphi\colon R/dR\to (R/dR)_{\fn_1}\times\dotsb\times (R/dR)_{\fn_t}.\]
Therefore, $N/dM\cong N_{\fn_1}/dM_{\fn_1}\oplus\dotsb\oplus N_{\fn_t}/dM_{\fn_t}$. By Lemma~\ref{criteriafghlocaltorsion}, $N/dM$ is finitely generated as a $\Lambda/d\Lambda$-module. Therefore, $N$ is finitely generated  as a $\Lambda$-module.
\end{Proof}

\begin{PDTh}[Localization of Submodules] \label{dealsubmodules}
Let $R$ be an $h$-local domain with field of fractions $Q$. Let $\Lambda$ be an $R$-algebra, and let $M$ be a finitely generated right $\Lambda$-module, which is torsion-free as an $R$-module. For each maximal ideal $\fm$ of $R$, let $X(\fm)$ be a $\Lambda_\fm$-submodule of $M_\fm$, which is torsion-free as an $R_\fm$-module, and such that $X(\fm)_Q=(M_\fm)_Q$. Then the following statements are equivalent
\begin{enumerate}
    \item[(i)] There is a $\Lambda$-submodule $N\subseteq M$, which is torsion-free as an $R$-module, and such that $N_\fm=X(\fm)$ for all maximal ideals $\fm$ of $R$.
    \item[(ii)] $X(\fm)=M_\fm$ for almost all maximal ideals $\fm$ of $R$.
\end{enumerate}
Moreover, if each of the $X(\fm)$ is finitely generated as a $\Lambda_\fm$-module, then $N$ is also finitely generated as a $\Lambda$-module.
\end{PDTh}
\begin{Proof}
$(ii)\Rightarrow(i)$. Suppose that $X(\fm)=M_\fm$ for almost all maximal ideals $\fm$ of $R$. Let $\mathcal{M}=\{\fn_1,\dotsc,\fn_t\}$ denote the finite set of maximal ideals of $R$ such that $X(\fm)\neq M_\fm$.

Fix $i\in\{1,\dotsc,t\}$. Since $M_{\fn_i}$ is torsion-free as an $R_{\fn_i}$-module, it can be seen as a $\Lambda_{\fn_i}$-submodule of $(M_{\fn_i})_Q=X(\fn_i)_Q$. Therefore, $M_{\fn_i}\subseteq X(\fn_i)_Q$, and since $M_{\fn_i}$ is finitely generated, by Lemma~\ref{d}, $d_iM_{\fn_i}\subseteq X(\fn_i)$ for some non-zero $d_i\in R$. Since $X(\fn_i)\neq M_{\fn_i}$, $d_i\in\fn_i$ (otherwise $d_i/1$ would be a unit in $R_{\fn_i}$ and, hence $X(\fn_i)=M_{\fn_i}$, a contradiction). 

Let $d=d_1\dotsb d_t\in\fn_1\cap\dots\cap\fn_t$. Since $R$ is $h$-local and $d$ is non-zero, $d$ is contained only in finitely many maximal ideals of $R$, say $\{\fm_1=\fn_1,\dots,\fm_t=\fn_t,\dotsc,\fm_k\}$. By Lemma~\ref{hlocal}, the canonical homomorphism   
    \[\varphi:R/dR\to(R/dR)_{\fm_1}\times\dotsb\times(R/dR)_{\fm_k}\]
is an isomorphism. Therefore, there are non-zero $b_1,\dotsc,b_t,b\in R$ such that $\varphi(\overline b_i)=(\bar0,\dotsc,\bar1^{(i)},\dotsc,\bar0)$ and $\varphi(\overline b)=(\bar0,\overset{(t}{\dotsc},\bar0,\bar1,\dots,\bar1)$. Note that $b_i/1$ is a unit in $R_{\fn_i}$ for every $i=1,\dotsc,t$, $b/1,b_j/1\in dR_{\fn_i}$ for every $j\neq i$, and $b/1$ is a unit in $R_{\fm_i}$ for every $i=t+1,\dotsc,k$.

For each $i\in\{1,\dotsc,t\}$, let $X_i'$ be the $\Lambda$-submodule of $M$ generated by the numerators of some set of $\Lambda_{\fn_i}$-generators of $X(\fn_i)$. Let $N=b_1X_1'+\dotsb+b_tX_t'+bM+dM$, which is a $\Lambda$-submodule of $M$. Since $dM_{\fn_i}\subseteq X(\fn_i)$, $N_{\fn_i}=X(\fn_i)+dM_{\fn_i}=X(\fn_i)$ for every $i=1,\dotsc,t$. On the other hand, $N_{\fm_i}=M_{\fm_i}+dM_{\fm_i}=M_{\fm_i}$ for every $i=t+1,\dotsc,k$. For any other maximal ideal $\fm$, $d/1$ is a unit in $R_\fm$, hence $N_\fm=M_\fm$, as claimed.

$(i)\Rightarrow(ii)$. Conversely, let $N$ be a $\Lambda$-submodule of $M\subseteq M_Q$, which is torsion-free as an $R$-module, and such that $N_\fm=X(\fm)$ for all maximal ideals $\fm$ of $R$. Since $M_Q=X(\fm)_Q=(N_\fm)_Q$, by Lemma~\ref{almostQ}, $M_\fm=N_\fm=X(\fm)$ for almost all maximal ideals $\fm$ of $R$. 

Finally, the last part of the statement follows from Lemma~\ref{lem:fg}. This finishes the proof of the theorem.
\end{Proof}

\begin{Cor} \label{existence}
Let $R$ be an $h$-local domain with field of fractions $Q$. For each maximal ideal $\fm$ of $R$, let $X(\fm)$ be a finitely generated torsion-free $R_\fm$-module of rank $n$. Then the following statements are equivalent
\begin{enumerate}
    \item[(i)] There is a finitely generated torsion-free $R$-module $N$ of rank $n$, such that $N_\fm\cong X(\fm)$ for all maximal ideals $\fm$ of $R$.
    \item[(ii)] $X(\fm)$ is a free $R_\fm$-module of rank $n$ for almost all maximal ideals $\fm$ of $R$.
\end{enumerate}
\end{Cor}
\begin{Proof}
$(ii)\Rightarrow(i)$. Suppose that $X(\fm)$ is a free $R_\fm$-module of rank $n$ for almost all maximal ideals $\fm$ of $R$. Let $\mathcal{M}=\{\fm_1,\dotsc,\fm_t\}$ be the finite set of maximal ideals of $R$ such that $X(\fm)$ is not free. 

Fix $i\in\{1,\dotsc,t\}$. Since $X(\fm_i)$ is finitely generated and torsion-free, we can see $X(\fm_i)$ as an essential submodule of its injective hull $Q^n$. On the other hand, $R_{\fm_i}^n$ has the same injective hull $Q^n$. Since $X(\fm_i)$ is finitely generated, by Lemma~\ref{d}, $d_iX(\fm_i)\subseteq R_{\fm_i}^n$ for some non-zero $d_i\in R$. 

Let $d=d_1\dots d_t\in R$ and let $M=d^{-1}R^n$. Note that $X(\fm)$ is an $R_\fm$-submodule of $M_\fm$. Since $R$ is $h$-local and $d$ is non-zero, $d$ is contained only in finitely many maximal ideals of $R$. For any other maximal ideal $\fm$, $d/1$ is a unit in $R_\fm$. Therefore, $M_\fm=R_\fm^n=X(\fm)$ for almost all maximal ideals $\fm$ of $R$. By Package Deal Theorem~\ref{dealsubmodules}, there is a finitely generated torsion-free $R$-module $N$, such that $N_\fm=X(\fm)$ for all maximal ideals $\fm$ of $R$.

$(i)\Rightarrow(ii)$. The converse is just Corollary~\ref{almost-free}. 
\end{Proof}

In \S~\ref{s:dsummands} we will prove some results that allow us to make the following considerations about the (non)uniqueness of the modules constructed in Theorem~\ref{dealsubmodules} and Corollary~\ref{existence}.

\begin{Rem}
It is well known that two finitely generated modules $M$ and $M'$ over a commutative domain $R$ that are  in the same genus (that is, with isomorphic localizations at maximal ideals of $R$) need not be isomorphic. 

If $M$ and $M'$ are also torsion-free, we will see in Proposition~\ref{summandk} that  $\add(M)=\add(M')$. If $R$ has finite character then, by Proposition~\ref{directNK},  $M$ is a direct summand of $M'\oplus M'$ and $M'$ is a direct summand of $M\oplus M$.  If $R$ is semilocal then $M\cong M'$ by Corollary~\ref{semilocalsummands}
\end{Rem}

\subsection{Package deal for traces of projective modules}\label{PDtraces}

Now we want to focus on localizations of trace ideals of countably generated projective modules over suitable algebras over $h$-local domains. Trace ideals of countably generated projective right modules were characterized by Whitehead in \cite{whitehead} and with more detail by Herbera and Příhoda in \cite{traces}. We recall here this characterization.

\begin{Prop}\emph{\cite[Proposition 2.4]{traces}}\label{traces}
Let $R$ be a ring, and let $I$ be a two-sided ideal of $R$. Then $I$ is the trace ideal of a countably generated projective right $R$-module if and only if there exists an ascending chain of finitely generated left ideals $(J_n)_{n\ge1}$ such that $J_{n+1}J_n=J_n$ and $I=\bigcup_{n\ge1} J_nR$.
\end{Prop}

It is useful to keep in mind the following Lemma, as it explains some modifications that can be made in the ascending chain in Proposition~\ref{traces}.

\begin{Lemma}\emph{\cite[Lemma 2.2]{traces}}\label{lem:traces}
Let $R$ be a ring. Let $J_1\subseteq J_2$ be finitely generated left ideals of $R$ satisfying that $J_2J_1=J_1$. For $i=1,2$, fix $A_i$ a finite set of generators of $J_i$. Let $X$ be a finite subset of $R$ such that $1\in X$. For $i=1,2$, set
	\[J_i'=\sum_{x\in X,\,a\in A_i} Rax.\]
Then $J_1'\subseteq J_2'$ and $J_2'J_1'=J_1'$. Moreover, for $i=1,2$, $J_i\subseteq J_iR=J_i'R$.
\end{Lemma}

\begin{Lemma}
\label{almosttrace}
Let $R$ be an $h$-local domain with field of fractions $Q$. Let $\Lambda$ be a torsion-free $R$-algebra such that $\Lambda_Q$ is a simple artinian ring. Let $I$ be a non-zero two-sided ideal of $\Lambda$. Then $I_\fm=\Lambda _\fm$ for almost all maximal ideals $\fm$ of $R$.
\end{Lemma}
\begin{Proof}
$I_Q$ is a non-zero two-sided ideal of $\Lambda_Q$, so $I_Q= \Lambda_Q $. Therefore, there is a non-zero $q\in R$ such that $q=\sum a_i\lambda_i$ for $a_i\in I$ and $\lambda _i \in \Lambda$. Since $R$ is $h$-local,  $q$ is invertible in almost all maximal ideals of $R$. The claim follows. 
\end{Proof}

\begin{PDTh}[Localization of Trace Ideals] \label{dealtraces}
Let $R$ be an $h$-local domain with field of fractions $Q$. Let $\Lambda$ be a torsion-free $R$-algebra such that $\Lambda_Q$ is a simple artinian ring. For each maximal ideal $\fm$ of $R$, let $I(\fm)$ be a non-zero two-sided ideal of $\Lambda_\fm$ which is the trace ideal of a countably generated projective right $\Lambda_\fm$-module. Then the following statements are equivalent
\begin{enumerate}
    \item[(i)] There is a two-sided ideal $I$ of $\Lambda$, which is the trace ideal of a countably generated projective right $\Lambda$-module, and such that $I_\fm=I(\fm)$ for all maximal ideals $\fm$ of $R$.
    \item[(ii)] $I(\fm)=\Lambda_\fm$ for almost all maximal ideals $\fm$ of $R$.
\end{enumerate}
Moreover, the ideal $I$ that satisfies the equivalent conditions $(i)$ and $(ii)$ is unique.
\end{PDTh}
\begin{Proof} By \cite[Lemmas~10.3 and 10.4]{wiegand} two trace ideals with isomorphic localizations at maximal ideals of $R$ are equal. This implies that the ideal $I$ in statement $(i)$ is unique. Now we proceed to prove the equivalence of the two statements.

$(ii)\Rightarrow(i)$. Let $\mathcal{M}=\{\fn_1,\dotsc,\fn_t\}$ and fix some $i\in\{1,\dotsc,t\}$. We may assume that $I(\fn_i)\neq\Lambda_{\fn_i}$. Since $I(\fn_i)$ is the trace ideal of a countably generated projective right $\Lambda_{\fn_i}$-module, by Proposition~\ref{traces}, there is an ascending chain of non-zero finitely generated left ideals $(J_{i,n})_{n\ge1}$ of $\Lambda_{\fn_i}$ such that $J_{i,n+1}J_{i,n}=J_{i,n}$ and $I(\fn_i)=\bigcup_{n\ge 1} J_{i,n}\Lambda_{\fn_i}$. 

Fix $n\ge1$. Let $A_{i,n}\subseteq\Lambda$ be a finite set of $\Lambda_{\fn_i}$-generators of $J_{i,n}$. Since $J_{i,n}$ is non-zero and $\Lambda_Q$ is simple artinian, $(J_{i,n})_Q\Lambda_Q=\Lambda_Q$. Since $\Lambda$ is torsion-free, it can be seen as a $\Lambda$-submodule of $\Lambda_Q=(J_{i,n})_Q\Lambda_Q$. Therefore, $\Lambda\subseteq (J_{i,n})_Q\Lambda_Q$, and by Lemma~\ref{d}, $d_{i,n}\Lambda_{\fn_i}\subseteq J_{i,n}\Lambda_{\fn_i}$ for some non-zero $d_{i,n}\in\fn_i$. Then,
	\[\frac{1}{d_{i,n}} \sum_{j=1}^{m_{i,n}}\sum_{a_{i,n}\in A_{i,n}} \frac{r_{j,i,n}}{1}\frac{a_{i,n}}{1}\frac{s_{j,i,n}}{1}=1_{\Lambda_Q}\in\Lambda_Q,\]
for some $r_{j,i,n},s_{j,i,n}\in\Lambda$ for every $j=1,\dotsc,m_{i,n}$.

Let $X_n=\{1_\Lambda\}\cup\bigcup_{i=1}^t\{s_{1,i,n},s_{2,i,n},\dotsc,s_{m_{i,n},i,n}\}\subseteq\Lambda$ and define 
	\[J_{i,n}'=\sum_{x_n\in X_n,\,a_{i,n}\in A_{i,n}} \Lambda_{\fn_i}a_{i,n}x_n.\]
By Lemma~\ref{lem:traces}, $J_{i,n}'\subseteq J_{i,n+1}'$, $J_{i,n+1}'J_{i,n}'=J_{i,n}'$, and $J_{i,n}\subseteq J_{i,n}\Lambda_{\fn_i}=J_{i,n}'\Lambda_{\fn_i}$. Let $d_n=d_{1,n}\dots d_{t,n}\in\fn_1\cap\dots\cap\fn_t$. Note that $d_n\Lambda_{\fn_i}\subseteq J_{i,n}'$ by construction. Define
	\[L_{i,n}=\sum_{x_n\in X_n,\,a_{i,n}\in A_{i,n}} \Lambda a_{i,n}x_n+\Lambda d_n,\]
which are left ideals of $\Lambda$. Note that $(L_{i,n})_{\fn_i}=J_{i,n}'$ because $d_n\Lambda_{\fn_i}\subseteq J_{i,n}'$. Since $R$ is of finite character and $d_n$ is non-zero, $d_n$ is contained only in finitely many maximal ideals of $R$. For any other maximal ideal $\fm$, $d_n/1$ is a unit in $R_\fm$, hence $(L_{i,n})_\fm=\Lambda_\fm$. Let $\mathcal{M}_n=\{\fm_1=\fn_1,\dotsc,\fm_t=\fn_t,\fm_{t+1},\dotsc,\fm_{k_n}\}$ be the finite set of maximal ideals containing $d_n$. By Lemma~\ref{hlocal}, the canonical homomorphism  
    \[\varphi_n:R/d_nR\longrightarrow(R/d_nR)_{\fm_1}\times\dotsb\times(R/d_nR)_{\fm_{k_n}}\]
is an isomorphism. Therefore, there are non-zero elements $b_{1,n},\dotsc,b_{t,n},b_n\in R$ such that $\varphi_n(\bar b_{i,n})=(\bar0,\dotsc,\bar1^{(i)},\dotsc,\bar0)$ and $\varphi_n(\bar b_n)=(\bar0,\overset{t)}{\dotsc},\bar0,\bar1,\dotsc,\bar1)$. 

Let $J_n=L_{1,n}b_{1,n}+\dots+L_{t,n}b_{t,n}+\Lambda b_n+\Lambda d_n$. Note that $(J_n)_{\fn_i}=J_{i,n}'$, and $(J_n)_\fm=\Lambda_\fm$ for every $\fm\neq\fn_i$. Consider the short exact sequence of left $\Lambda$-modules
    \[\begin{tikzcd}
        0\rar & J_{n+1}\rar & J_{n+1}+J_n\rar & (J_{n+1}+J_n)/J_{n+1}\rar & 0.
    \end{tikzcd}\]
Since localization is an exact functor,
    \[\begin{tikzcd}
        0\rar & (J_{n+1})_\fm\rar & (J_{n+1}+J_n)_\fm\rar & (J_{n+1}+J_n)_\fm/(J_{n+1})_\fm\rar & 0
    \end{tikzcd}\]
is also a short exact sequence of left $\Lambda_\fm$-modules. Note that, since $J_{i,n}'\subseteq J_{i,n+1}'$, $(J_{n+1}+J_n)_{\fn_i}/(J_{n+1})_{\fn_i}=(J_{i,n+1}'+J_{i,n}')/J_{i,n+1}'=0$ for every $i=1,\dotsc,t$. On the other hand, $(J_{n+1}+J_n)_\fm/(J_{n+1})_\fm=(\Lambda_\fm+\Lambda_\fm)/\Lambda_\fm=0$ for every $\fm\neq\fn_i$. Therefore, $(J_{n+1}+J_n)/J_{n+1}=0$ and we deduce that $(J_n)_{n\ge1}$ is an increasing sequence of left ideals of $\Lambda$. 

Now, $(J_{n+1})_{\fn_i}(J_n)_{\fn_i}=J_{i,n+1}'J_{i,n}'=J_{i,n}'=(J_n)_{\fn_i}$ for every $i=1,\dotsc,t$, and $(J_{n+1})_\fm(J_n)_\fm=\Lambda_\fm^2=\Lambda_\fm=(J_n)_\fm$ for every $\fm\neq\fn_i$. Therefore, $J_{n+1}J_n=J_n$ for every $n\ge1$, and by Proposition~\ref{traces}, $I=\bigcup_{n\ge1} J_n\Lambda$ is a two-sided ideal of $\Lambda$ such that it is the trace ideal of a countably generated projective right $\Lambda$-module, $I_{\fn_i}=I(\fn_i)$ and $I_\fm=\Lambda_\fm$ otherwise.

$(i)\Rightarrow(ii)$. The converse is just Lemma~\ref{almosttrace}.
\end{Proof}

\section{The global case}\label{s:global}

In this section we study when the class $\mathcal{F}$, of direct sums of finitely generated torsion-free modules over an $h$-local domain $R$, is closed under direct summands. In the next Proposition we prove that this property \emph{localizes} at maximal ideals of $R$, so the results from Section~\ref{s:local} hold for $R_\mathfrak{m}$ for any maximal ideal $\mathfrak{m}$ of $R$.

It is worth mentioning that we do not know how to prove directly that if $\mathcal{F}$ is closed under direct summands then the same holds for $R_\mathfrak{m}$ for any maximal ideal $\mathfrak{m}$ of $R$. The property that passes to the localization of $R$ is that, for any finitely generated torsion-free $R$-module, every direct summand of $R^{(\omega)}\oplus X$ is a direct sum of finitely generated modules.

\begin{Prop}\label{propertylocalizes}
Let $R$ be an $h$-local domain of Krull dimension $1$. Let $X$ be a finitely generated torsion-free $R$-module, and assume that every direct summand of $R^{(\omega)}\oplus X$ is a direct sum of finitely generated modules. Then, for any maximal ideal $\fm$ of $R$, $\End_{R_\fm}(X_\fm)$ is a semiperfect ring, so that $X_\fm$ satisfies the equivalent conditions of Proposition~\ref{localAddRM}.
\end{Prop}
\begin{Proof}
Fix a maximal ideal $\fm$ of $R$, and set $M'=X_\fm$. By Lemma~\ref{endfg}(v), $\End_{R_\fm}(M')$ is a semilocal ring. We want to prove that it is semiperfect. We decompose $M'=F\oplus M$, where $F$ is a finitely generated free $R_\fm$-module and $M$ is a finitely generated module with no projective direct summands. Notice that $\End_{R_{\fm}}(M')$ is semiperfect if and only if so is $\End_{R_{\fm}}(M)$.

By Lemma~\ref{localnil}(i) (see also Corollary~\ref{localontobij}), to prove that $\End_{R_{\fm}}(M)$ is semiperfect, it suffices to show that the idempotents of $\End_{R_{\fm}}(M)/M\Hom_{R_{\fm}}(M,R_\fm)$ can be lifted to idempotents of $\End_{R_{\fm}}(M)$. Therefore, we may assume that $M$ is non-zero and that $J=M\Hom_{R_{\fm}}(M,R_\fm)$ is a proper ideal of $\End_{R_{\fm}}(M)$ contained in $J(\End_{R_{\fm}}(M))$.

Let $e\in  \End_{R_{\fm}}(M)/J$ be a non-trivial idempotent. By Lemma~\ref{liftingidempotents} and Proposition~\ref{equivalencia}, to show that $e$ can be lifted to an idempotent of $\End_{R_{\fm}}(M)$ we only need to find a decomposition $M\cong A_1\oplus A_2$ such that, if we set $P_i=\mathrm{Hom}_{R_{\fm}} (M,A_i)$ for $i=1,2$, then $P_1/P_1J\cong e\left(\End_{R_{\fm}}(M)/J\right)$ and $P_2/P_2J\cong (1-e)\left(\End_{R_{\fm}}(M)/J\right)$.

As $M\cong A_1\oplus A_2$, it is easy to see that, for $i=1,2$, 
$P_iJ = \{f\in \mathrm{Hom}_{R_{\fm}} (M,A_i)\mid \mbox{ $f$ factors throught a free module}\} =A_i\mathrm{Hom}_{R_{\fm}} (M,R_{\fm})$. Therefore, we need to prove that  $M\cong A_1\oplus A_2$
with $\mathrm{Hom}_{R_{\fm}} (M,A_1)/A_1\mathrm{Hom}_{R_{\fm}} (M,R_{\fm})\cong e\left(\End_{R_{\fm}}(M)/J\right)$ and $\mathrm{Hom}_{R_{\fm}} (M,A_2)/A_2\mathrm{Hom}_{R_{\fm}} (M,R_{\fm})\cong (1-e)\left(\End_{R_{\fm}}(M)/J\right)$.

By Lemma~\ref{endfg} there is an extension of rings $R/I\subseteq \End_R(X)/X\Hom_R(X,R)$, where $0\neq I=R\cap X\Hom_R(X,R)$. Since $R$ is $h$-local, $I$ is contained only in finitely many maximal ideals $\{ \fm=\fm _1,\dots ,\fm_\ell\}$ of $R$. So  we have an isomorphism 
\[\End_R(X)/X\Hom_R(X,R)\to \prod _{i=1}^\ell \left( \End_R(X)/X\Hom_R(X,R)\right) _{\fm _i}\cong \] \[\cong \prod _{i=1}^\ell \End_{R_{\fm _i}}(X_{\fm _i})/X_{\fm _i}\Hom_{R_{\fm _i}}(X_{\fm _i},R_{\fm _i}) \]
given by taking localization at $\fm _i$ at each component (cf. Lemma~\ref{criteriafghlocaltorsion}).

Hence, there exists $\widetilde{e}^2=\widetilde{e}\in \End_R(X)/X\Hom_R(X,R)$ such that $\lambda _\fn (\widetilde{e})=0$ for any maximal ideal $\fn \neq \fm$ and $\lambda _\fm (\widetilde{e})=e$, where 
    \begin{align*}
        \lambda _\fn \colon \End_R(X)/X\Hom_R(X,R) &\to \left(\End_R(X)/X\Hom_R(X,R) \right )_\fn\cong \\
        &\qquad\cong \End_{R_\fn}(X_\fn)/X_\fn\Hom_{R_\fn}(X_\fn ,R _\fn)
    \end{align*}
denotes the localization morphism.

By Corollary~\ref{stable}, $R^{(\omega)}\oplus X= X_1\oplus X_2$ with 
    \begin{align*}
        \Hom_{R}(R\oplus X, X_1)/X_1\Hom_{R}(R\oplus X,R)&\cong \widetilde{e} \left( \End_R(X)/X\Hom_R(X,R)\right) \\ &\cong e\left( \End_{R_{\fm}}(M)/J\right)
    \end{align*}
and 
    \begin{gather*}
        \Hom_{R}(R\oplus X, X_2)/X_2\Hom_{R}(R\oplus X,R)\cong (1-\widetilde{e}) \left( \End_R(X)/X\Hom_R(X,R)\right)\cong \\
        \cong (1-e)\left( \End_{R_{\fm}}(M)/J\right)\times \left( \prod _{i=2}^\ell \End_{R_{\fm _i}}(X_{\fm _i})/X_{\fm _i}\Hom_{R_{\fm _i}}(X_{\fm _i},R_{\fm _i})\right).
    \end{gather*}
Therefore, 
    $$R_\fm^{(\omega)}\oplus M \cong R_\fm^{(\omega)}\oplus X_\fm= (X_1)_\fm\oplus (X_2)_\fm ,$$ 
$$\Hom_{R_\fm}(R_\fm\oplus X_\fm, (X_1)_\fm)/(X_1)_\fm\Hom_{R_\fm}(R_\fm\oplus X_\fm,R_\fm)\cong  e\left( \End_{R_{\fm}}(M)/J\right)$$
and 
$$\Hom_{R_\fm}(R_\fm\oplus X_\fm, (X_2)_\fm)/(X_2)_\fm\Hom_{R_\fm}(R_\fm\oplus X_\fm,R_\fm)\cong  (1-e)\left( \End_{R_{\fm}}(M)/J\right).$$

By hypothesis, for $i=1,2$, $X_i$ is a direct sum of finitely generated $R$-modules; hence,  $N_i=(X_i)_\fm$ is also a direct sum of finitely generated $R_{\fm}$-modules. By Lemma~\ref{exchange}, we obtain a decomposition  $M=A_1\oplus A_2$ such that, for each $i\in\{1,2\}$, $$\Hom_{R_{\fm}}(R_\fm\oplus M, N_i)/N_i\Hom_{R_{\fm}}(R_\fm\oplus M,R_{\fm})\cong  \Hom_{R_{\fm}}(M,A_i)/A_i\Hom_{R_{\fm}}(M,R_{\fm}).$$

Therefore $\Hom_{R_{\fm}}(M,A_1)/A_1\Hom_{R_{\fm}}(M,R_{\fm})\cong e\left(\End_{R_{\fm}}(M)/J\right)$ and $$\Hom_{R_{\fm}}(M,A_2)/A_2\Hom_{R_{\fm}}(M,R_{\fm})\cong (1-e)\left(\End_{R_{\fm}}(M)/J\right),$$ as we wanted to see. So we can conclude that $e$ can be lifted to an idempotent of $\End_{R_{\fm}}(M)$ by an application of Lemma~\ref{liftingidempotents}.
\end{Proof}

Next result shows that, over an $h$-local domain $R$, to have the class $\mathcal{F}$ closed under direct summands we need   that the ranks of indecomposable modules over different localizations at maximal ideals are coprime. 

This somewhat surprising Theorem~\ref{coprime} is the extension of \cite[Lemma~4]{P3} to our  setting. For its proof, we need our versions of the  Package Deal Theorems~\ref{dealsubmodules} and \ref{dealtraces}, as well as the results from \S~\ref{s:lifting}.

\begin{Th} \label{coprime}
Let $R$ be an $h$-local domain with at least two different maximal ideals $\fm_1,\fm_2$. For each $i=1,2$, let $M_i$ be a finitely generated, indecomposable, torsion-free $R_{\fm_i}$-module with local endomorphism ring and with rank $r_i$. Then, 
\begin{enumerate}
    \item[(1)] for $i=1,2$, there exists a finitely generated, indecomposable, torsion-free $R$-module $X_i$ such that $(X_i)_{\fm_i}\cong M_i$ and $X_\fm\cong R_\fm^{r_i}$ for every maximal ideal $\fm$ different from $\fm_i$; 
    \item[(2)] if $r_1$ and $r_2$ are not coprime then $\Add (X_1\oplus X_2)$ has elements that are not a direct sum of finitely generated modules.
\end{enumerate}
\end{Th}
\begin{Proof}
Corollary~\ref{existence} ensures the existence of finitely generated, indecomposable, torsion-free $R$-modules $X_1$ and $X_2$ with the claimed properties in statement $(1)$.   

To prove $(2)$, let $d=\gcd(r_1,r_2)$. We show that if $d>1$, there exists a module in $\Add(X_1\oplus X_2)$ which is not a direct sum of finitely generated torsion-free $R$-modules.

Let $\Lambda=\End_R(X_1\oplus X_2)$. By Lemma~\ref{endfg}(i) and (iii), $\Lambda$ is a torsion-free $R$-module and $\Lambda_Q$ a simple artinian ring. Recall from Remark~\ref{mplusM} that $\Lambda$ can be identified with the matrix ring:
    \[\Lambda=\begin{pmatrix}
        \End_R(X_1) & \Hom_R(X_2,X_1) \\
        \Hom_R(X_1,X_2) & \End_R(X_2)
    \end{pmatrix}\]
By Lemma~\ref{isofg}, $\Lambda_\fm\cong\End_{R_\fm}((X_1)_\fm\oplus (X_2)_\fm)$ for every maximal ideal $\fm$ of $R$. Let 
    \[I_1=\Lambda_{\fm_1} \begin{pmatrix}1&0\\0&0\end{pmatrix} \Lambda_{\fm_1}\qquad\text{and}\qquad I_2= \Lambda_{\fm_2} \begin{pmatrix}0&0\\0&1\end{pmatrix} \Lambda_{\fm_2}.\]
 By Package Deal Theorem~\ref{dealtraces}, there is a non-zero two-sided ideal $I$ of $\Lambda$,  which is the trace of a countably generated projective right $\Lambda$-module,  such that $I_{\fm_i}=I_i$ ($i=1,2$) and $I_\fm=\Lambda_\fm$ for every maximal ideal $\fm$ of $R$ different from $\fm_i$.

By Lemma~\ref{endfg}(iv), $I\cap R\neq\{0\}$, so $\Lambda/I$ is a torsion $R$-module. Note that the only maximal ideals containing $I\cap R$ are $\fm_1$ and $\fm_2$ (otherwise, if there exists another maximal ideal $\fn$ containing $I\cap R$, $I_\fn\cap R_\fn\neq R_\fn$, and then $I_\fn\neq\Lambda_\fn$, a contradiction). By Lemma~\ref{endfg}(vi), the canonical homomorphism 
    \[\textstyle\Lambda/I\to\bigoplus_{\fm\in\mSpec R}(\Lambda/I)_\fm=(\Lambda/I)_{\fm_1}\times(\Lambda/I)_{\fm_2}\]
is an isomorphism.

From Remark~\ref{mplusM}, 
    \[I_1=\begin{pmatrix} \End_{R_{\fm_1}}(M_1) & \Hom_{R_{\fm_1}}(M_1,R_{\fm_1}^{r_2})\\ \Hom_{R_{\fm_1}}(R_{\fm_1}^{r_2},M_1)& \Hom_{R_{\fm_1}}(M_1,R_{\fm_1}^{r_2})\Hom_{R_{\fm_1}}(R_{\fm_1}^{r_2},M_1) \end{pmatrix}.\]
Therefore, $(\Lambda/I)_{\fm_1}\cong\Lambda_{\fm_1}/I_1\cong M_{r_2}(R_{\fm_1})/J$, where $$J=\Hom_{R_{\fm_1}}(M_1,R_{\fm_1}^{r_2})\Hom_{R_{\fm_1}}(R_{\fm_1}^{r_2},M_1).$$ Notice that $J\neq\End_{R_{\fm_1}}(R_{\fm_1}^{r_2})$ because this would imply that the identity map of $R_{\fm_1}^{r_2}$ is in $J$, so that $R_{\fm_1}^{r_2}$ is a direct summand of $M_1^n$ for some $n\in\N$. Since $M_1$ has local endomorphism ring, the Krull-Schmidt Theorem implies that $M_1\cong R_{\fm_1}$, but we are assuming that $r_1>1$, a contradiction.

Therefore, there is an isomorphism
    \[\textstyle\varphi\colon\Lambda/I\to M_{r_2}(R_{\fm_1})/J\times(\Lambda/I)_{\fm_2}.\]
Then there is an idempotent element $e\in \Lambda/I$ such that $\varphi(e)=(E_{11}+J,0)$, where $E_{11}$ is the idempotent matrix with 1 in the 1,1-entry and zeros elsewhere. Since $M_{r_2}(R_{\fm_1}) E_{11} M_{r_2}(R_{\fm_1})=M_{r_2}(R_{\fm_1})$, $E_{11}\notin J$. 

By Theorem~\ref{liftingproj}(i), there is a countably generated projective right $\Lambda$-module $Q$ such that $Q/QI\cong e(\Lambda/I)$. We claim that such $Q$ is neither finitely generated nor a direct sum of finitely generated modules.

Recall that $\Lambda_{\fm_1}\cong\End_{R_{\fm_1}}(M_1\oplus R_{\fm_1}^{r_2})$. Since $M_1$ has local endomorphism ring, there are only two finitely generated indecomposable projective $\Lambda_{\fm_1}$-modules up to isomorphism, namely $P_{1a}=\Hom_{R_{\fm_1}}(M_1\oplus R_{\fm_1}^{r_2},M_1)$ and $P_{1b}=\Hom_{R_{\fm_1}}(M_1\oplus R_{\fm_1}^{r_2},R_{\fm_1})$. Notice that $P_{1a}$ and $P_{1b}$ are not isomorphic because $M_1$ has rank $r_1>1$. In addition, by the Krull-Schmidt Theorem, all projective $\Lambda_{\fm_1}$-modules can be written in a unique way as a direct sum of copies of $P_{1a}$ and $P_{1b}$.

By the definition of $P_{1a}$, $\Tr(P_{1a})=I_1$. On the other hand, $(Q/QI)_{\fm_1}\cong (e(\Lambda/I))_{\fm_1}\cong (E_{11}+J) M_{r_2}(R_{\fm_1})$, so $Q_{\fm_1}\cong P_{1a}^{(\kappa_1)}\oplus P_{1b}$, where $\kappa_1$ is, at most, a countable cardinal.

Similarly, there are two finitely generated indecomposable projective $\Lambda_{\fm_2}$-modules up to isomorphism $P_{2a}=\Hom_{R_{\fm_2}}(R_{\fm_2}^{r_1}\oplus M_2,M_2)$ and $P_{2b}=\Hom_{R_{\fm_2}}(R_{\fm_2}^{r_1}\oplus M_2,R_{\fm_2})$. By the definition of $P_{2a}$, $\Tr(P_{2a})=I_2$. Therefore, since $(Q/QI)_{\fm_2}\cong (e(\Lambda/I))_{\fm_2}=\{0\}$, $Q_{\fm_2}\cong P_{2a}^{(\kappa_2)}$, where $\kappa_2$ is, at most, a countable cardinal. 

Let $L=Q\otimes_\Lambda (X_1\oplus X_2)$, which is in $\Add(X_1\oplus X_2)$. By Proposition~\ref{equivalencia}, $L$ is a direct summand of $(X_1\oplus X_2)^{(\omega)}$ and $\Hom_R(X_1\oplus X_2,L)\cong Q$. Moreover, by Proposition~\ref{equivalencia}, 
    \[L_{\fm_1}\cong M_1^{(\kappa_1)}\oplus R_{\fm_1} \quad(*) \qquad\text{and}\qquad
        L_{\fm_2}\cong M_2^{(\kappa_2)}\quad (**)\]
If $Q$ is finitely generated, $\kappa_1$ and $\kappa_2$ are finite. But then, by the isomorphism~$(*)$, the rank of $L$ is congruent to 1 modulo $d$ and, by the isomorphism~$(**)$, the rank of $L$ is divisible by $d$, a contradiction.

To prove the second statement in the claim, assume that $Q = \bigoplus_{i \in \N} Q_i$, where the $Q_i$'s are finitely generated. Since $Q/QI \cong e(\Lambda/I)$ is finitely generated, $Q_i = Q_iI$ for almost 
all $i \in \N$. If $I_0 := \{i \in \N \mid Q_i \neq Q_iI\}$, then $Q_f = \bigoplus_{i \in I_0} Q_i$ is a finitely generated projective right $\Lambda$-module such that 
$Q_f/Q_fI \cong e(\Lambda/I)$. By the previous part of the proof, such $Q_f$ does not exist, hence $Q$ is not a direct sum of finitely generated modules. 
By Proposition~\ref{equivalencia}, $L \in \Add (X_1 \oplus X_2)$ is not a direct sum of finitely generated modules. This finishes the proof of $(2)$.
\end{Proof}

Next Proposition reviews a construction, that goes back to an idea of Bass \cite{bass}, of an indecomposable finitely generated torsion-free module of rank two over any local domain that has a finitely generated ideal that cannot be generated by two elements. 

This construction is particularly relevant to us in view of Theorem~\ref{coprime}. It is going to imply that if $\mathcal{F}$ is closed under direct summands, then finitely generated ideals of $R_\fm$ are two-generated for all maximal ideals of $R$ except maybe one, see Theorem~\ref{th2generated}.

\begin{Prop} \label{bassmodule} 
Let $R$ be a commutative local domain, with field of fractions $Q$, and maximal ideal $\fm$. Let $a,b,c$ be elements in $R$ such that the ideal $K=aR+bR+cR$ cannot be generated by two elements. Set $\alpha=(a,b,c)\in R^3$ and let $H=R^3\cap Q\alpha$. Then:
\begin{enumerate}
    \item[(i)] $M=R^3/H$ is  an indecomposable, finitely generated torsion-free module of rank $2$.
    \item[(ii)] If $\overline{R}$ is local, then so is $\End_R(M)$.
\end{enumerate}
\end{Prop}
\begin{Proof}
$(i)$. Let $k=R/\fm$ be the residue field of $R$. Let $M=A_1/H\oplus A_2/H$ be a nontrivial decomposition with $A_i$ $R$-submodules of $R^3$ containing $H$. Consider $M\otimes _R k\cong \left(A_1/H \otimes _R k\right) \oplus \left( A_2/H\otimes _R k \right)$. Then,  since $M\otimes _R k$ is at most a three-dimensional vector space, one of the direct summands is one-dimensional. We may assume it is $A_2/H\otimes _R k $, which in turn implies that $A_2/H\cong R^3/A_1$ is isomorphic to $R$. Therefore, the exact sequence
    \[\begin{tikzcd}
        0\rar & A_1\rar & R^3 \rar & R^3/A_1\rar & 0
    \end{tikzcd}\]
splits, and we deduce that $A_1$ is a $2$-generated free module. So that if $v_1$ and $v_2$ is a basis of $A_1$, there exists $r_1, r_2\in R$ such that $(a,b,c)=v_1r_1+v_2r_2$. This implies that $K\subseteq r_1R+r_2R$. Since $v_1$ and $v_2$ can be completed to a basis of $R^3$, it follows that $K=r_1R+r_2R$, which contradicts the assumption that $K$ cannot be generated by two elements. This finishes the proof that $M$ is indecomposable.

$(ii)$. The module $M$ fits into the exact sequence
    \[\begin{tikzcd}
        0\rar & H\rar & R^3\rar & M\rar & 0.
    \end{tikzcd}\]
Let $T=\{f\in \End_R(R^3)\mid f(H)\subseteq H\}$, which is a subring of $\End_R(R^3)$. Then there is an onto ring endomorphism  $\varphi \colon T\to \End_R(M)$. We will prove that $T$ is a local ring, and then so is $\End_R(M)$.

Restriction to $H$ induces a ring morphism $\psi \colon T\to \End_R(H)$ with kernel  $I=\{f\in T\mid f(H)={0}\}$. 

We claim that the embedding $T\hookrightarrow \End_R(R^3)$ is local. Indeed, let $f\in T$ be such that $f$ is invertible in $\End_R(R^3)$. Then there is a commutative diagram
\[\begin{tikzcd}[column sep=huge,row sep=huge]
	0\rar&	H\rar\dar{\psi(f)} & R^3\dar{f}\rar&M\dar{\varphi (f)}\rar &0\\
	0\rar&	H\rar & R^3\rar&M\rar &0
	\end{tikzcd}\]
Since $f$ is invertible, it is onto and then so is $\varphi(f)$. Since $M$ is a finitely generated module over a commutative ring, we deduce that $\varphi(f)$ is bijective. Now the Snake Lemma implies that $\psi(f)$ is also bijective, and then  we can deduce that $f^{-1}\in T$.

Now we prove that $I\subseteq J(T)$. Notice that if $f\in I$, then it induces a module homomorphism $\tilde{f}\colon M\to R^3$ such that $\im f=\im \tilde{f}$. Since $M$ is indecomposable, and $R$ is local $\im \tilde{f}\subseteq R^3\fm$. Hence, $1-f$ is invertible in  $\End_R(R^3)$, and since the embedding $T\hookrightarrow \End_R(R^3)$ is local, we deduce that $1-f$ is invertible in $T$. This proves that the two-sided ideal $I\subseteq J(T)$.

The image of the ring morphism $\psi$ is the ring   
    $$S=\{g\in \End_R(H)\mid g\mbox{ can be extended to an endomorphism of $R^3$}\}.$$
Hence $T/I\cong S$, and $T$ is local provided $S$ is local. Notice that for any $g\in S$ its extension to $R^3$ satisfies a monic polynomial of degree $3$ with coefficients in $R$, so $g$ also satisfies that polynomial. This is to say that the extension $R\hookrightarrow S$ is integral.

Since $H$ is a torsion-free module of rank $1$, $\End_R(H)$ can be identified with a subring of $Q$. Then $S$ is a subring of $\overline R$. Since $\overline R$ is local, so is $S$. This finishes the proof of the result.
\end{Proof}

The next Corollary shows that the conclusions of Proposition~\ref{bassmodule} can be extended to any  ring $S$ between $R$ and its integral closure.

\begin{Cor} \label{bassmoduleextensions}
Let $R$ be a commutative local domain, with field of fractions $Q$ and maximal ideal $\fm$ such that $\overline{R}$ is local. Let $S$ be an intermediate ring $R\subseteq S\subseteq \overline R$.
Let $a,b,c$ be elements in $S$ such that the ideal $K=aS+bS+cS$ cannot be generated by two elements. Set $\alpha=(a,b,c)\in S^3$ and let $H=S^3\cap Q\alpha$. Then:
\begin{enumerate}
    \item[(i)] $M=S^3/H$ is a torsion-free $R$-module of rank $2$ and $\End_R(M)=\End_S(M)$ is a local ring. In particular, $M$ is indecomposable.
    \item[(ii)] Consider the ring $T=R[a,b,c]$ which is a finite integral extension of $R$. Let $H'=H\cap T^3$. Then $M'=T^3/H'$ is an indecomposable, finitely generated $R$-module of rank $2$ with local endomorphism ring.
\end{enumerate}
\end{Cor}
\begin{Proof}
Since $\overline{R}$ is local, so is any intermediate ring $R\subseteq S\subseteq \overline{R}$. 

Notice also that an $S$-module $M$ has finite rank $n$ as $S$-module if and only if it has finite rank $n$ as $R$-module because $M\otimes _S Q=M\otimes _S S\otimes _RQ=M\otimes _RQ$. This implies, in addition, that $\End_S(M)\cong \{A\in M_n(Q)\mid AM\subseteq M\}$ coincides with $\End_R(M)$.

Therefore, the result is a consequence of Proposition~\ref{bassmodule} applied to $S$ in $(i)$ and to $T$ in $(ii)$.
\end{Proof}

The following lemma is a variation for domains of finite character of a well-known fact about bounds on the number of generators of finitely generated ideals.

\begin{Lemma} \label{local2generated} \cite[Proposition~1.4]{basstf}
Let $R$ be a commutative domain of finite character. Let $k\ge 2$. If $I$ is a non-zero finitely generated ideal of $R$ such that $IR_\fm$ is $k$-generated for every maximal ideal $\fm$ of $R$, then $I$ is also $k$-generated.
\end{Lemma}
\begin{Proof}
Let $I$ be a non-zero ideal of $R$. Since $R$ is of finite character, $I$ is contained only in finitely many maximal ideals $\{\fm_1,\dotsc,\fm_n\}$ of $R$. 
First, observe that there exists $\alpha \in I$ such that $\frac{\alpha}{1} \not \in \fm_i I_{\fm_i}$ for every $i = 1,\dots,n$. Indeed, 
let $\alpha_i \in I$ be such that $\frac{\alpha_i}{1} \not \in \fm_i I_{\fm_i}$ and let $e_1,\dots,e_n$ be such that 
$e_i - 1 \in \fm_i$ for every $i=1,\dotsc,n$ and $e_i \in \fm_j$ whenever $1 \leq  i \neq j\leq n$. Then let $\alpha = \sum_{i = 1}^n e_i \alpha_i$.

Since $R$ is of finite character, $\alpha$ is contained only in finitely many maximal ideals of $R$, say  $\mathcal M = \{\fm_1,\dotsc,\fm_n,\fm_{n+1},\dots,\fm_{n+\ell}\}$.

For each $i = 1,\dots,n$, there are $\alpha_{i,1},\dots,\alpha_{i,k-1} \in I$ such that $I_{\fm_i} = \frac{\alpha}{1}R_{\fm_i} + \sum_{j = 1}^{k-1} \frac{\alpha_{i,j}}{1}R_{\fm_i}$.
For each $i = n+1,\dots,n + \ell,j = 1,\dots,k-1$ let $\alpha_{i,1} \in I \setminus \fm_{i}$ and $\alpha_{i,j} = 0$ if $j>1$.
Notice that $I_{\fm_i} = \frac{\alpha}{1}R_{\fm_i} + \sum_{j = 1}^{k-1} \frac{\alpha_{i,j}}{1}R_{\fm_i}$ for every $i = 1,\dots,n+\ell$.
As before, consider $f_1,\dots,f_{n+\ell} \in R$ such that $f_i - 1 \in \fm_i$ and $f_i \in \fm_j$ if $1 \leq j \neq i \leq n + \ell$.
Set $\beta_{i,j} = f_i \alpha_{i,j}$ for every $i = 1,\dots,n+\ell$ and $j = 1,\dots,k-1$.
For $j = 1,\dots,k-1$ set $\alpha_j := \sum_{i = 1}^{n+\ell} \beta_{i,j}$. We claim that $I =  \alpha R +  \sum_{j = 1}^{k-1} \alpha_j R$.
It is sufficient to verify this equality locally. If $\fm \not \in \mathcal M$ then both sides localize to $R_{\fm}$.
Also, $\alpha R_{\fm_i} + \sum_{j = 1}^{k-1}\alpha_jR_{\fm_i} + \fm_iI_{\fm_i} = \alpha R_{\fm_i} + \sum_{j = 1}^{k-1} \alpha_{i,j}R_{\fm_i} + \fm_iI_{\fm_i}$.
Since $I$ is finitely generated, we can conclude by Nakayama's lemma. 
\end{Proof}

Now we are ready to prove the main result of the section.

\begin{Th} \label{th2generated}
Let $R$ be a commutative domain of finite character, and of Krull dimension $1$. Assume that for any finitely generated torsion-free $R$-module $X$, every element of $\Add (X)$ is a direct sum of finitely generated modules. Then $R$ satisfies the following properties:
\begin{enumerate}
    \item[(1)] For any maximal ideal $\fm$ of $R$ and any ring $S$ such that $R_\fm\subseteq S \subseteq \overline{R_\fm}$ we have that finitely generated indecomposable torsion-free $S$-modules have local endomorphism ring;
    \item[(2)] for any maximal ideal $\fm$ of $R_\fm$, except maybe one maximal ideal $\fm _0$, all finitely generated ideals of $R_\fm$ are at most two-generated. Then $\overline{R_\fm}$ is a valuation ring for any maximal ideal $\fm \neq \fm _0$.
    \item[(3)] for each maximal ideal $\fm$ of $R$, there is a unique  maximal ideal of $\overline{R}$ lying over $\fm$. In particular, $\overline{R}$ has also finite character.
\end{enumerate}
If $\overline{R_\fm}$ satisfies the two-generated property for any maximal ideal $\mathfrak{m}$, then $\overline{R}$ is a Pr\"ufer domain of Krull dimension $1$, so any finitely generated ideal of $R$ is, at most, two-generated. 
\end{Th}

\begin{Proof}
By Corollary~\ref{generalrank} and Proposition~\ref{propertylocalizes}, we deduce that for any maximal ideal $\fm$ of $R$, any ring $S$ such that $R_\fm\subseteq S \subseteq \overline{R_\fm}$ satisfies that every finitely generated indecomposable torsion-free $S$-module has local endomorphism ring. This shows~(1).

Notice  that $(1)$ implies that $\overline{R}$ is also a ring of finite character because there is just one maximal ideal of $\overline{R}$ lying over each maximal ideal of $R$. This shows $(3)$.

By Theorem~\ref{coprime}, the ranks of two finitely generated, indecomposable, torsion-free modules over different localizations of $R$ at maximal ideals must be coprime. In view of Corollary~\ref{bassmoduleextensions}, we deduce that finitely generated ideals of $R_\fm$ are $2$-generated for any maximal ideal $\fm$ of $R$  except maybe for one  that we will denote by $\fm _0$. By \cite[Proposition~III.1.11]{fuchssalce}, $\overline{R_\fm}$ is a valuation ring for any maximal ideal $\fm \neq \fm _0$. This proves statement $(2)$.

If the maximal ideal $\fm _0$ does not exist,  then $\overline{R}$ is a Pr\"ufer domain of Krull dimension $1$. Moreover, by Lemma~\ref{local2generated}, any finitely generated ideal of $R$ is, at most, $2$-generated.
\end{Proof}

\section{The integrally closed case}\label{s:ic}

In this section, we will use the results developed until now to show that the converse of Theorem~\ref{th2generated} is true for integrally closed $h$-local domains of Krull dimension $1$. This will be done in Corollary~\ref{integrallyclosed}.

We recall the following definition, which will be needed for this section.

\begin{Def}
Let $R$ be a ring. Let $M$ be a right $R$-module. A submodule $N$ of $M$ is called \emph{relatively divisible} or an \emph{$\mathrm{RD}$-submodule} if $Nr=N\cap Mr$ for each $r\in R$.
\end{Def}

In particular, this definition tells us that if for some $r\in R$, $a\in N$, the equation $xr=a$ has a solution in $M$, then it has a solution in $N$ as well. In this sense, pure submodules are a generalization of $\mathrm{RD}$-submodules to the case where instead of an equation, we have a system of linear equations which has a solution in $N^{n}$ whenever it has a solution in $M^{n}$, for some positive integer $n\ge1$. 

For a discussion of the basic properties of $\mathrm{RD}$-submodules, the reader is referred to \cite[Chapter~I.7]{fuchssalce}.

\begin{Prop} \label{prop:ic}
Let $R$ be an $h$-local domain of Krull dimension $1$, and with field of fractions $Q$. Assume that $R_\fm$ is a valuation domain for all maximal ideals $\fm$ of $R$ except maybe one maximal ideal $\fm _0$. 

Then the class of $R$-modules that are direct sums of finitely generated torsion-free $R$-modules is closed under direct summands if and only if the class of direct sums of finitely generated torsion-free $R_{\fm_0}$-modules is closed under direct summands. 
\end{Prop}

\begin{Proof} 
By Proposition~\ref{propertylocalizes}, if the class of $R$-modules that are direct sums of finitely generated torsion-free $R$-modules is closed under direct summands, then so is the corresponding class for $R_{\fm _0}$. We need to prove the converse result.

Let $\{N_i\mid i \in \N\}$ be a countable family of non-zero finitely generated torsion-free $R$-modules and set $N=\bigoplus_{i \in \N} N_i$. Let $A$ be a direct summand of $N$. Then $A_{\fm_0}$ is a direct summand of $N_{\fm_0}$. By hypothesis, $A_{\fm _0}$ is a direct sum of  finitely generated $R_{\fm _0}$-modules, say $A_{\fm _0} = \bigoplus_{i \in \N} X_i$. For $i \in \N$, let $A_i := \{a \in A \mid \frac{a}{1} \in \bigoplus_{1 \leq j \leq i} X_j\}$. We claim  that $A_i$ is an $\mathrm{RD}$-submodule of $A$. Let $a \in A$ and $r\in R\setminus \{0\}$ be such that $\frac a1 \frac r1 \in \bigoplus_{1 \leq j \leq i} X_j$, since $\frac r1$ is not a zero divisor of $X_\ell$ for any $\ell$, we deduce that $\frac a1\in A_i$. This finishes the proof of the claim. 

Note also that there exists an $\ell\ge1$ such that $\bigoplus_{1 \leq j \leq i} X_j$ is a submodule of $\bigoplus_{i=1}^\ell  (N_i)_{\fm _0}$. Hence, since all modules involved are torsion-free, $A_i\le \bigoplus_{i=1}^\ell  N_i$. We claim that this is also an $\mathrm{RD}$-embedding. Assume that $n\in \bigoplus_{i=1}^\ell  N_i$ is such that there exists $r\in R$ such that $nr\in A_i\le A$. As $A$ is an $\mathrm{RD}$-submodule of $N$, we deduce that $n\in A$, and since $A_i$ is an $\mathrm{RD}$-submodule of $A$ we deduce that $n\in A_i$.

Now we shall prove  that $A_i$ is finitely generated. Since $A_i$ is a submodule of a free module $R^{r}$, by Lemma~\ref{almostfg} there exist $a_1,\dots,a_s \in A_i$ and a finite set of maximal ideals of the form $\mathcal{S}=\{\fm _0,\fm _1, \dots ,\fm _t\}$ such that $(A_i)_\fm =\sum _{j=1}^s\frac{a_j}1R_\fm$ for any maximal ideal $\fm\notin\mathcal{S}$. In particular, $(A_i)_Q=\sum _{j=1}^s\frac{a_j}1Q$ and in the  exact sequence,
    \[\begin{tikzcd}
        0\rar & \displaystyle\sum _{j=1}^s{a_j}R \rar & A_i \rar & A_i/\Big(\displaystyle\sum _{j=1}^s{a_j}R\Big)\rar & 0
    \end{tikzcd}\]
the module $A_i/(\sum _{j=1}^s{a_j}R)$ is torsion. Since $R$ is an $h$-local domain of Krull dimension $1$,
    \[A_i/\Big(\sum _{j=1}^s{a_j}R\Big)=\bigoplus_{\fm \in \mSpec(R)} \Big( A_i/\Big(\sum _{j=1}^s{a_j}R\Big)\Big)_\fm =\bigoplus_{\fm \in \mathcal{S}} \Big( A_i/\Big(\sum _{j=1}^s{a_j}R\Big)\Big)_\fm.\]
Hence, to prove the claim, it is enough to show that $(A_i)_{\fm _j}$ is finitely generated as a $\Lambda _{\fm _j}$-module for any $j=0,\dots ,t$, cf. Lemma~\ref{criteriafghlocaltorsion}.

Note that $(A_i)_{\fm _0}=\bigoplus_{1 \leq j \leq i} X_j$ is finitely generated. Assume now that $j\ge 1$. Since $A_i$ is an $\mathrm{RD}$-submodule of $\bigoplus_{k=1}^\ell  N_k$, we deduce that $(A_i)_{\fm _j}$ is an $\mathrm{RD}$-submodule of $\bigoplus_{k=1}^\ell  (N_k)_{\fm _j}$. Since $R_{\fm _j}$ is a valuation domain, and the module $\bigoplus_{k=1}^\ell  (N_k)_{\fm _j}/(A_i)_{\fm _j}$ is finitely generated and torsion-free, it is projective. Therefore, $(A_i)_{\fm _j}$ is a direct summand of $\bigoplus_{k=1}^\ell  (N_k)_{\fm _j}$, so it is finitely generated, as claimed. 

Now, $A$ is a union of a chain of finitely generated modules $A_1\subseteq A_2 \subseteq A_3\subseteq \cdots$. For any $i \in \N$, the exact sequence $0 \to A_i \to A_{i+1} \to A_{i+1}/A_i \to 0$ splits upon localization at $\fm _0$ by construction. It also splits when localized at other maximal ideals since it is $\mathrm{RD}$-exact, $R_{\fm}$ is a  valuation domain, and all the involved modules are finitely generated and torsion-free, and therefore projective.  

By Lemma~\ref{locallysplit}, the sequence $0 \to A_i \to A_{i+1} \to A_{i+1}/A_i \to 0$ splits in $R$ for any $i\in \N$. Therefore, $A \cong A_1\oplus (\bigoplus _{i\in \N} A_{i+1}/A_i)$ is a direct sum of finitely generated modules.
\end{Proof}

\begin{Cor} \label{integrallyclosed} 
Let $R$ be an $h$-local domain of Krull dimension $1$ that is integrally closed in its field of fractions $Q$. Then the following statements are equivalent:
\begin{enumerate}
    \item[(i)] The class of $R$-modules that are direct sums of finitely generated torsion-free $R$-modules is closed under direct summands.
    \item[(ii)] For any finitely generated torsion-free $R$-module $X$, every element in $\Add(X)$ is a direct sum of finitely generated modules.
    \item[(iii)] The ring $R$ satisfies one of the following conditions:
    \begin{itemize}
        \item[(1)] $R$ is a Pr\"ufer domain; or
        \item[(2)] there is a maximal ideal $\fm _0$ such that for every maximal ideal $\fm\neq\fm_0$, $R_\fm$ is a valuation domain, $R_{\fm _0}$ is an integrally closed domain that is not a valuation domain, and every indecomposable finitely generated torsion-free $R_{\fm_0}$-module has local endomorphism ring.
    \end{itemize}
\end{enumerate}
\end{Cor}
\begin{Proof}
It is clear that $(i)$ implies $(ii)$. Theorem~\ref{th2generated} shows that $(ii)$ implies $(iii)$.

Assume $(iii)$ holds. If $R$ is a Pr\"ufer domain, then it is semihereditary, so finitely generated torsion-free modules are projective, and all projective modules are the direct sum of finitely generated ideals. So $(i)$ is trivially satisfied. 

Now assume that $R$ satisfies the condition  $(2)$. By Corollary~\ref{generalrank}, $R_{\fm _0}$ satisfies that the class of finitely generated torsion-free modules is closed under direct summands. Then $(i)$ follows from Proposition~\ref{prop:ic}.
\end{Proof}

\begin{Cor} \label{integralclosure} 
Let $R$ be an $h$-local domain of Krull dimension $1$. Assume that the class of $R$-modules that are direct sums of finitely generated torsion-free modules is closed under direct summands. Then its integral closure also satisfies this property.
\end{Cor}
\begin{Proof}
By Theorem~\ref{th2generated}, $\overline{R}$ satisfies conditions $(1)$ and $(2)$ of Corollary~\ref{integrallyclosed} so we can deduce that the class of finitely generated torsion-free $\overline{R}$-modules is closed by direct summands.
\end{Proof}

\section{Infinite direct sums are determined by the genus}\label{s:genus}

\begin{Def}
Let $R$ be a commutative ring, and let $\Lambda$ be an $R$-algebra. Let $M$ and $N$ be right $\Lambda$-modules. We say that $M$ and $N$ are \textit{in the same genus} if $M_\fm$ is isomorphic to $N_\fm$ for every maximal ideal $\fm$ of $R$.
\end{Def}

\begin{Lemma} \label{genus}
Let $R$ be a commutative ring, and let $\Lambda$ be a module-finite $R$-algebra. Let $M,N,X$ be right $\Lambda$-modules. If $X$ is finitely generated over $R$ and $M\oplus X$ and $N\oplus X$ are in the same genus, then $M$ and $N$ are in the same genus.
\end{Lemma}
\begin{Proof}
Recall that if a module $Y$ has semilocal endomorphism ring, then it has the cancellation property, that is, $Y\oplus A\cong Y\oplus B$ implies $A\cong B$ (see, for example, \cite[Theorem 4.5]{libro}). 

By Proposition~\ref{fgendomorphism}, $X_{\fm}$ has semilocal endomorphism ring, so if $M_{\fm} \oplus X_{\fm} \cong N_{\fm} \oplus X_{\fm}$ then $M_{\fm} \cong N_{\fm}$ for every maximal ideal $\fm$ of $R$.
\end{Proof}

\begin{Lemma}\label{fg.dsummands}
Let $R$ be a commutative domain of finite character with field of fractions $Q$. Let $\Lambda$ be a module-finite $R$-algebra such that $\Lambda_Q$ is a simple artinian ring. Let $M=\bigoplus_{i\in I} M_i$ be an infinite direct sum of non-zero finitely generated right $\Lambda$-modules which are torsion-free as $R$-modules, and let $N$ be a finitely generated right $\Lambda$-module which is torsion-free as an $R$-module. If $N_\fm$ is a direct summand of $M_\fm$ for every maximal ideal $\fm$ of $R$, then $N$ is a direct summand of $M$. 
\end{Lemma}
\begin{Proof}
Suppose that $N_\fm$ is a direct summand of $M_\fm$ for every maximal ideal $\fm$ of $R$. Since $M_Q$ is a free $\Lambda_Q$-module and $N$ is finitely generated, there is some finite subset $I_1\subseteq I$ such that $N_Q$ is a direct summand of $\bigoplus_{i\in I_1}(M_i)_Q$. By Lemma~\ref{almostallmaximals}, there is a $\Lambda$-module homomorphism $f\colon\bigoplus_{i\in I_1} M_i\to N$ such that the induced homomorphism $f_\fm\colon\bigoplus_{i\in I_1}(M_i)_\fm\to N_\fm$ is a splitting epimorphism for almost all maximal ideals $\fm$ of $R$.

Let $\mathcal{M}$ be the finite set of maximal ideals $\fm$ such that $f_\fm$ is not a splitting epimorphism. Since $N$ is finitely generated, there is some finite subset $I_2\subseteq I$ containing $I_1$ such that $N_\fm$ is a direct summand of $\bigoplus_{i\in I_2}(M_i)_\fm$ for every $\fm\in\mathcal{M}$. Then $N_\fm$ is a direct summand of $\bigoplus_{i\in I_2}(M_i)_\fm$ for all $\fm\in\mSpec R$. 

Since $\Lambda_Q$ is simple artinian and $N$ is finitely generated, there exists $I_3\subseteq I\setminus I_2$ finite such that $N_Q$ is a direct summand of $\bigoplus_{i\in I_3}(M_i)_Q$. Then, by Proposition~\ref{directNK}, $N$ is a direct summand of $\bigoplus_{i\in I_0} M_i$, where $I_0=I_2\cup I_3$. Hence, $N$ is a direct summand of $M$.
\end{Proof}

\begin{Th}\label{countablesum}
Let $R$ be a commutative domain of finite character with field of fractions $Q$. Let $\Lambda$ be a module-finite $R$-algebra such that $\Lambda_Q$ is a simple artinian ring. Let $M=\bigoplus_{i\in\No} A_i$ and $N=\bigoplus_{i\in\No} B_i$ be direct sums of non-zero finitely generated right $\Lambda$-modules which are torsion-free as $R$-modules. If $M$ and $N$ are in the same genus, then there are decompositions
    \[\textstyle M=\bigoplus_{i\in\No} M_i\qquad\text{and}\qquad N=\bigoplus_{i\in\No} N_i\]
such that both $M_i$ and $N_i$ are finitely generated, and $M_i\cong N_i$ for every $i\in\No$. In particular, $M$ and $N$ are isomorphic.
\end{Th}

\begin{Proof}
Let $\{\fm_\alpha\}_{\alpha\in\Omega}$ denote the set of maximal ideals of $R$. Let $\{a_i\}_{i\in\No}$ and $\{b_j\}_{j\in\No}$ be countable sets of generators for $M$ and $N$, respectively. By Lemma~\ref{fgendomorphism}, $(A_i)_{\fm_\alpha}$ and $(B_i)_{\fm_\alpha}$ have semilocal endomorphism rings for every $\alpha\in\Omega$, and $i\in\N$. Hence, both $A_i$ and $B_i$ satisfy the cancellation property locally for every $i\in \N$. 

We claim that there exist $m_0,n_0\in\N$ and decompositions 
    \[\textstyle M=M_0\oplus Z_0\oplus \bigoplus_{i=m_0+1}^\infty A_i\quad\text{and}\quad N=N_0\oplus \bigoplus_{j=n_0+1}^\infty B_j\]
for some finitely generated modules $M_0,N_0,Z_0$ such that $M_0$ is isomorphic to $N_0$, $a_0\in M_0$, $b_0\in N_0$, and $Z_0\oplus\bigoplus_{i=m_0+1}^\infty A_i$ and $\bigoplus_{j=n_0+1}^\infty B_j$ are in the same genus. 

Let $m_0'=\min\{m\in\No\mid a_0\in\bigoplus_{i=0}^m A_i\}$, and consider $X_0=\bigoplus_{i=0}^{m_0'} A_i$. By Lemma~\ref{fg.dsummands} applied to $X_0$, there is $n_0'\in\No$ such that $X_0$ is isomorphic to a direct summand of $\bigoplus_{j=0}^{n_0'} B_j$. Therefore, we can write 
    \[\textstyle M=X_0\oplus \bigoplus_{i=m_0'+1}^\infty A_i\qquad\text{and}\qquad N=Y_0\oplus Z_0'\oplus \bigoplus_{j=n_0'+1}^\infty B_j\]
with $X_0\cong Y_0$ and $Y_0 \oplus Z_0' = \bigoplus_{j = 0}^{n_0'} B_j$. By Lemma~\ref{genus}, $\bigoplus_{i=m_0+1}^\infty A_i$ and $Z_0'\oplus \bigoplus_{j=n_0'+1}^\infty B_j$ are in the same genus. 

Let $n_0=\max\{n_0'+1,\min\{n\in\No\mid b_0\in\bigoplus_{j=0}^n B_j\}\}$, and consider $Y_1=Z_0'\oplus\bigoplus_{j=n_0'+1}^{n_0} B_j$. By Lemma~\ref{fg.dsummands} applied to $Y_1$, there is $m_0>m_0'$ such that $Y_1$ is isomorphic to a direct summand of $\bigoplus_{i=m_0'+1}^{m_0} A_i$. Then
    \[\textstyle M=X_0\oplus X_1\oplus Z_0\oplus\bigoplus_{i=m_0+1}^\infty A_i\qquad\text{and}\qquad N=Y_0\oplus Y_1\oplus \bigoplus_{j=n_0+1}^\infty B_j\]
with $X_0\cong Y_0$ and $X_1\cong Y_1$. By Lemma~\ref{genus}, $Z_0\oplus\bigoplus_{i=m_0+1}^\infty A_i$ and $\bigoplus_{j=n_0+1}^\infty B_j$ are in the same genus. Take $M_0=X_0\oplus X_1$ and $N_0=Y_0\oplus Y_1$. We can repeat this argument $\omega$ times to obtain ascending chains of positive integers $m_0<\dotsb< m_\ell<\dotsb$ and $n_0<\dotsb< n_\ell<\dotsb$, and decompositions
    \[\textstyle M=(\bigoplus_{i=0}^\ell M_i)\oplus Z_\ell \oplus (\bigoplus_{i=m_\ell+1}^\infty A_i)\quad\text{and}\quad N=(\bigoplus_{i=0}^\ell N_i)\oplus (\bigoplus_{j=n_\ell+1}^\infty B_j)\]
such that $M_i$ is finitely generated and isomorphic to $N_i$ for every $i=0,\dotsc,\ell$, $a_0,\dots,a_\ell\in\bigoplus_{i=0}^\ell M_i$, $b_0,\dots,b_\ell\in\bigoplus_{i=0}^\ell N_i$, and $Z_\ell \oplus (\bigoplus_{i=m_\ell+1}^\infty A_i)$ and $\bigoplus_{j=n_\ell+1}^\infty B_j$ are in the same genus. Therefore, we obtain two families of finitely generated $\Lambda$-modules which are torsion-free as $R$-modules $\{M_i\}_{i\in\No}\subseteq M$ and $\{N_i\}_{i\in\No}\subseteq N$ such that $M_i$ is finitely generated and isomorphic to $N_i$, $\{a_i\}_{i\in\No}\subseteq \bigoplus_{i\in\No} M_i$, and $\{b_i\}_{i\in\No}\subseteq \bigoplus_{i\in\No} N_i$. We deduce that $M=\bigoplus_{i\in\No} M_i$ and $N=\bigoplus_{i\in\No} N_i$. Therefore, $M$ is isomorphic to $N$. This finishes the proof of the proposition.
\end{Proof}

The following proposition generalizes Lemma~\ref{fg.dsummands}. 

\begin{Prop}
Let $R$ be a commutative domain of finite character with field of fractions $Q$. Let $\Lambda$ be a module-finite $R$-algebra such that $\Lambda_Q$ is a simple artinian ring. Let $M=\bigoplus_{i\in I} M_i$ be an infinite direct sum of non-zero finitely generated right $\Lambda$-modules which are torsion-free as $R$-modules, and let $A$ be a direct summand of $M$ of infinite rank. Let $N$ be a finitely generated right $\Lambda$-module which is torsion-free as an $R$-module. If $N_{\fm}$ is isomorphic to a direct summand of $A_{\fm}$ for each maximal ideal $\fm$ of $R$, then $N$ is isomorphic to a direct summand of $A$.
\end{Prop}
\begin{Proof}
Since $A$ is a direct summand of $M$, there are $\Lambda$-homomorphisms ${\iota\colon A\hookrightarrow M}$ and $\pi\colon M\to A$ such that $\pi \iota = 1_A$. For every subset $J \subseteq I$, let $\pi_J \colon M \to \bigoplus_{j \in J} M_j$ and $\iota_J \colon \bigoplus_{j \in J} M_j \hookrightarrow M$ denote the corresponding canonical projection and canonical embedding, respectively. 

Since $N_Q$ is a direct summand of $A_Q$, there are $\Lambda_Q$-homomorphisms ${f\colon N_Q\to A_Q}$ and $g\colon A_Q\to N_Q$ such that $g f = 1_{N_Q}$. By Lemma~\ref{isofg}, there is a $\Lambda$-module homomorphism $f_0\colon N\to A$ such that $f=f_0/s$ for some non-zero $s\in R$. Let $J_0$ be a finite subset of $I$ such that ${\rm Im}\ f_0 \subseteq \bigoplus_{j \in J_0} M_j$. Hence $g(\iota_{J_0}\pi_{J_0}\iota)_Qf = 1_{N_Q}$. 

Again, by Lemma~\ref{isofg}, there is a $\Lambda$-homomorphism $g'\colon \bigoplus_{j \in J_0} M_j\to N$ such that $g(\iota_{J_0})_Q = g'/t$ for some non-zero $t\in R$. Let $g_0 = g'\pi_{J_0} \iota $, and note that $g_0/t\colon A_Q\to N_Q$ is a $\Lambda_Q$-homomorphism satisfying that $(g_0/t)f = 1_{N_Q}$. Hence, we may assume that there exists a $\Lambda$-homomorphism $g_0\colon A\to N$ that factors through $\pi_{J_0}\iota$, and such that $g$ is of the form $g_0/t$ for some non-zero $t\in R$. In particular, $g_0f_0 = st 1_{N}$. Also note that $g_0(A \cap \bigoplus_{i \in I \setminus J_0} M_i) = 0$ and ${\rm Im}\ f_0 \subseteq \bigoplus_{j \in J_0} M_j$. 

Let $r_0 = st\in R$. Since $R$ is of finite character, $r_0$ is contained only in finitely many maximal ideals of $R$, say $\mathcal M=\{\fm_1,\fm_2,\dots,\fm_k\}$. A similar argument as above shows that there are $\Lambda$-homomorphisms $f_1,\dots,f_k\colon N\to A$,  $g_1,\dots,g_k\colon A\to N$, and non-zero elements $r_1,\dots,r_k \in R$ such that $g_if_i = r_i 1_{N}$ and $r_i \notin \fm_i$, for each $i=1,\dotsc,k$. These morphisms can be chosen such that there are finite subsets $J_1,J_2,\dots,J_k \subseteq I$ such that $g_j(A \cap \bigoplus_{i \in I \setminus J_j}M_i) = 0$ and ${\rm Im}\ f_j \subseteq \bigoplus_{i \in J_j} M_i$.

Let $J = \bigcup_{i = 1}^k J_i$, and let $J' := I \setminus J$. Consider an exact sequence 
    $$\begin{tikzcd}
        0 \rar & X \rar{\nu} & A \rar{\pi_{J}|_{A}} & \bigoplus_{j \in J} M_j,
    \end{tikzcd}$$ 
where $X = A \cap \bigoplus_{j \in J'} M_j$. Since $A$ is of infinite rank, $X$ has infinite rank as well, and hence there is a $\Lambda_Q$-monomorphism $f'\colon N_Q\to X_Q$. Then $\nu_{Q}f'$ splits because $\Lambda_Q$ is simple artinian. Let $g'\colon A_Q\to N_Q$ be a $\Lambda_Q$-homomorphism such that $g'\nu_{Q}f' = 1_{N_Q}$. As in the first part of the proof, we may assume that $g'$ factors through $(\pi_K\iota)_Q$ for some finite subset $K \subseteq J'$. Therefore, there are $\Lambda$-homomorphisms $f_{\infty}\colon N\to A$ and $g_{\infty}\colon A\to N$ such that $ g_{\infty}f_{\infty} = r_{\infty} 1_N$ for some non-zero $r_{\infty} \in R$. Moreover, ${\rm Im} f_{\infty} \subseteq X$ and $g_{\infty}(A \cap \bigoplus_{j \in J} M_j) = 0$.

Since $R$ is of finite character, $r_\infty$ is contained only in finitely many maximal ideals of $R$, say $\mathcal N=\{\fn_1,\fn_2,\dots,\fn_{\ell}\}$. By the Chinese Remainder Theorem, we can find $e_1,\dots,e_{\ell} \in R$ such that $e_{j}\equiv0\mod{\fn}_i$ whenever $j \neq i$ and $e_j\equiv1\mod\fn_j$, for every $j = 1,\dots, \ell$. Note that any maximal ideal of $R$ either does not contain $r_0$ or is in $\mathcal M$. Hence, by the definition of the $r_1,\dotsc,r_k$, no maximal ideal of $R$ contains  the whole set $\{r_0,\dots,r_k\}$. 

For every $j = 1,\dots,\ell$, let $i(j) \in \{0,\dots,k\}$ be such that $r_{i(j)} \not \in \fn_{j}$. Let $g\colon A\to N$ be the $\Lambda$-homomorphism given by $g = g_{\infty} + \sum_{j = 1}^\ell e_j g_{i(j)}$.

We claim that $g$ is a locally split epimorphism, so by Lemma~\ref{locallysplit} it is a split epimorphism. Note that $g f_{\infty} = g_{\infty} f_{\infty} = r_{\infty} 1_N$, so $g_{\fn}$ splits if  $\fn$ does not contain $r_{\infty}$. On the other hand, 
    $$gf_{i(t)} = \sum_{j = 1}^{\ell} e_{j}g_{i(j)}f_{i(t)} = e_{t}r_{i(t)}1_{N} + \sum_{\substack{1 \leq j \leq \ell\\ j \neq t}} e_{j} g_{i(j)}f_{i(t)}.$$
When we localize at $\fn_t$, the first summand is an invertible element, while the second summand is in $J({\End}_{\Lambda_{\fn_t}}(N_{\fn_t}))$ by Lemma~\ref{jacobson}. Hence, $g_{\fn}$ splits for any maximal ideal of $R$, and we conclude with the proof of the statement.
\end{Proof}

Now we are going to give an example showing that Lemma~\ref{fg.dsummands} is not true for direct summands that are not finitely generated. We start with the following well-known lemma that will provide us a source  of noetherian local domains of Krull dimension $1$ with local integral closure. So that indecomposable finitely generated torsion-free modules have local endomorphism ring (see, for example, Corollary~\ref{local}).

\begin{Lemma} \label{nummonoids} 
Let $\alpha _1,\dots ,\alpha _n$ be (non-zero) coprime elements of $\N$. Let $K$ be a field,  $R=K[t^{\alpha _1 },\dots ,t^{\alpha _n}]$ and $\fm =t^{\alpha _1}R+\cdots + t^{\alpha _n}R$. Set $\Sigma =R\setminus \fm$. Then:
\begin{itemize}
    \item[(i)] $R$ is a noetherian domain of Krull dimension $1$ and  field of quotients  $K(t)$. The integral closure of $R$ into its field of fractions  is $K[t]$.
    \item[(ii)] The integral closure of $R_\fm$ into its field of fractions is $K[t]_\Sigma =K[t]_{(t)}$. In particular, it is a local ring.
    \item[(iii)] Every finitely generated, indecomposable, torsion-free  $R_\fm$-module has local endomorphism ring.
\end{itemize}
\end{Lemma}

\begin{Proof}
$(i)$. The field of fractions of $R$ is a subfield of $K(t)$ that contains $R$ and also $t$ (because $\alpha _1,\dots ,\alpha _n$ are coprime) so it coincides with $K(t)$. Since $K[t]$ is a PID, it is integrally closed, and being an integral extension of $R$, it is the integral closure.

$(ii).$ By $(i)$ the extension $R_\fm \subseteq K[t]_\Sigma$ is integral. If $\fn$ is a maximal ideal of $K[t]_\Sigma$ then $\fn \cap R=\fm$. Therefore $t\in \fn$, so that $\fn =tK[t]_\Sigma$, and $K[t]_\Sigma$ is a local ring.

The statement $(iii)$ follows from $(ii)$ and Corollary~\ref{local}.
\end{Proof}

\begin{Ex}\label{counterex}
\textit{There is a semilocal noetherian domain $R$ of Krull dimension $1$, and $R$-modules $M$ and $N$ that can be written as infinite direct sums of non-zero finitely generated torsion-free $R$-modules such that $N_\fm$ is a direct summand of $M_\fm$ for every maximal ideal $\fm$ of $R$ but $N$ is not a direct summand of $M$.}
\end{Ex}
\begin{Proof}
Let $K$ be an infinite field. Let $R_1=K[t^2,t^3]_{(t^2,t^3)}$, and let $R_2=K[t^3,t^7]_{(t^3,t^7)}$. So Lemma~\ref{nummonoids} applies to $R_1$ and $R_2$.

The local domain $R_1$ is a Bass domain with just two indecomposable finitely generated torsion-free modules (or rank one) up to isomorphism. Namely, $X=R_1$ and $Y=\fm _1$ where $\fm _1$ denotes the maximal ideal of $R_1$. Note that, by Lemma~\ref{nummonoids}, $X$ and $Y$ have local endomorphism ring.

The domain $R_2$ has infinitely many indecomposable finitely generated torsion-free modules of all ranks $\ge 2$, because it fails to satisfy the Drozd-Roiter conditions (cf. \cite[Theorem~4.2]{leuschke}). More precisely, if we denote by $\fm _2$ the maximal ideal of $R_2$, the $R_2$-module $(\fm _2K[t]_{(t)}+R_2)/R_2$  is not cyclic because a minimal set of generators is  $t^4$ and $t^5$. 

We single out an infinite family $Z_1,Z_2,\dotsc$ of indecomposable finitely generated, torsion-free $R_2$-modules of rank $2$. Note that, by Lemma~\ref{nummonoids}, such modules have local endomorphism ring.

Let $\varphi\colon K(t) \to K(t)$ be the automorphism that fixes $K$ and $\varphi (t)=t+1$. Let $R$ be the ring that fits in the  pull-back diagram
    \[\begin{tikzcd}[column sep=huge,row sep=huge]
        R\rar\dar & R_1\dar[hook] \\
        R_2\rar{\varphi'} & K(t)
    \end{tikzcd}\]
where $\varphi '\colon R_2\hookrightarrow K(t) \stackrel{\varphi}{\to}K(t)$.

By \cite[Theorem~4.4]{WW}, $R$ is a noetherian domain of Krull dimension $1$ with exactly two maximal ideals $\fm$ and $\fn$ and satisfying that $R_{\fm}\cong R_1$ and $R_{\fn}\cong R_2$.

Apply the results on \cite{WW} (or just  Corollary~\ref{existence}) to  construct two sequences of finitely generated torsion-free $R$-modules $M_1,M_2,\dotsc,N_1,N_2,\dots$ such that
    \[(M_i)_\fn=(N_i)_\fn=Z_i,\qquad (M_i)_\fm=X\oplus Y,\qquad (N_i)_\fm=X\oplus X.\]
Call $M = \bigoplus_{i\in \N} M_i$, $N = \bigoplus  _{i\in \N} N_i$. Then $N$ is locally a direct summand of $M$, but $N$ is not isomorphic to a direct summand of $M$.   Indeed, if $N\oplus N'\cong M$ then $\left(\bigoplus_{i\in \N} Z_i\right)\oplus N' _\fn\cong \bigoplus_{i\in \N} Z_i$. Since $\{ Z_i\}_{i\in \N}$ are non-isomorphic and have local endomorphism rings, $N'_\fn=\{0\}$. Hence, as $N'$ is torsion-free, $N'=\{0\}$. But  $N_\fm \not \cong M_\fm$, a contradiction. 
\end{Proof}

\begin{Remark}
By Proposition~\ref{semilocal}, the ring $R$ in Example~\ref{counterex} satisfies that the class of direct sums of finitely generated torsion-free modules is closed under direct summands.
\end{Remark}

\begin{Prop}\label{genusdecomposition}
Let $R$ be a commutative domain of finite character. Let $M$ be a torsion-free $R$-module of countable rank. Let $\fm_0$ be a maximal ideal of $R$. Assume that
\begin{enumerate}
    \item[(a)] $M_\fm$ is a direct sum of finitely generated torsion-free modules of rank one with local endomorphism ring for every maximal ideal $\fm\neq\fm_0$, and
    \item[(b)] $M_{\fm_0}$ is a direct sum of finitely generated torsion-free modules.
\end{enumerate}
Then $M$ is in the same genus as a direct sum of finitely generated torsion-free modules if and only if
\begin{enumerate}
    \item[(i)] $M_\fm$ is a free module for all but countably many maximal ideals $\fm$ of $R$.
	\item[(ii)] If $\mathcal{M}=\{\fm\in\mSpec R\mid M_\fm\text{ is not free}\}$, and $r_\fm$ denotes the number of free direct summands in the decomposition of $M_{\fm}$, then for every $b\in \N$ the set $\{\fm \in \mathcal{M} \mid r_{\fm} \leq b\}$ is finite.
\end{enumerate}
\end{Prop}
\begin{Proof}
Let $N=\bigoplus_{i\in\N} N_i$, where each $N_i$ is a finitely generated torsion-free module, and assume that $M$ and $N$ are in the same genus. (i) follows from Corollary~\ref{almost-free-count-ds}. If $\mathcal{M}$ is finite, (ii) is clear. Assume $\mathcal{M}$ is infinite and let $b\in\N$. Then there are only finitely many maximal ideals $\fm\in\mathcal{M}$ such that at least one of $(N_1)_{\fm},\dots,(N_b)_{\fm}$ is not free, so $r_{\fm} \geq b$ for almost all $\fm \in \mathcal{M}$. This proves (ii).

Conversely, assume that (i) and (ii) are satisfied and let $M_\fm=\bigoplus_{i\in\N} M_{\fm,i}$, where each $M_{\fm,i}$ is a finitely generated torsion-free module and such that $M_{\fm,i}$ has rank one for every maximal ideal $\fm\neq\fm_0$. For each $i\in\N$, let $d_i$ denote the rank of the module $M_{\fm_0,i}$.

First, assume that $\mathcal{M}$ is finite. For each $i\in\N$, apply Package Deal Theorem~\ref{dealsubmodules} to $X_i(\fm_0)=M_{\fm_0,i}$ and $X_i(\fm)=\bigoplus_{i=d_0+\dotsb+d_{i-1}+1}^{d_0+\dotsb+d_i} M_{\fm,i}$ for each $\fm\neq\fm_0$ (note that, since $\mathcal{M}$ is finite, $X_i(\fm)$ is free for almost all maximal ideals $\fm$ of $R$). Then, for each $j\in\N$, there is an $R$-module $N_j$ such that $(N_j)_\fm\cong X_j(\fm)$ for every maximal ideal $\fm$ in $R$. Therefore, $\bigoplus_{j\in\N} N_j$ is in the same genus as $M$.

Now, assume that $\mathcal{M}$ is infinite and consider the case where $r_\fm$ is finite for every maximal ideal $\fm$ of $R$. First, write the elements of $\mathcal{M}$ in a sequence $\fm_1,\fm_2,\dotsc$ in such a way that $r_{\fm_1}\le r_{\fm_2}\le\dotsb$ is non-decreasing, and assume that the direct summands in the decompositions of $M_\fm$ are indexed such that $M_{\fm,1},\dotsc,M_{\fm,r_\fm}$ are the free direct summands. For each $i\in\N$, apply Package Deal Theorem~\ref{dealsubmodules} to $X_i(\fm_0)=M_{\fm_0,i}$ and $X_i(\fm)=\bigoplus_{j=d_0+\dotsb+d_{i-1}+1}^{d_0+\dotsb+d_i} M_{\fm,j}$ for each $\fm\neq\fm_0$ (note that, with this reordering and considering (i) and (ii) with $b=d_0+\dotsb+d_i$, $M_{\fm,i}$ is free for almost all maximal ideals). Then, for each $j\in\N$, there is an $R$-module $N_j$ such that $(N_j)_\fm\cong X_j(\fm)$ for every maximal ideal $\fm$ in $R$. Therefore, $\bigoplus_{j\in\N} N_j$ is in the same genus as $M$.

Finally, assume that $\mathcal{M}$ is infinite and there is at least one maximal ideal $\fm$ with $r_\fm$ infinite. 
Let $\fm_1,\fm_2,\dots$ be a (finite or infinite) list of elements in $\{\fm \in \mathcal{M} \mid r_{\fm}$ is finite$\}$ and
let $\fn_1,\fn_2,\dots$ be a (finite or infinite) list of elements in $\{\fn \in \mathcal{M} \mid r_{\fn}$ is infinite$\}$.
The first sequence is chosen such that $r_{\fm_1}\le r_{\fm_2}\le\dotsb$ is non-decreasing, and again we assume 
the direct summands in the decompositions of $M_{\fm_i}$ are indexed such that $M_{\fm,1},\dotsc,M_{\fm,r_{\fm_i}}$ are the free direct summands.
Moreover, we assume that the direct summands in the decompositions of $M_{\fn_i}$ are indexed such that $M_{\fn_i,1},\dotsc,M_{\fn_i,i}$ are free. 

Let us check that each $M_{\fm,j}$ is free for almost all $\fm \in \mathcal{M}$. If $M_{\fn_i,j}$ is not free, then $i<j$. If $M_{\fm,j}$ is not free for infinitely many $\fm$'s with finite $r_{\fm}$, then $b = j-1$ would contradict (ii). Hence, as before, we find finitely generated torsion-free $R$-modules $N_j$ such that $(N_j)_{\fm} \cong X_j(\fm)$ for every maximal ideal $\fm$ of $R$. Therefore, $\bigoplus_{j\in\N} N_j$ is in the same genus as $M$.
\end{Proof}

\section{The noetherian case}\label{s:noeth}

The results of Section~\ref{s:genus} allow us to show in this section that the converse of Theorem~\ref{th2generated} is true for  semilocal noetherian domains of Krull dimension $1$ (cf. Proposition~\ref{semilocal}), and this result will allow us also to prove the converse of Theorem~\ref{th2generated} for noetherian domains of Krull dimension $1$ and with finitely generated integral closure in Theorem~\ref{converse}.

Let $\mathcal{C}$ be a pre-additive category, let $M$ be an object of $\mathcal{C}$, and let $I$ be a two-sided ideal of the ring $\End_\mathcal{C}(M)$. Recall that the \emph{ideal of $\mathcal{C}$ associated to $I$} is the ideal $\mathcal{A}_I$ of the category $\mathcal{C}$ defined as follows. A morphism $f\colon X\to Y$ belongs to $\mathcal{A}_I(X,Y)$ if and only if $\beta f\alpha\in I$ for every pair of morphisms $\alpha\colon M\to X$ and $\beta\colon Y\to M$ in the category $\mathcal{C}$. 

Notice that if $\mathcal{C}'$ is a full subcategory of $\mathcal{C}$. Then the restriction of $\mathcal{A}_I$ to $\mathcal{C}'$ gives an ideal of the category $\mathcal{C}'$.

\begin{Remark}\label{rem:maximal}
Let $R$ be a commutative ring, and let $\Lambda$ be a module-finite $R$-algebra. Let $M$ be a non-zero finitely generated right $\Lambda$-module with endomorphism ring $S=\End_\Lambda(M)$. Let $\varphi\colon R\to S$ denote the canonical homomorphism. Let $\fn$ be a two-sided maximal ideal of $S$, and let $\fm=\varphi^{-1}(\fn)$. By Lemma~\ref{maximal}(i), $\fm$ is a maximal ideal of $R$.

Let $X,Y$ be two objects in $\Mod\Lambda$. If $f\in\mathrm{Hom}_\Lambda (X,Y)$ is such that $f(X)\subseteq Y\fm$, then $f\in\mathcal A_\fn$. Indeed, $\beta f\alpha(M) \subseteq M\fm$ for every pair of morphisms $\alpha\colon M\to X$, $\beta\colon Y\to M$. Hence, by Lemma~\ref{maximal}(ii), $\beta f \alpha \in \fn$.

If $\mathcal C$ is a full subcategory of $\Mod\Lambda$, we will still denote by $\mathcal A_\fn$ the restriction to $\mathcal{C}$ of the ideal defined in the whole module category.

\end{Remark}

\begin{Lemma} \label{aux}
Let $R$ be a commutative ring, and let $\fm$ be a maximal ideal of $R$. Let $f\colon X\to Y$ be a homomorphism of $R$-modules such that $f_\fm(X_\fm)\subseteq Y_\fm \fm R_\fm$. Then $f(X)\subseteq Y\fm$.
\end{Lemma}
\begin{Proof}
Let $x\in X$. Then there exist $r_i\in\fm, y_i \in Y, s \in R \setminus \fm$ such that
    \[\frac{f(x)}{1}=\sum_{i=1}^n \frac{y_ir_i}{1} \frac1s =\frac ys,\qquad y\in Y\fm\]
Therefore, there exists $t\notin\fm$ such that $f(x)st=yt$. Then there exists $u\notin\fm$ such that $1-stu\in\fm$ and $f(x)=f(x)(1-stu)+f(x)stu=f(x)(1-stu)+ytu\in Y\fm$.
\end{Proof}

\begin{Prop} \label{local_to_factor}
Under the assumptions of Remark~\ref{rem:maximal}, let $N=\bigoplus_{i\in\N} M_i$ be a direct sum of finitely generated right $\Lambda$-modules which are torsion-free as $R$-modules, and consider $A,A'$ two direct summands of $N$. Assume that $\mathcal C$ contains $M,A,A'$ and every $M_i$, for each $i\in\N$. If $A_\fm\cong A'_\fm$, then $A$ is isomorphic to $A'$ in the factor category $\mathcal C/\mathcal A_\fn$.
\end{Prop}

\begin{Proof}
We have to check that there are homomorphisms $f \colon A \to A'$ and $g \colon A' \to A$ such that $\mathrm{Id}_A - gf,\mathrm{Id}_{A'} - fg \in \mathcal{A}_{\fn}$. By Remark~\ref{rem:maximal} and Lemma~\ref{aux}, it is sufficient to find $f \colon A \to A'$
and $g \colon A' \to A$ such that $(\mathrm{Id}_A - gf)_{\fm}(A_{\fm}) \subseteq A_{\fm} \fm R_{\fm}$ and  $(\mathrm{Id}_{A'} - fg)_{\fm} (A_\fm') \subseteq A'_{\fm}\fm R_{\fm}$.

Assume that $A_\fm\cong A'_\fm$. Let $\alpha\colon A_\fm\to A'_\fm$ and $\beta\colon A'_\fm\to A_\fm$ denote mutually inverse isomorphisms. Then $\gamma=\iota_{A'_\fm}\alpha\pi_{A_\fm}$ and $\delta=\iota_{A_\fm}\beta\pi_{A'_\fm}$ are elements in $\End_{R_\fm}(\bigoplus_{i\in\N}(M_i)_\fm)$ such that the following diagram commutes
    \[\begin{tikzcd}[column sep=huge,row sep=huge]
        \bigoplus_{i\in\N}(M_i)_\fm\rar[shift left]{\gamma}\dar[shift left,two heads]{\pi_{A_\fm}} 
        & \bigoplus_{i\in\N}(M_i)_\fm\lar[shift left]{\delta}\dar[shift left,two heads]{\pi_{A'_\fm}} \\
        A_\fm\rar[shift left]{\alpha}\uar[shift left,hook]{\iota_{A_\fm}}
        & A'_\fm\lar[shift left]{\beta}\uar[shift left,hook]{\iota_{A'_\fm}}
    \end{tikzcd}\]
By Lemma~\ref{isofg}, we can consider $\gamma$ as a column-finite matrix such that the $i$-th column represents $\Hom_{\Lambda_\fm}((M_i)_\fm,\bigoplus_{j\in I}(M_j)_\fm)\cong \Hom_\Lambda(M_i,\bigoplus_{j\in I} M_j)\otimes_R R_\fm$. 

For each $i\in\N$, let $s_i\notin\fm$ be the product of the denominators in the $i$-th column, and define $\tau=\bigoplus_{i\in\N} s_i\mathrm{Id}_{(M_i)_\fm}\in\mathrm{Aut}_{\Lambda_\fm}(N_\fm)$. By Lemma~\ref{isofg}, there exists $h \in \End_\Lambda(N)$ such that $\gamma\circ\tau = h_{\fm}$. For each $i\in\N$, take $t_i\notin\fm$ such that $1-t_is_i\in\fm$ and define $\theta=\bigoplus_{i\in\N} t_i\mathrm{Id}_{M_i}\in\End_\Lambda(N)$. Similarly, define $\tau'$ and $\theta'$ starting with the endomorphism $\delta$, i.e., 
$\delta \circ \tau' = h'_{\fm}$ for some $h' \in \End_\Lambda(N)$. Note that $\tau\theta_\fm=\mathrm{Id}_{N_\fm} - r$ for some $r \in {\End}_{\Lambda_{\fm}}(N_\fm)$ with $\im\ r \subseteq (N_{\fm})\fm R_{\fm}$. Similarly, there exists $r' \in {\End}_{\Lambda_{\fm}}(N_\fm)$ with $\im\ r' \subseteq (N_{\fm})\fm R_{\fm}$ such that $\tau'\theta'_\fm=\mathrm{Id}_{N_\fm} - r'$.

Set $f=\pi_{A'} h \theta\iota_A$ and $g=\pi_A h' \theta'\iota_{A'}$. Then
    \begin{align*}
        (\mathrm{Id}_A-g f)_\fm&=(\mathrm{Id}_A-\pi_A h' \theta'\iota_{A'}\pi_{A'} h \theta\iota_A)_\fm \\
        &=\mathrm{Id}_{A_\fm}-\pi_{A_\fm} \delta\tau'\theta'_\fm\iota_{A'_\fm}\pi_{A'_\fm} \gamma \tau\theta_\fm\iota_{A_\fm} \\
        &=\mathrm{Id}_{A_\fm}-\pi_{A_\fm}\iota_{A_\fm}\beta\pi_{A'_\fm}\tau'\theta'_\fm\iota_{A'_\fm}\pi_{A'_\fm}\iota_{A'_\fm}\alpha\pi_{A_\fm}\tau\theta_\fm\iota_{A_\fm} \\
        &=\mathrm{Id}_{A_\fm}-\pi_{A_\fm}\iota_{A_\fm}\beta\pi_{A'_\fm}(\mathrm{Id}_{N_\fm}-r')\iota_{A'_\fm}\pi_{A'_\fm}\iota_{A'_\fm}\alpha\pi_{A_\fm}(\mathrm{Id}_{N_\fm}-r) \iota_{A_\fm} \\
        & = \mathrm{Id}_{A_\fm}-\pi_{A_\fm}\iota_{A_\fm}\beta\pi_{A'_\fm} 1 \iota_{A'_\fm}\pi_{A'_\fm}\iota_{A'_\fm}\alpha\pi_{A_\fm} 1 \iota_{A_\fm} + s,
    \end{align*}
where $s \in \End_{\Lambda_{\fm}}(A_{\fm})$ satisfies that $\im\ s \subseteq A_{\fm} \fm R_{\fm}$.
A symmetric computation shows that $\im\ (\mathrm{Id}_{A'} - fg)_\fm \subseteq A'_{\fm} \fm R_{\fm}$.
\end{Proof}

\begin{Remark}\label{rem:finpres}
The statement in Proposition~\ref{local_to_factor} is also true if we assume that $N$ is a direct sum of finitely presented right $\Lambda$-modules instead of a direct sum of finitely generated right $\Lambda$-modules which are torsion-free as $R$-modules using \cite[Theorem~3.18]{reiner} instead of Lemma~\ref{isofg}.
\end{Remark}

\begin{Cor} \label{local_to_iso}
Let $R$ be a commutative ring, and let $\Lambda$ be a module-finite $R$-algebra. Let $M=\bigoplus_{i\in\N} M_i$ be a direct sum of non-zero finitely presented right $\Lambda$-modules with semilocal endomorphism rings. If $A$ and $A'$ are direct summands of $M$, then $A \cong A'$ if and only if $A_{\fm} \cong A'_{\fm}$ for every maximal ideal $\fm$ of $R$.
\end{Cor}
\begin{Proof}
Consider $\mathcal{C} = \Add(\bigoplus_{i \in \N} M_i)$. We aim to apply \cite[Proposition~3.6]{P4}. We have to check if $\fn$ is a maximal ideal of $\End_\Lambda(M_i)$
for some $i \in \N$, then $A$ and $A'$ are isomorphic objects of the factor category $\mathcal{C}/\mathcal{A}_{\fn}$. Note that the preimage of $\fn$ in the canonical homomorphism $R \to \End_\Lambda(N_i)$ is a maximal ideal of $R$ by Lemma~\ref{maximal}(i), so we can apply Proposition~\ref{local_to_factor} and Remark~\ref{rem:finpres}.
\end{Proof}

From now on, we only consider the case when $R$ is a noetherian domain of Krull dimension $1$. We can try to find appropriate generalizations reflecting the results in the previous parts of the paper. For the semilocal case below, we need Krull dimension $1$ because of Corollary~\ref{local}.

\begin{Prop} \label{semilocal}
Let $R$ be a semilocal noetherian domain of Krull dimension $1$. Let $\fm_0,\fm_1,\dots,\fm_k$ be the list of maximal ideals of $R$. Assume that
\begin{enumerate}
    \item[(i)] $R_{\fm_i}$ has the two-generator property for $i = 1,\dots,k$.
    \item[(ii)] The normalization of $R_{\fm_i}$ is a discrete valuation domain for every $i = 0,1,\dots,k$.
\end{enumerate}
Then every pure projective torsion-free $R$-module is a direct sum of finitely presented modules. 
\end{Prop}
\begin{Proof}
By Kaplansky's Theorem \cite[Theorem 1]{kaplansky}, every pure projective module is a direct sum of countably generated modules. Therefore, we may consider $A$ to be a direct summand of $\bigoplus_{i \in \N} M_i$, where each
$M_i$ is a non-zero finitely presented torsion-free $R$-module, say $A \oplus A' = \bigoplus_{i \in \N} M_i$. Also, we may assume that $A$ is not finitely generated.

If $\fm$ is a maximal ideal of $R$, then $A_{\fm} \oplus A'_{\fm} = \bigoplus_{i \in \N}(M_i)_{\fm}$, and by (ii) and Corollary~\ref{local},  $A_{\fm}$ is a direct sum of finitely generated $R_{\fm}$-modules with local endomorphism ring. 

For each $i = 0,\dots,k$, let $A_{\fm_i} = \bigoplus_{j \in \N} N_{i,j}$, where $N_{i,j}$ is a finitely generated indecomposable $R_{\fm_i}$-module. In \cite[Theorem~4.3]{rush}, Rush proved that over a noetherian ring with the two-generator property any finitely generated indecomposable torsion-free module has rank $1$,  thus $(i)$ implies that $N_{i,j}$ is a module of rank $1$ if $i \geq 1$. 

Therefore, the module $A$ fulfills the hypothesis of Proposition~\ref{genusdecomposition}, and since the number of maximal ideals of $R$ is finite, we can deduce that $A$ is in the same genus as a module $B$ that is a direct sum of finitely generated torsion-free modules. Since $A$ and $B$ are direct summands of $B\oplus \left( \bigoplus_{i \in \N} M_i\right)$, we deduce from Corollary~\ref{local_to_iso} that $A\cong B$. \end{Proof}

The following lemma is a variation of Proposition~\ref{prop:ic} adapted to the noetherian setting.

\begin{Lemma}\label{Pontryagin}
Let $R$ be a noetherian domain of Krull dimension $1$, and with module-finite normalization $\overline{R}$. Let $\mathcal{M}=\{\fm_1,\dots,\fm_k\}$ be the list of maximal ideals of $R$ such that $R_{\fm}$ is a principal ideal domain for any maximal ideal $\fm \notin\mathcal{M}$, and let $\Sigma = R \setminus \bigcup_{i = 1}^{k} \fm_i$. Further, let $M \subseteq \overline{R}^{(\omega)}$ be such that $M_{\Sigma}$ is a direct sum of finitely generated $R_{\Sigma}$-modules. Then $M$ is a direct sum of finitely generated modules. 
\end{Lemma} 
\begin{Proof}
We may assume that the rank of $M$ is infinite. Let $M_{\Sigma} = \bigoplus_{i = 1}^{\infty} M_i$, where each $M_i$ is a non-zero finitely generated $R_{\Sigma} $-module.
For each $i = 1,2,\dots$ let $N_i := \{ m \in M \mid \frac{m}{1} \in \bigoplus_{j = 1}^i M_j\}$. Note that $N_i \subseteq \overline{R}^{\ell}$ for some $\ell$, hence $N_i$ is a finitely generated submodule of $M$. Further, we claim that $N_i$ is an $\mathrm{RD}$-submodule of $M$. Consider $m \in M, n\in N_i$ and $0\neq r \in R$ such that $mr = n$. Then $\frac{m}{1}\frac{r}{1} = \frac{n}{1}$ in $M_{\Sigma}$. Since $\frac{r}{1}$ is not a zero-divisor on each $M_t$, $\frac{m}{1}$ is an element of $\bigoplus_{j = 1}^i M_j$. 

We claim that for every $i \in \N$ the inclusion $\iota_{i} \colon N_i \to N_{i+1}$ splits. By Lemma~\ref{locallysplit}, it is enough to check that $\iota_i \otimes_{R} R_{\fm}$ splits for every maximal ideal $\fm$ of $R$. Since $\iota_i \otimes_{R} R_{\fm}  \cong \iota_i \otimes_{R}  R_{\Sigma} \otimes_{R_{\Sigma}} R_{\fm}$ for any $\fm \in \mathcal{M}$ and $\iota_{i} \otimes_{R} R_{\Sigma}$ splits,  $\iota_{i} \otimes_{R} R_{\fm}$ splits for any $\fm \in \mathcal{M}$.

Now suppose that $\fm \notin \mathcal{M}$. Then it is easy to check that $(N_i)_{\fm}$ is an $\mathrm{RD}$-submodule of $(N_{i+1})_{\fm}$ and since $R_{\fm}$ is a valuation domain, also a pure submodule of $(N_{i})_{\fm}$. The pure exact sequence
    \[\begin{tikzcd}
        0 \rar & (N_i)_{\fm} \rar{\iota \otimes_R R_{\fm}} & (N_{i+1})_{\fm} \rar & (N_{i+1}/N_{i})_{\fm} \rar & 0
    \end{tikzcd}\]
splits, since the right-hand term is finitely presented. From this, we conclude that $\iota_i$ splits.

Overall, $M$ is a union of the chain $N_1 \subseteq N_2 \subseteq N_3 \subseteq \cdots$ where each $N_i$ splits in $N_{i+1}$, i.e., there exists a submodule $D_{i+1} \subseteq N_{i+1}$ such that $N_{i+1} = N_i \oplus D_{i+1}$. Further let $D_1:=N_1$.
Then $M = \bigoplus_{i = 1}^{\infty} D_i$ is a direct sum of finitely generated $R$-modules.
\end{Proof}

\begin{Th} \label{converse}
Let $R$ be a noetherian domain of Krull dimension $1$, and with module-finite normalization $\overline{R}$.  Then the following statements are equivalent:
\begin{itemize}
    \item[(1)] Every pure projective torsion-free $R$-module is a direct sum of finitely presented modules. 
    \item[(2)] For any finitely generated, torsion-free $R$-module $X$, every element in $\Add (X)$ is a direct sum of finitely generated modules.
    \item[(3)] $R$ satisfies the following two conditions: 
    \begin{enumerate}
        \item[(i)] there exists at most one maximal ideal $\fm _0$ of $R$ such that $R_{\fm _0}$ is not a Bass domain, and
        \item[(ii)] the normalization of $R_{\fm}$ is a discrete valuation domain for every maximal ideal $\fm$ of $R$.
    \end{enumerate}
\end{itemize} 
\end{Th}
\begin{Proof} 
It is clear that $(1)$ implies $(2)$. Statement $(2)$ implies $(3)$ by Theorem~\ref{th2generated}. We only need to prove that $(3)$ implies $(1)$.

Note that, if $R$ is integrally closed, then $R$ is a Dedekind domain, and pure projective torsion-free modules are projective. Assume that $R \neq \overline{R}$, and let $\mathcal{M}=\{\fm_1,\dotsc,\fm_k\}$ be a finite set of maximal ideals of $R$ such that $R_{\fm}$ is integrally closed whenever $\fm \not \in \mathcal{M}$. Let $\Sigma := R \setminus \bigcup_{i=1}^k \fm_i$. Therefore, $R_{\Sigma}$ is a semilocal ring satisfying the assumptions of Proposition~\ref{semilocal}. Hence, any pure projective torsion-free $R_{\Sigma}$-module is a direct sum of finitely presented modules. 

Let $A$ be a countably generated pure projective torsion-free $R$-module, i.e., a direct summand of $\bigoplus_{i \in \N} M_i$, where each $M_i$ is a finitely generated torsion-free module. Note that each $M_i$ can be considered as a submodule of $M_i\overline{R}$ (the $\overline{R}$-submodule of $(M_i)_{0}$ generated by $M_i$) which is a projective $\overline{R}$-module. Hence, we may consider $A$ as a submodule of $\overline{R}^{(\omega)}$. 

Since $A_{\Sigma}$ is a pure projective torsion-free $R_{\Sigma}$-module, it has to be a direct sum of finitely generated modules. By Lemma~\ref{Pontryagin}, $A$ is a direct sum of finitely generated modules. 
\end{Proof}

\section{A family of  examples} \label{s:example}

In this section, we provide a family of examples of $h$-local domains of Krull dimension $1$, not necessarily noetherian, that exemplifies well the situations described in Theorem~\ref{th2generated} and also in Corollary~\ref{integrallyclosed} as we show that they satisfy that direct summands of finitely generated torsion-free modules are direct sums of finitely generated modules. This family of examples was suggested to us by Carmelo Antonio Finocchiaro and Paolo Zanardo \cite{FZ}.

Let $K\subseteq L$ be a  field extension, and consider the ring $R=K+xL[x]$ and its localization $R_\fm$ at the maximal ideal $\fm=xL[x]$. The field of fractions of $R$ is $L(x)$, and  $R$ fits in the pullback diagram of rings
    \[\begin{tikzcd}[column sep=huge,row sep=huge]
        R=K+xL[x]\rar[hook]{\text{incl.}}\dar[swap,two heads]{\mathrm{ev}_0} & L[x]\dar[two heads]{\mathrm{ev}_0} \\
        K\rar[hook]{\text{incl.}} & L
    \end{tikzcd}\qquad (\mathrm{C1})\]
where $\mathrm{ev}_0$ denotes the evaluation at 0. In fact, this diagram is a conductor square, with conductor ideal $\mathfrak{c}=\fm$ since $R/\fm\cong K$ and $L[x]/\fm\cong L$. 

Note that $R_\fm=K+xL[x]_\fm$, since every element in $R_\fm$ can be written as
    \[\frac{a+xp(x)}{1+xq(x)}=a+\frac{x(p(x)-aq(x))}{1+xq(x)},\]
with $a\in K$, $p(x),q(x)\in L[x]$.

The localization of $(*)$ at $\fm$ gives  also a conductor square for $R_\fm$
\[\begin{tikzcd}[column sep=huge,row sep=huge]
        R_\fm =K+xL[x]_\fm\rar[hook]{\text{incl.}}\dar[swap,two heads]{\mathrm{ev}_0} & L[x]_\fm\dar[two heads]{\mathrm{ev}_0} \\
        K\rar[hook]{\text{incl.}} & L
    \end{tikzcd}\qquad (\mathrm{C2})\]
with conductor $\mathfrak{c}_\fm =xL[x]_\fm$.

\begin{Lemma} \label{KLbasic}
Let $K\subseteq L$ be a  field extension, let $R$ denote the ring $K+xL[x]$, and $\fm=xL[x]$. Then
\begin{enumerate}
    \item[(i)] If $I$ is an ideal of $R$, then $IL=L[x]$ if and only if $I=R$.
    \item[(ii)] Every maximal ideal $\fn\neq\fm$ of $R$ is generated by a polynomial $q(x)\in R$ that is irreducible in $L[x]$ with $q(0)=1$.
    \item[(iii)] $R_\fn=L[x]_\fn$, that is, $R_\fn$ is a discrete valuation ring  for every maximal ideal $\fn\neq\fm$ in $R$.
    \item[(iv)] $R$ is an $h$-local domain of Krull dimension $1$.
\end{enumerate}
\end{Lemma}
\begin{Proof}
$(i)$. Suppose that $IL=L[x]$, so $y_0+y_1a_1+\dotsb+y_na_n=1$ for some elements $y_i\in I$ and $a_i\in L\setminus K$. Then, for every  polynomial $q(x)\in L[x]$, $y_0xq(x)+y_1a_1xq(x)+\dotsb+y_na_nxq(x)=xq(x)\in I$. Therefore $xL[x]\subseteq I$. Since $xL[x]$ is a maximal ideal of $R$ and $I\not \subseteq xL[x]$, $I=K+xL[x]$.

$(ii)$. Let $\fn \neq\fm$ be a maximal ideal of $R$. By $(i)$ and because $\fn$ is maximal, $\fn L\neq L[x]$ and $\fn L\cap R=\fn$. There exists  $p(x)\in L[x]$ with $p(0)=1$ such that $p(x)L[x]=\fn L$. Since $p(x) \in \fn L\cap R$, $p(x)\in \fn$. Notice that if $q(x)\in L[x]$ and $p(x)q(x)\in R$ then also $q(x)\in R$. Therefore $p(x)R=\fn$. If $p(x)=q_1(x)q_2(x)$ in  $L[x]$, and since $p(0)=1$, we can assume that $q_1(0)=q_2(0)=1$ so that it is also a decomposition in $R$. Since $\fn$ is maximal, either $q_1(x)$ or $q_2(x)\in \fn$ which implies that either $q_1(x)$ or $q_2(x)$ is equal $1$.   This shows that $p(x)$ is irreducible in $L[x]$.

$(iii)$. Since the element $x\in R\setminus\fn$ is invertible in $R_\fn$, we have $a=ax/x\in R_\fn$ for every element $a\in L$, so we deduce that $L\subseteq R_\fn$, that is, $L[x]_\fn\subseteq R_\fn$.

$(iv)$.  Let $p(x)$ be a non-zero polynomial in $R$. If $p(0)\neq 0$, $p(x)$ can be written as $p(x)=kq_1(x)\dotsb q_n(x)$, where $k=p(0)\in K^*$ and $q_i(x)$ are irreducible polynomials in $L[x]$ with $q_i(0)=1$ for every $i=1,\dots,n$ so, in particular, $q_i(x)\in R$. This decomposition is unique up to a unit (an element in $K$) because $L[x]$ is a UFD, and in this case, $p(x)$ is only contained in the maximal ideals generated by the irreducible polynomials $q_i(x)$ for $i=1,\dots,n$.

If $p(0)=0$, then $p(x)$ can be written as $p(x)=x^mq(x)$, where $m\ge1$ and $q(x)\in L[x]$ with $q(0)\neq 0$. Then $q(x)$ can be written as $q(x)=lq_1(x)\dotsb q_n(x)$, where $l=q(0)\in L^*$ and $q_i(x)$ are irreducible polynomials in $L[x]$ with $q_i(0)=1$ for every $i=1,\dots,n$, so again, $q_i(x)\in R$. Therefore, $p(x)$ is only contained in the maximal ideals generated by the irreducible polynomials $q_i(x)$ for $i=1,\dots,n$ and it is also contained in $\fm$. Therefore, $R$ has finite character.

Now we prove that $R_\fm$ has Krull dimension $1$. Let $\mathfrak{p}$ be a non-zero prime ideal of $R_\fm$ and let $0\neq g\in\mathfrak{p}$. We can assume, up to a unit of $R_\fm$, that $g$ is of the form $x^na$ with $a\in L\setminus\{0\}$ and $n\ge 1$. If $n>1$, then $x^na=x\cdot (x^{n-1}a)$ which implies that either $x$ or $x^{n-1}a$ are  in $\mathfrak{p}$. By induction on $n$ we can deduce that $\mathfrak{p}$ contains an element of the form $xa$ with $a\in L\setminus \{0\}$ so it also contains $x^2b=xa\cdot xa^{-1}b$ for any $b\in L$. As $(xb)^2\in \mathfrak{p}$ for any $b\in L$, we deduce that $xb\in \mathfrak{p}$. Hence $\mathfrak{p}=\fm $.

Since the Krull dimension of $R$ localized at any maximal ideal is $1$, we deduce that $R$ has Krull dimension $1$ and, being of finite character, it is $h$-local.
\end{Proof}

From (iii) it follows that finitely generated torsion-free $R_\fn$-modules are free for every maximal ideal $\fn\neq\fm$ in $R$.

\begin{Remark} 
Assume $K\varsubsetneq L$.
\begin{enumerate}
    \item \textit{$R$ is integrally closed if and only if the extension $K\subseteq L$ is purely transcendental.} If $\alpha \in L\setminus K$ is algebraic over $K$ then $\alpha$ is integral over $R$. So if $R$ is integrally closed, $K\subseteq L$ is a purely transcendental extension.
    
    To prove the converse, let $f\in L(x)$ satisfy a polynomial equation $f^n+a_{n-1}f^{n-1}+\dotsb+a_0=0$ with  $a_i \in R$ for $i= 0,\dots, n-1$. Since $L[x]$ is already integrally closed, we can assume that $f\in L[x]$. Let $f_0$ denote the evaluation at 0 of $f$. Then $f_0$ satisfies a polynomial equation in $K$, so $f_0$ is algebraic over $K$. Since $L$ is purely transcendental, we deduce that $f_0\in K$, so $f\in R$.
    
    \item \textit{$R_\fm$ is  noetherian if and only if the extension $K\subseteq L$ has finite degree.} Indeed, let $(a_i)_{i\in A}$ be a basis of $L$ as $K$-vector space.  The maximal ideal $\fm=xL[x]$ is generated by $\mathcal{A}=\{xa_i\}_{i\in A}$ and no proper subset of $\mathcal{A}$ generates $\fm$. Then $R_\fm$ is noetherian if and only if $A$ is finite.
    
    \item \textit{$R_\fm$ is not a valuation ring.} Note that every element $a\in L\setminus K$ is in its field of fractions $L(x)$, but neither $a$ nor $a^{-1}$ belong to $R_\fm$.
\end{enumerate}
\end{Remark}

Now we will study torsion-free modules over the ring $R$ and over the ring $R_ \fm$. The inclusion of $R$ into the principal ideal domain $L[x]$ allows us to use the techniques of the \emph{conductor square and artinian pairs} (see \cite{leuschke}). First, we fix some notation.

We shall describe a class of torsion-free modules over the ring $T$ that can be either $R$ of $R_ \fm$. The ring $T$ is included in a principal ideal domain $S$, where  $S=L[x]$ if $T=R$ and $S=L[x]_\fm$ if $T=R_\fm$. The field of fractions of $T$  coincides with the one of $S$ and it is $Q=L(x)$. We will denote by $\mathfrak{c}$ the conductor of  both conductor squares (C1) and (C2).

We denote by $\lambda \colon T\to Q$ and $\lambda '\colon S\to Q$ the corresponding localization maps, and by $\varepsilon \colon T\to S$ the ring inclusion. Then if $M_T$ is a torsion-free $T$-module, there is a commutative diagram
    \[\begin{tikzcd}[row sep=huge,column sep=huge]
        M \rar[hook]{M\otimes \lambda}  \drar[hook][swap]{M\otimes \varepsilon} & M\otimes _TQ  \\
        & M\otimes _TS \uar[swap]{M\otimes \lambda '}
    \end{tikzcd}\]
Therefore, we can identify $M_T$ with an essential submodule of $Q^{(A)}$ for a suitable set $A$, and then, we set $(M\otimes \lambda ')(M\otimes _TS)=MS$.

In general, $M\otimes \lambda '$ is not an injective map; it is when $M_T$ is projective or, more generally, when $M_T$ is flat.

\begin{Lemma}\label{basicmoduleskl} Let $M_T$ be a torsion-free module over $T$. Then:
\begin{itemize}
\item[(i)] $MS=ML$;
\item[(ii)] $ML\cap M\supseteq M\mathfrak{c}$;
\item[(iii)] As $K$-vector space, $M=V\oplus M\mathfrak{c}$, where $V$ is any complement of $M\mathfrak{c}$ in $M$. Moreover, as $L$-vector space, $ML=VL\oplus M\mathfrak{c}$. Setting $W=VL$, it follows that $\mathrm{dim} _K(V)\ge \mathrm{dim} _L(W)$.
\item[(iv)] As $L$-vector space, $ML=W\oplus M\mathfrak{c}$, where $W$ is any complement of $M\mathfrak{c}$ in $M$. Moreover, $V=W\cap M$ is a $K$-vector space such that $M=V\oplus M\mathfrak{c}$ and $W=VL$.
\item[(v)] For any pair $V$, $W$ chosen as in $(iii)$ or $(iv)$, $V\cong M/M\mathfrak{c}$ and $W\cong ML/M\mathfrak{c}$. In addition,  there is a pullback diagram 
\[\begin{tikzcd}[row sep=huge,column sep=huge]
        M\rar[hook]{\text{incl.}}\dar[two heads] & ML\dar[two heads] \\
        V\rar[hook]{\text{incl.}} & W
    \end{tikzcd} \qquad(*)\]
\end{itemize}
\end{Lemma}
\begin{Proof} 
Statement $(i)$ follows because $TS=TL$, cf. Lemma~\ref{KLbasic}(i). Statement $(ii)$ is clear.

The ring $T$ fits in an exact sequence of $T$-modules
    \[\begin{tikzcd}
        0\rar & \mathfrak{c}\rar & T\rar & K\rar & 0
    \end{tikzcd}\] 
    
shows that $(*)$ is a pull-back diagram.
\end{Proof}

\begin{Def} \label{modartinianpair} The category $\mathcal{B}$ of modules over the artinian pair $K\hookrightarrow L$ has as objects the inclusions $V\hookrightarrow W$ where $V$ is a $K$-subspace of the $L$-vector space $W$ satisfying that $VL=W$. 

If $V\hookrightarrow W$ and $V'\hookrightarrow W'$ are two objects of $\mathcal{B}$ a morphism between them consists of one $K$-linear map $f\colon V\to V'$ and one $L$-linear map $g\colon W\to W'$ making the diagram
\[\begin{tikzcd}[row sep=huge,column sep=huge]
        V\rar[hook]{\text{incl.}}\dar{f} & W\dar{g} \\
        V'\rar[hook]{\text{incl.}} & W'
    \end{tikzcd} \]
commutative.    
\end{Def}

\begin{Cor} \label{functor} As usual, $T$ denotes the ring $R=K+xL[x]$ or $R_\fm$ where $\fm =xL[x]$. Let $\mathcal{A}_T$ be the category of torsion-free modules over $T$. Let $\mathcal{B}$ be the category of modules over the artinian pair $K\hookrightarrow L$. Then there is a functor $F_T$ that assigns to each object $M_T$ in $\mathcal{A}_T$ the object of $\mathcal{B}$, $M/M\mathfrak{c}\hookrightarrow ML/M\mathfrak{c}$.

Moreover, $F_R (M)=F_{R_\fm}(M_\fm)$ for any object $M$ of $\mathcal{A}_R$.
\end{Cor}
\begin{Proof}
Lemma~\ref{basicmoduleskl}(v) and the fact that a morphism between torsion-free modules over $T$ induces a morphism between the modules over the artinian pair, imply the existence of such a functor.

If $M$ is an object in  $\mathcal{A}_R$ then $M/M\mathfrak{c}$ is an $R$-module and also an $R_{\fm}$-module. Also, $ML/M\mathfrak{c}$ is an $L[X]$ module as well as an $L[x]_\fm$-module. So, the second part of the statement is clear.
\end{Proof}

For a general torsion-free module $M$ the pull-back diagram of Lemma~\ref{basicmoduleskl}(v) can be trivial. Take, for example, $M=Q$. As  $M=ML=M\mathfrak{c}$ and $V=W=0$. So, in the notation of Corollary~\ref{functor}, $F_T(Q)=0$.

To get a better correspondence, we need to restrict the class of torsion-free modules we are interested in. We  shall consider the class
\[\mathcal{C}_T=\{ M\mid M \mbox{ is a torsion-free $T$-module such that $ML$ is a free $S$-module}\}\]

Notice that $\mathcal{C}_T$ contains all finitely generated torsion-free modules, and it is closed by arbitrary direct sums and direct summands. In addition,  $M=S\in \mathcal{C}_T$.  

Assume that $ML$ is a free $S$-module, then we can fix an $S$-basis $\mathcal{B}=\{v_i\}_{i\in A}$. Then if $W$ is the $L$-vector space generated by $\mathcal{B}$, $M\mathfrak{c}=W\mathfrak{c}$ and $ML=W\oplus W\mathfrak{c}$. Therefore, by the modular law  and by Lemma~\ref{basicmoduleskl}, $M=V\oplus W\mathfrak{c}$ where $V=W\cap M$.

If $M_1$ and $M_2$ are two modules in $\mathcal{C}_T$ then, for $i=1,2$, $M_iL=W_i\oplus W_i\mathfrak{c}$  and  $M_i=V_i\oplus W_i\mathfrak{c}$ where $V_i=W_i\cap M_i$. Therefore, if $(f,g)$ is a morphism between the artinian pairs $V_1\hookrightarrow W_1$ and $V_2\hookrightarrow W_2$ or, equivalently, $f$ is a $K$-linear map and $g$ is an $L$-linear map, and there is a commutative diagram
    \[\begin{tikzcd}[column sep=huge,row sep=huge]
        V_1\rar[hook]{\text{incl.}}\dar{f} & W_1\dar{g} \\
        V_2\rar[hook]{\text{incl.}} & W_2
    \end{tikzcd} \qquad(*)\]
then $f$ can be extended to an $S$-linear map $$\widetilde{g}\colon M_1L=W_1\oplus W_1\mathfrak{c}\to M_2L=W_2\oplus W_1\mathfrak{c}$$ by setting $\widetilde{g}(wx^n)=g(w)x^n$ for any $w\in W_1$ and any $n\ge 0$. Notice that $g$ is an isomorphism if and only if $\widetilde{g}$ is an isomorphism because $\widetilde{g^{-1}}=(\widetilde{g})^{-1}$.

As a consequence, we have

\begin{Cor} \label{functorc}
The functor $F_T$ described in Corollary~\ref{functor} is full when restricted to the category $\mathcal{C}_T$ and it reflects isomorphisms.

Therefore,
\begin{itemize}
    \item[(i)] two objects $M_1$ and $M_2$  of $\mathcal{C}_R$ are   isomorphic if and only if $(M_1)_\fm$ and $(M_2)_\fm$ are isomorphic;
    \item[(ii)] $M_T$ is projective if and only if $M\in \mathcal{C}_T$ and the inclusion $ M/M\mathfrak{c} \hookrightarrow ML/M\mathfrak{c}$ sends $K$-basis of $M/M\mathfrak{c}$ to $L$-basis of $ML/M\mathfrak{c}$ if and only if $M_T$ is a free $T$-module.
 \end{itemize}
\end{Cor}
\begin{Proof}
The remarks before the statement prove that the functor $F_T$ restricted to the category $\mathcal{C}_T$ is full and reflects isomorphisms.

The statement $(i)$ follows because, by Corollary~\ref{functor}, $F_R (M)=F_{R_\fm}(M_\fm)$ for any object $M$ of $\mathcal{C}_R$ and, by the first part of the statement, two modules in $\mathcal{C}_R$ are isomorphic if and only if their corresponding artinian pairs are isomorphic if and only if their localizations at $\fm$ are isomorphic. 

To prove $(ii)$, notice that since over $R_\fm$ all projective modules are free, it follows by $(i)$ that all projective modules over $R$, since they are modules in $\mathcal{C}_R$, are isomorphic to a free module.  It is easy to check that this happens if and only if a $K$-basis of $V\cong M/M\mathfrak{c}$ is also an $L$-basis of $W=ML/M\mathfrak{c}$.
\end{Proof}

Now we summarize some results on finitely generated torsion-free modules over $R=K+xL[x]$. As observed before, these are always in $\mathcal{C}_T$ as well as infinite direct sums of them.

\begin{Lemma} \label{isoindRm}
Let $K\subseteq L$ be a  field extension, let $R$ denote the ring $K+xL[x]$, and $\fm=xL[x]$. Let $M$, $N$ and $\{M_i\}_{i\in I}$ be finitely generated, torsion-free $R$-modules. Then
\begin{enumerate}
    \item[(i)] Every finitely generated indecomposable torsion-free $R_\fm$-module has local endomorphism ring.
    \item[(ii)]  $M\cong N$ if and only if $M_\fm\cong N_\fm$.
    \item[(iii)] $M$ is indecomposable if and only if $M_\fm$ is indecomposable.
    \item[(iv)] If $M$ and $\{M_i\}_{i\in I}$ are indecomposable and $M$ is a direct summand of $\bigoplus_{i\in I} M_i$, then $M\cong M_i$ for some $i\in I$.
\end{enumerate}
\end{Lemma}
\begin{Proof} $(i)$ This follows from Lemma~\ref{tecnical_extension}, taking $T=L[x]_\fm$. Note that $\fm$ is the conductor, so it is different from zero as required in the hypothesis of the Lemma.

Statement $(ii)$ is included in Corollary~\ref{functorc}.

$(iii)$. Suppose that $M_\fm=A\oplus B$, and $M$ is indecomposable. Let $A',B'$ be $R$-submodules of $A$ and $B$ generated by some finite set of $R_\fm$-generators of $A$ and $B$, respectively. Then $(A'\oplus B')_\fm=A\oplus B=M_\fm$ and by (ii), we deduce that $A'\oplus B'\cong M$. Therefore, $A'=0$ or $B'=0$, that is, $A=0$ or $B=0$, and $M_\fm$ is indecomposable. The other implication is clear.

$(iv)$. Suppose that $M$ is a direct summand of $\bigoplus_{i\in I} M_i$. Then $M_\fm$ is also a direct summand of $(\bigoplus_{i\in I} M_i)_\fm\cong\bigoplus_{i\in I}(M_i)_\fm$. Since $(M_i)_\fm$ has local endomorphism ring, it satisfies the Krull-Schmidt property, that is, $M_\fm\cong(M_i)_\fm$ for some $i\in I$. By (ii), we deduce that $M\cong M_i$ for some $i\in I$.
\end{Proof}

\begin{Cor}
Let $K\subseteq L$ be a  field extension, let $R$ denote the ring $K+xL[x]$. Then the class of modules that are direct sums of finitely generated torsion-free $R$-modules is closed under direct summands. 
\end{Cor}
\begin{Proof}
By Lemmas~\ref{KLbasic} and \ref{isoindRm}, the result follows from Proposition~\ref{prop:ic}.
\end{Proof}

Now we want  to explicitly construct finitely generated, indecomposable, torsion-free $T$-modules where $T$ denotes either $R=K+xL[x]$ or $R_\fm$. First, we will specify better  how these modules and their endomorphism rings can look like.

Recall that if $T=R=K+xL[x]$ then $S=L[x]$, and if $T=R_\fm$ then $S=L[x]_\fm$.

\begin{Remark} \label{Kendomorphism}
Let $M_T$ be a finitely generated torsion-free $T$-module of rank $n$. As $ML$ is a finitely generated free $S$-module,   we may assume that $M$ is a $T$-submodule of $S^n$ such that $ML = S^n$, in 
particular $(xS)^n \subseteq M$. Therefore, $M_T = V + (xS)^n$, where $V$ is a $K$-subspace of $L^n$ satisfying $VL = L^n$. Notice that, this gives us a very explicit construction of a $T$-module $M\in \mathcal{C}$ such that   $F(M)$ (cf. Corollary~\ref{functor})  has as an image the module over the artinian pair $V\hookrightarrow W$.

Now we can also identify the endomorphism ring of $M_T$ with a subring of $M_n(Q)$,  in fact $\End_T(M_T)=\{A\in M_n(Q)\mid AM\subseteq M\}.$ 
Since $(xS)^n \subseteq M$, if $A \in M_n(Q)$ represents 
an endomorphism of $M_T$, then $A \in M_n(S)$, i.e., $\End_T(M_T)=\{A\in M_n(S)\mid AM\subseteq M\}.$ 
Every matrix $A\in M_n(S)$ can be uniquely decomposed as the sum of a matrix $B\in M_n(L)$ and a matrix $C\in M_n(xS)$. Since $M_n(xS) \subseteq \End_T(M)$, $A \in \End_T(M)$ if and only if 
$BV \subseteq V$. This is to say that if $A\in \End_T(M_T)$, then $F(A)=B$.
\end{Remark}

\begin{Lemma} 
Let $K\subsetneq L$ be a field extension, let $R$ denote the ring $K+xL[x]$, and $\fm=xL[x]$. Let $M$ be a finitely generated torsion-free right $R_\fm$-module of rank $n$. Then $M_n(xL[x]_\fm)\subseteq J(\End_{R_\fm}(M))$.
\end{Lemma}
\begin{Proof}
Note that $J(R_\fm)=\fm R_\fm$. Then, for every $A\in M_n(xL[x]_\fm)$, we have $M=AM+(I-A)M\subseteq J(R_\fm)M+(I-A)M$. 
By Nakayama's Lemma, $M=(I-A)M$ and, by Lemma~\ref{basicfg}(i), $I-A$ is bijective. Hence $A\in J(\End_{R_\fm}(M))$.
\end{Proof}

The following construction is a modification of \cite[Construction 3.13]{leuschke} to build indecomposable modules over artinian pairs. Using Remark~\ref{Kendomorphism}, this immediately yields the existence of   indecomposable finitely generated $T$-modules of arbitrary finite rank $n\ge 2$.

\begin{Construction}
Let $n\ge 2$ be a fixed positive integer, and suppose we have chosen $\alpha,\beta\in L$ with $\{1,\alpha,\beta,\alpha^2,\alpha\beta,\beta^2\}$ linearly independent over $K$. Let $I$ be the identity $n\times n$ matrix and $H$ be the nilpotent $n\times n$ matrix with $1$ below the diagonal and $0$ elsewhere. For $t\in K$, we consider the $n\times 2n$ matrix,
    \[\Psi_t:=[I\mid\alpha I+\beta (tI+H)]=\left[\begin{array}{cccc|cccc}
        1 & 0 & \cdots & 0 & \alpha+t\beta & 0 & \cdots & 0 \\
        0 & 1 & \cdots & 0 & \beta & \alpha+t\beta & \cdots & 0 \\
        \vdots & \vdots &\ddots &\vdots & \vdots & \vdots &\ddots &\vdots \\
        0 & 0 & \cdots & 1 & 0 & 0 & \cdots & \alpha+t\beta
    \end{array}\right]\]
Let $V_t$ be the $K$-subspace of $L^n$ spanned by the columns of $\Psi_t$. Let $A\in\Hom_K(V_t,V_u)$, where
    \[\Hom_K(V_t,V_u)=\{A\in M_n(L)\mid AV_t\subseteq V_u\}.\]
The condition $AV_t\subseteq V_u$ implies that there is a $2n\times 2n$ matrix $\theta\in M_{2n}(K)$ such that $A\Psi_t=\Psi_u\theta$. Write $\theta=\begin{bmatrix} C & D \\ P & Q \end{bmatrix}$, where $C,D,P,Q\in M_n(K)$. Then, using the condition $A\Psi_t=\Psi_u\theta$, we have the following two equations: 
    \[\begin{cases} A=C+\alpha P+\beta(uI+H)P \\\alpha A+\beta A(tI+H)=D+\alpha Q+\beta(uI+H)Q\end{cases}\]
Substituting the first equation into the second and combining terms, we get the following:
    \begin{align*}
        -D&+\alpha(C-Q)+\beta(tC-uQ+CH-HQ)+\alpha^2P \\
        &+\alpha\beta(tP+uP+HP+PH)+\beta^2(tuP+HPH+tHP+uPH)=0.
    \end{align*}
From the linear independence of $\{1,\alpha,\beta,\alpha^2,\alpha\beta,\beta^2\}$, we have
    \[D=P=0,\qquad A=C=Q,\qquad (t-u)A+AH=HA.\]
In particular, we deduce that $A\in M_n(K)$, and if $A$ is an isomorphism and $t\neq u$, then the third equation above gives a contradiction since the left side is invertible and the right side is not. Thus, if $V_t\cong V_u$, then $t=u$. To see that $V_t$ is indecomposable, we take $u=t$ and suppose that $A$ is idempotent. But $AH=HA$, and it follows that $A$ is in $K[H]$, which is a local ring. Therefore $A=0$ or $I$, as desired. 

By Remark~\ref{Kendomorphism}, $M_t=V_t+(xS)^n$, where $V_t$ is the $K$-subspace of $L^n$ constructed above, is an indecomposable $T$-module of rank $n$ (as usual, $T=R=K+xL[x]$ or $T=R_\fm$). Notice that, in view of Corollary~\ref{functor}, there are, at least, $|K|$ different isomorphism classes of such modules.
\end{Construction}

\bibliographystyle{amsplain}
\bibliography{references}

\end{document}